\documentclass[11pt]{article}
\usepackage[margin=1in]{geometry}

\usepackage{amsmath, amssymb, amsthm, mathtools}
\usepackage{graphicx}
\usepackage{placeins}
\usepackage{hyperref}
\usepackage{bookmark}

\newtheorem{theorem}{Theorem}
\newtheorem{lemma}[theorem]{Lemma}
\newtheorem{proposition}[theorem]{Proposition}
\newtheorem{corollary}[theorem]{Corollary}
\theoremstyle{definition}
\newtheorem{definition}[theorem]{Definition}
\theoremstyle{remark}
\newtheorem{remark}[theorem]{Remark}

\newcommand{\bracketN}{4000}
\newcommand{\bracketEps}{0.15}
\newcommand{\bracketA}{0.032}

\newcommand{\expEigWall}{1.00}
\newcommand{\errWallMaxPct}{3}
\newcommand{\seedsRates}{5}

\newcommand{\expTau}{1.00}
\newcommand{\expCert}{0.56}
\newcommand{\expProjScaling}{1.00}
\newcommand{\expProjSphere}{1.00}
\newcommand{\aCrossLo}{0.064}
\newcommand{\aCrossHi}{0.128}
\newcommand{\lFocus}{3}
\newcommand{\projN}{4000}
\newcommand{\projManifold}{S^1}
\newcommand{\projEps}{0.15}
\newcommand{\expCertIn}{0.52}

\newcommand{\refineA}{0.0208}
\newcommand{\refineAFrac}{0.5}
\newcommand{\refineSeeds}{10}
\newcommand{\refineL}{3}
\newcommand{\refineNLo}{522}
\newcommand{\refineNHi}{5236}
\newcommand{\refineEpsHi}{0.240}
\newcommand{\refineEpsLo}{0.024}
\newcommand{\refineDeltaHi}{0.00250}
\newcommand{\refineDeltaLo}{0.000250}
\newcommand{\refineDegree}{40}
\newcommand{\refineWallExp}{-0.05}
\newcommand{\refineScalingExp}{-1.03}
\newcommand{\refineWorstPctLo}{14}
\newcommand{\refineWorstPctHi}{16}
\newcommand{\refineCertExp}{-0.05}
\newcommand{\refineCertLo}{0.54}
\newcommand{\refineCertHi}{0.60}
\newcommand{\refineProjGapLo}{79}
\newcommand{\refineProjGapHi}{311}

\newcommand{\sTwoA}{0.0125}
\newcommand{\sTwoAFrac}{0.5}
\newcommand{\sTwoSeeds}{10}
\newcommand{\sTwoL}{4}
\newcommand{\sTwoNLo}{1389}
\newcommand{\sTwoNHi}{8000}
\newcommand{\sTwoEpsHi}{0.24}
\newcommand{\sTwoEpsLo}{0.10}
\newcommand{\sTwoDeltaHi}{0.0015}
\newcommand{\sTwoDeltaLo}{0.000625}
\newcommand{\sTwoDegree}{20}
\newcommand{\sTwoWallExp}{-0.03}
\newcommand{\sTwoScalingExp}{-0.52}
\newcommand{\sTwoWorstPctLo}{18}
\newcommand{\sTwoWorstPctHi}{20}
\newcommand{\sTwoCertExp}{-0.08}
\newcommand{\sTwoCertLo}{0.60}
\newcommand{\sTwoCertHi}{0.69}
\newcommand{\sTwoProjGapLo}{126}
\newcommand{\sTwoProjGapHi}{202}

\usepackage{authblk}
\newcommand{\ukdate}{\number\day\ifcase\day\or st\or nd\or rd\or th\or th\or
  th\or th\or th\or th\or th\or th\or th\or th\or th\or th\or th\or th\or th\or
  th\or th\or st\or nd\or rd\or th\or th\or th\or th\or th\or th\or th\or st\fi
  \space\ifcase\month\or January\or February\or March\or April\or May\or June\or
  July\or August\or September\or October\or November\or December\fi
  \space\number\year}
\title{Spectral Stability and Convergence of the Unnormalized Graph
    Laplacian on Noisy Point Clouds}
\author{Daniel Ageev}
\author{Yariv Aizenbud}
\affil{School of Mathematical Sciences, Tel Aviv University, Israel}
\date{\ukdate}
\hypersetup{
  pdftitle={Spectral Stability and Convergence of the Unnormalized Graph
    Laplacian on Noisy Point Clouds},
  pdfauthor={Daniel Ageev, Yariv Aizenbud}}

\begin{document}
\maketitle

\begin{abstract}
The first eigenvalues and eigenvectors of the unnormalized graph Laplacian of
clean samples from a manifold, suitably scaled, approximate those of a weighted
Laplace--Beltrami operator. However, the samples are usually observed with
noise and therefore need not lie on the manifold. We build the clean and noisy
graph Laplacians, with bandwidth $\varepsilon$, on the clean samples and on the
noisy observations, where the noise is arbitrary but bounded by $\delta$. We
show that, with high probability, the first eigenvalues and eigenvectors of the
noisy graph Laplacian are close to those of the clean one: the eigenvalues
within $O(\delta/\varepsilon)$, and the subspaces spanned by the first
eigenvectors within $O(\sqrt{\delta/\varepsilon})$, whether or not the
eigenvalues are simple.
Consequently, if $\delta$ is small enough relative to $\varepsilon$, the noisy
graph Laplacian approximates the weighted Laplace--Beltrami operator at the same
rate as the clean one. We also show that, for large noise, we can get the same
rate as the clean graph Laplacian by first using a denoiser, and then building
the graph Laplacian on the denoised observations.
\end{abstract}

\section{Introduction}\label{sec:intro}

A point cloud is often analysed through a weighted graph on its points, in
which the weight of the edge between two points is given by a kernel, a
function of their distance and a bandwidth $\varepsilon$. The unnormalized
graph Laplacian of this graph is the diagonal matrix of degrees minus the
matrix of weights, where the degree of a point is the sum of the weights of
its edges (see, e.g., \cite{von_luxburg}). The first (smallest) eigenvalues of this
matrix, or of a normalized version of it, and their eigenvectors are used
to cluster the points \cite{von_luxburg}, to describe the geometry of the
point cloud with diffusion maps \cite{coifman_lafon}, and to embed the
points in a low-dimensional space with Laplacian eigenmaps
\cite{belkin_niyogi}.

When the points are i.i.d.\ samples from a density on a manifold, the
unnormalized graph Laplacian approximates a weighted Laplace--Beltrami
operator \cite{calder_gt,gt_gerlach_hein_slepcev}, a continuum operator
determined by the manifold and the density (see, e.g.,
\cite{gt_li_venkatraman}). This approximation has been studied through pointwise limits
\cite{belkin_niyogi_2008,gine_koltchinskii,hein_audibert_von_luxburg,singer},
spectral convergence on domains in
$\mathbb R^d$ \cite{gt_slepcev}, and quantitative spectral estimates for
discretizations of a manifold and samples from it
\cite{burago_2014,gt_gerlach_hein_slepcev}. Calder and
Garc\'ia Trillos \cite{calder_gt} prove that, with high
probability, the first eigenvalues of the graph Laplacian approximate those of
the weighted Laplace--Beltrami operator at rate $O(\varepsilon)$ if
$\varepsilon$ tends to zero slowly enough as the number of points grows. They
also prove, under the same conditions, that the corresponding unit
eigenvectors approximate eigenfunctions of this operator at the same rate,
even for a multiple eigenvalue. For a simple
eigenvalue and slightly different normalizations, Garc\'ia Trillos, Li and
Venkatraman \cite{gt_li_venkatraman} improve both rates to
$O(\varepsilon^2\log(1/\varepsilon))$ at the same bandwidth.

In practice the points are observed with noise, which changes the distances
between them and may move them off the manifold. We ask how much the first
eigenvalues of the graph Laplacian and their eigenvectors can change when
each point is moved by at most $\delta$.
We consider the unnormalized graph Laplacians of the points and of their
noisy observations, built with the same kernel and the same bandwidth, and
call them the clean and noisy graph Laplacians. For arbitrary noise bounded by
$\delta$, we prove that, for each fixed $l$, if $\delta/\varepsilon$ is
small enough, then with high probability, the $l$-th eigenvalue of the noisy
graph Laplacian differs by $O(\delta/\varepsilon)$ from that of the clean
one. If moreover the $l$-th and $(l+1)$-st eigenvalues of the weighted
Laplace--Beltrami operator are distinct, then the orthogonal projection onto
the subspace spanned by the first $l$ eigenvectors of the noisy graph
Laplacian differs by $O(\sqrt{\delta/\varepsilon})$ from the corresponding
projection for the clean one. Combined with the rate $O(\varepsilon)$ of
Calder and Garc\'ia Trillos \cite{calder_gt} for clean points, our bounds
show that the first eigenvalues of the noisy graph Laplacian approximate those of the
weighted Laplace--Beltrami operator at rate
$O(\varepsilon+\delta/\varepsilon)$, and that the first eigenvectors of the
noisy graph Laplacian approximate the eigenfunctions of this operator at
rate $O(\varepsilon+\sqrt{\delta/\varepsilon})$.

Noise that is not small compared with the bandwidth can sometimes be
handled through its structure. Several results concern random noise in high
ambient dimension. El Karoui \cite{el_karoui} shows that, if the noise is
spherical, the kernel matrix of the noisy points is, to first order, the
kernel matrix of the clean points with a different kernel. Ding and Wu
\cite{ding_wu} study how close the eigenvalues of the kernel matrix of the
noisy points are to those of the kernel matrix of the clean points,
depending on the bandwidth and on the variance of the clean points relative
to that of the noise. Landa, Coifman and Kluger \cite{landa_coifman_kluger}
show that a Gaussian kernel matrix with the doubly stochastic normalization
changes little under noise whose covariance may differ from point to point.
Cheng and Landa \cite{cheng_landa} prove pointwise consistency under outlier
noise for a graph Laplacian with an approximate version of this
normalization.

Noise that is not small compared with the bandwidth can also be handled by
denoising before building the graph
\cite{aizenbud,gt_sanz_alonso_yang,hein_maier,puchkin_spokoiny,sober_levin}.
Each denoised point may depend on all the observations, so the errors that
remain after denoising need not be independent. Our results allow such
dependence, since they hold for arbitrary noise bounded by a deterministic
$\delta$.
Section~\ref{subsec:intro-denoising} explains how our results apply after
denoising.

Section~\ref{subsec:intro-comparison} compares our results with those of
Garc\'ia Trillos, Sanz-Alonso and Yang \cite{gt_sanz_alonso_yang}, of
Burago, Ivanov and Kurylev \cite{burago_2019}, and of Divol \cite{divol}.

\subsection{Setting and main results}

We use the setting of Section~\ref{sec:setup}. The clean points
$x_1,\dots,x_n$ are i.i.d.\ samples from a density $\rho$ on an
$m$-dimensional manifold $\mathcal M\subset\mathbb R^d$. For a
non-increasing Lipschitz kernel profile $\eta$ supported in $[0,1]$, the
scaled unnormalized graph Laplacian $L$ acts on $u\in\mathbb R^n$ by
\[
  (Lu)_i=\frac{1}{n\varepsilon^{m+2}}\sum_{j=1}^n\omega_{ij}(u_i-u_j),
  \qquad \omega_{ij}=\eta\!\left(\frac{\|x_i-x_j\|}{\varepsilon}\right).
\]
The noisy graph Laplacian $\widetilde L$ is built in the same way from
$y_i=x_i+\delta_i$, where $\|\delta_i\|\leq\delta$ for a deterministic
bound $\delta$, possibly depending on $n$. The displacements are otherwise
arbitrary and may depend on the clean points.

Write $\lambda_l$, $\widetilde\lambda_l$ and $\lambda_l^{\mathcal M}$ for
the $l$-th smallest eigenvalues of $L$, of $\widetilde L$ and of the weighted
Laplace--Beltrami operator $\Delta_\rho$ defined in
Section~\ref{sec:setup}.
Fix $l$ and $p\geq1$, and choose
$\varepsilon=C_\varepsilon(\log n/n)^{1/(m+4)}$ with $C_\varepsilon$ large
enough, as in \cite{calder_gt}. For $\delta/\varepsilon$ small enough and
$n$ large enough, the following bounds hold with probability at least
$1-O(n^{-p})$.
Theorem~\ref{thm:main} shows that $\widetilde\lambda_l$ is close to
$\lambda_l$:
\[
  |\lambda_l-\widetilde\lambda_l|=O(\delta/\varepsilon).
\]

Let $P_l$ and $\widetilde P_l$ be the
orthogonal projections onto the subspaces spanned by the first $l$
eigenvectors of $L$ and of $\widetilde L$. If
$\lambda_l^{\mathcal M}<\lambda_{l+1}^{\mathcal M}$, whether or not
$\lambda_l^{\mathcal M}$ is simple, Theorem~\ref{thm:projector-headline}
shows that $\widetilde P_l$ is close to $P_l$:
\[
  \|P_l-\widetilde P_l\|_{\mathrm{HS},n}=O(\sqrt{\delta/\varepsilon}),
\]
where $\|\cdot\|_{\mathrm{HS},n}$ is the Hilbert--Schmidt norm, equal to
the Frobenius norm (Section~\ref{sec:results}).

For each bound, the event depends only on the clean points, so the bound
holds for all displacements with $\|\delta_i\|\leq\delta$ at once
(Lemma~\ref{lem:denoise-uniform}).

Corollary~\ref{cor:continuum-eig} and
Theorem~\ref{thm:eigenvector-headline} combine these bounds with the
clean-sample results of \cite{calder_gt} to give the eigenvalue and
eigenvector convergence rates stated above.

\subsection{Idea of the proofs}\label{subsec:intro-idea}

We bracket both graph Laplacians between two graph Laplacians on the clean
points. Garc\'ia Trillos, Sanz-Alonso and Yang
\cite[Lemma~18]{gt_sanz_alonso_yang} use this idea for the indicator kernel.
Moving each point by at most $\delta$ changes the distance between
two points by at most $2\delta$. Since $\eta$ is non-increasing, shifting its argument
$\|x_i-x_j\|/\varepsilon$ up or down by $2\delta/\varepsilon$, with
negative arguments replaced by zero, gives lower and upper bounds for both
the clean and the noisy weights. The graph Laplacians on the clean points
with these weights, scaled as $L$, are the comparison operators $A$ and
$B$, and
\[
  A\preceq L\preceq B,
  \qquad
  A\preceq\widetilde L\preceq B,
\]
for every point cloud and all displacements with $\|\delta_i\|\leq\delta$,
where $S\preceq T$ means
$\langle Su,u\rangle_n\leq\langle Tu,u\rangle_n$ for every $u\in\mathbb R^n$,
and $\langle u,v\rangle_n=n^{-1}\sum_iu_iv_i$ is the empirical inner
product, with norm $\|u\|_n$. By the
Courant--Fischer principle \cite[Theorem~4.2.6, p.~236]{horn_johnson}, the
$l$-th eigenvalues of $L$ and $\widetilde L$ therefore both lie between
those of $A$ and $B$.

Write $\omega_{ij}^+\leq\omega_{ij}^-$ for the weights of $A$ and $B$. These
weights differ by $O(\delta/\varepsilon)$, and only for pairs of clean
points at distance at most $\varepsilon+2\delta$; with high probability each
point has $O(n\varepsilon^m)$ such neighbours. For every vector $u$,
\[
  \langle(B-A)u,u\rangle_n
  =\frac{1}{2n^2\varepsilon^{m+2}}
   \sum_{i,j}\bigl(\omega_{ij}^--\omega_{ij}^+\bigr)(u_i-u_j)^2,
\]
so $\langle(B-A)u,u\rangle_n=O(\delta/\varepsilon)$ if $\|u\|_n=1$ and the
entries $u_i$ and $u_j$ differ by $O(\varepsilon)$ whenever
$\|x_i-x_j\|\leq\varepsilon+2\delta$. For a unit
eigenvector $u$ of a graph Laplacian on i.i.d.\ samples, with eigenvalue below
a fixed cutoff, the Lipschitz
estimate of Calder, Garc\'ia Trillos and Lewicka
\cite[Theorem~2.3]{calder_lip} bounds $|u_i-u_j|$ by a constant multiple of
the geodesic distance between $x_i$ and $x_j$ plus the bandwidth. The operator
$A$ is built on the clean points, its shifted kernels satisfy the
hypotheses of this estimate, and the constants of the estimate can be chosen
uniformly in the shift (Section~\ref{subsec:aux-lipschitz} and
Appendix~\ref{app:cl-audit}). With high probability, therefore, the entries
of each of the first $l$ unit eigenvectors of $A$ differ by $O(\varepsilon)$
whenever $\|x_i-x_j\|\leq\varepsilon+2\delta$. Since $l$ is fixed, the same
holds for every unit vector $u$ in their span, so
$\langle(B-A)u,u\rangle_n=O(\delta/\varepsilon)$ for every such $u$. The
Courant--Fischer principle
\cite[Theorem~4.2.6, p.~236]{horn_johnson} then bounds the
difference between the $l$-th eigenvalues of $A$ and $B$, and hence
$|\lambda_l-\widetilde\lambda_l|$, by $O(\delta/\varepsilon)$
(Section~\ref{sec:eigenvalue}).

For projectors, let $P_l^+$ be the orthogonal projection onto the subspace
spanned by the first $l$ eigenvectors of $A$. The squared distances
$\|P_l-P_l^+\|_{\mathrm{HS},n}^2$ and
$\|\widetilde P_l-P_l^+\|_{\mathrm{HS},n}^2$ are bounded by a constant multiple
of the maximum of $\langle(B-A)u,u\rangle_n$ over the unit vectors $u$ in this
subspace, divided by the spectral gap between the $l$-th and $(l+1)$-st
eigenvalues of $A$. The maximum is $O(\delta/\varepsilon)$ by the preceding
paragraph. If
$\lambda_l^{\mathcal M}<\lambda_{l+1}^{\mathcal M}$, this spectral gap is
bounded away from zero with high probability. Taking square roots and applying
the triangle inequality give Theorem~\ref{thm:projector-headline}
(Section~\ref{sec:projector}); we do not claim that the square root is sharp.
For $l\geq2$, the proof of Theorem~\ref{thm:eigenvector-headline} applies the
projector bound at the spectral gaps immediately below and above
$\lambda_l^{\mathcal M}$ and combines the results with the proof of
\cite[Theorem~2.7(2)]{calder_gt}.

\subsection{Comparison with previous work}\label{subsec:intro-comparison}

Coifman and Lafon \cite{coifman_lafon} observe that, for a smooth kernel, a
diffusion map remains a valid approximation as long as the bandwidth is larger
than the size of the noise, and refer to spectral perturbation theory for the
eigenvalues and eigenfunctions. We compare our results with later quantitative
bounds, first for the eigenvalues and then for the eigenvectors, where our main
contribution lies.

For the eigenvalues, Burago, Ivanov and Kurylev \cite[Lemma~5.6]{burago_2019}
treat two metrics on the same set with the same measure, the situation of clean
and noisy points. For the Laplacian with the indicator kernel, normalized by the
measure of the ball around each point, they show that, under conditions on the
measure, every eigenvalue changes by a relative amount $O(\delta/\varepsilon)$ if
the distances change by at most $2\delta$. Divol \cite{divol} builds on their
notion of closeness of measures, which holds when mass is moved over distances at
most $\delta$. He replaces their indicator kernel by continuous kernels and their
conditions on the measure by a coarse Poincar\'e inequality, which he shows
empirical measures inherit with high probability. Under his geometric and
smallness hypotheses, he bounds the relative change of the first eigenvalues
between two close measures \cite[Proposition~3.3]{divol}.

Garc\'ia Trillos, Sanz-Alonso and Yang \cite{gt_sanz_alonso_yang} show that, for
i.i.d.\ centred noise normal to the manifold and bounded by $\delta$, averaging
nearby observations can reduce the error in the distances between nearby points
from $2\delta$ to order $\delta^{3/2}$, which allows a smaller bandwidth for a
given accuracy. For the indicator kernel, they bound each eigenvalue of the graph
on the averaged points between multiples of those of two graphs on the clean
points, at a smaller and a larger bandwidth \cite[Lemma~18]{gt_sanz_alonso_yang}.
Combining this bracket with a rate for clean points from
\cite{gt_gerlach_hein_slepcev}, they bound the eigenvalues of the graph on the
averaged points against those of a weighted Laplace--Beltrami operator
\cite[Theorem~3]{gt_sanz_alonso_yang}. Corollary~\ref{cor:eig-bracket} carries
this bracket over from the indicator kernel to non-increasing Lipschitz kernels,
and together with the estimate of $\langle(B-A)u,u\rangle_n$ above gives
Theorem~\ref{thm:main}.

For the eigenvectors, Burago, Ivanov and Kurylev \cite[Section~6]{burago_2019}
observe that a small perturbation generically splits a multiple eigenvalue. For
eigenvalues much smaller than $\varepsilon^{-2}$, they describe each eigenvector,
after transport to the perturbed space, as close to a linear combination of
eigenvectors with nearby eigenvalues. Divol \cite[Proposition~6.8]{divol} gives
an eigenfunction bound under his geometric and smallness hypotheses, for an
eigenvalue separated from all the others by a relative gap, and hence simple.
The bound controls the distance from an eigenfunction of one measure to the
corresponding eigenspace of a close measure, and involves the difference of the
two Laplacians applied to the eigenfunction. Garc\'ia Trillos, Sanz-Alonso and
Yang remark that statements about the convergence of eigenvectors, ``or better
yet, spectral projections'', are possible \cite[Remark~19]{gt_sanz_alonso_yang}.

Our main contribution is a bound on the change of the
subspace spanned by the first $l$ eigenvectors under arbitrary bounded noise, with
multiple eigenvalues allowed. If $\delta/\varepsilon$ is small enough and
$\lambda_l^{\mathcal M}<\lambda_{l+1}^{\mathcal M}$, whether or not
$\lambda_l^{\mathcal M}$ is simple, Theorem~\ref{thm:projector-headline} shows
that, with high probability, $\widetilde P_l$ is within
$O(\sqrt{\delta/\varepsilon})$ of $P_l$. This bound depends on the noise and the
bandwidth only through their ratio $\delta/\varepsilon$, as in the observation of
Coifman and Lafon \cite{coifman_lafon}. For clean points against the continuum,
eigenvector bounds that allow a multiple eigenvalue are known for discretizations
of a manifold \cite[Theorem~4]{burago_2014} and for samples \cite{calder_gt}.
Combined with the results of \cite{calder_gt}, our bound gives a rate for the
eigenvectors of the noisy graph Laplacian against the eigenfunctions of
$\Delta_\rho$ (Theorem~\ref{thm:eigenvector-headline}).

\subsection{Application to denoising}\label{subsec:intro-denoising}

Our results need only one property of the noise: every noisy point lies
within a deterministic distance $\delta$ of its clean point, with
$\delta/\varepsilon$ small enough. They therefore apply to the output of
any denoising method that guarantees this with high probability, even when
the noise in the observations is larger than the bandwidth. Our bounds are
$O(\delta/\varepsilon)$ for the eigenvalues and $O(\sqrt{\delta/\varepsilon})$
for the projectors. So if $\delta=O(\varepsilon^2)$, the eigenvalues of the
graph Laplacian on the denoised points converge at the rate $O(\varepsilon)$
that Calder and Garc\'ia Trillos \cite{calder_gt} prove for the clean graph,
and if $\delta=O(\varepsilon^3)$, so do the eigenvectors.

Appendix~\ref{app:denoising} works this out for the projection algorithm of
Aizenbud and Sober \cite{aizenbud}, for observations uniform in a tube of
fixed radius around $\mathcal M$ (Theorem~\ref{thm:denoised-spectrum}).
Puchkin and Spokoiny \cite{puchkin_spokoiny} bound, with high probability,
the largest distance from a denoised point to its clean point, for bounded
noise whose size tends to zero. A related approach is the local averaging of
Garc\'ia Trillos, Sanz-Alonso and Yang \cite{gt_sanz_alonso_yang}, which
reduces the error in the distances between nearby points
(Section~\ref{subsec:intro-comparison}).

\subsection{Organization of the paper}

Section~\ref{sec:setup} fixes the setting and Section~\ref{sec:results}
states the results. Section~\ref{sec:bracket} constructs the bracket,
Section~\ref{sec:eigenvalue} proves the eigenvalue bounds and
Section~\ref{sec:projector} the projector and eigenvector bounds.
Section~\ref{sec:numerics} reports numerical experiments on the circle $S^1$ and
the unit sphere $S^2$. Appendix~\ref{app:cl-audit}
verifies that the constants of the Lipschitz estimate do not depend on the
shift, Appendix~\ref{app:denoising} contains the denoising application, and
Appendix~\ref{app:cgt-hypotheses} states and verifies the hypotheses of the
results of \cite{calder_gt} that we use.
\section{Setup}\label{sec:setup}

Let $\mathcal M\subset\mathbb R^d$ be a compact, connected, orientable,
smooth $m$-dimensional manifold without boundary, equipped with the
Riemannian structure induced by the ambient space, and define $d_{\mathcal M}$
as its geodesic distance. Let
$\operatorname{Vol}_{\mathcal M}$ be its volume measure and let
$\mu=\rho\,d\operatorname{Vol}_{\mathcal M}$ be a probability measure,
where $\rho\in C^{2,\vartheta}(\mathcal M)$ for some
$\vartheta\in(0,1]$, that is, $\rho$ is twice continuously differentiable and
its second derivatives are $\vartheta$-H\"older continuous, and
\[
  0<\rho_{\min}\leq\rho\leq\rho_{\max}<\infty.
\]
Let $\{x_i\}_{i=1}^n$ be sampled i.i.d.\ from the probability measure $\mu$.
Then the vector
$(x_1,\dots,x_n)$ has distribution $\mathbb P_n:=\mu^{\otimes n}$, the product
measure on the $n$-fold product
$\mathcal M^n:=\mathcal M\times\cdots\times\mathcal M$. For an event
$\mathcal E\subseteq\mathcal M^n$ we write
$\mathcal E^c:=\mathcal M^n\setminus\mathcal E$ for its complement. Let
$\{y_i\}_{i=1}^n \subset \mathbb{R}^d$ be the noisy points, with
$y_i = x_i + \delta_i$, where the displacements $\delta_i\in\mathbb R^d$ satisfy
$\|\delta_i\| \leq \delta$ for every $i$, for some deterministic $\delta\geq0$.
The displacements are otherwise arbitrary; in particular, they may depend on the
clean points.

The clean and noisy graph Laplacians defined below, and the auxiliary graph
Laplacians of Section~\ref{sec:bracket}, are all built from a fixed kernel
profile $\eta : [0,\infty) \to [0,\infty)$ that satisfies
\begin{equation}\label{eq:kernel-hyp}
\begin{aligned}
&\text{(K1) } \mathrm{supp}(\eta) \subseteq [0,1] \text{ and } \eta(1)=0;\\
&\text{(K2) } \eta \text{ is non-increasing};\\
&\text{(K3) } \eta \text{ is } C_\eta\text{-Lipschitz on } [0,1];\\
&\text{(K4) } \int_{B(0,1)\subset\mathbb{R}^m} \eta(|w|)\,dw \;=\; 1;\\
&\text{(K5) } \eta(1/2)>0.
\end{aligned}
\end{equation}
In (K3), $C_\eta$ is a positive constant, and $\eta$ is $C_\eta$-Lipschitz on
$[0,1]$ if $|\eta(s)-\eta(t)|\leq C_\eta|s-t|$ for all $s,t\in[0,1]$. The
kernel is the radial function $w\mapsto\eta(|w|)$ on $\mathbb R^m$, which
appears in (K4). With the hypotheses on $\mathcal M$ and $\rho$ above,
(K1)--(K3) and (K5) imply the hypotheses that \cite[Section~2]{calder_gt}
places on the manifold, the density and the kernel, so its results apply to the graph
Laplacian $L$ defined below. The Lipschitz estimate of
\cite{calder_lip} is applied in Section~\ref{sec:bracket} to kernel profiles
built from $\eta$, and its hypotheses are stated and checked there. We define
the second moment of the kernel in the direction of $e_1$, the first standard
basis vector of $\mathbb R^m$, by
\[
  \sigma_\eta \;:=\; \int_{B(0,1)} \langle w, e_1\rangle^2\, \eta(|w|)\,dw,
\]
and $\sigma_\kappa$ in the same way for any other kernel profile $\kappa$.
By (K2) and (K5), $\eta\geq\eta(1/2)>0$ on $[0,1/2]$. Hence $\eta(0)>0$, and
$\sigma_\eta>0$ because the integrand of $\sigma_\eta$ is positive on
$B(0,1/2)$ off the hyperplane $\{\langle w,e_1\rangle=0\}$.

Define the weighted Laplace--Beltrami operator on a smooth function $f$ by
\[
  \Delta_\rho f
  \;:=\; -\frac{1}{2\rho}\operatorname{div}
                 \bigl(\rho^2\nabla f\bigr),
\]
where the divergence and gradient are intrinsic to $\mathcal M$ (see,
e.g., \cite[equation~(2.12)]{calder_gt}). The operator
$\Delta_\rho$ is positive semidefinite on $L^2(\mathcal M,\mu)$, and its
spectrum consists of eigenvalues of finite multiplicity that tend to infinity
\cite[Section~2.2.1]{calder_gt}. We
denote by $\lambda_l^{\mathcal M}$ the $l$-th eigenvalue of $\Delta_\rho$.
Eigenvalues are indexed here and below by
$l=1,2,\dots$ in non-decreasing order and repeated according to multiplicity,
so that $\lambda_1^{\mathcal M}\leq\lambda_2^{\mathcal M}\leq\cdots$. We write
$\gamma_l^{\mathcal M}:=\lambda_{l+1}^{\mathcal M}-\lambda_l^{\mathcal M}$ for
the $l$-th spectral gap of $\Delta_\rho$.

For a bandwidth $\varepsilon > 0$ and every pair $1\leq i,j\leq n$,
define the clean weights $\omega_{ij}$ and the noisy weights
$\widetilde{\omega}_{ij}$ by
\[
  \omega_{ij} \;=\; \eta\!\left(\tfrac{\|x_i - x_j\|}{\varepsilon}\right),
  \qquad
  \widetilde{\omega}_{ij} \;=\; \eta\!\left(\tfrac{\|y_i - y_j\|}{\varepsilon}\right).
\]
The corresponding (unnormalized) graph Laplacians are defined by
\[
  L \;=\; \tfrac{1}{n\varepsilon^{m+2}}(D - W),
  \qquad
  \widetilde{L} \;=\; \tfrac{1}{n\varepsilon^{m+2}}(\widetilde{D} - \widetilde{W}),
\]
where the weight matrix $W=(\omega_{ij})_{i,j=1}^n$ and the diagonal degree
matrix $D=\operatorname{diag}(\sum_{j=1}^n\omega_{ij})_{i=1}^n$ are $n\times n$
matrices, and $\widetilde W$, $\widetilde D$ are defined in the same way from
the $\widetilde\omega_{ij}$. Define $X_n:=\{x_1,\dots,x_n\}$ and identify
each $u\in\mathbb R^n$ with the function $x_i\mapsto u_i$ on $X_n$. For
$f:\mathcal M\to\mathbb R$, define $f|_{X_n}$ as the vector with entries $f(x_i)$.
We denote the empirical $L^2$-norm, that is, the norm of $L^2(X_n)$,
by $\|u\|_n^2 := \tfrac{1}{n}\sum_{i=1}^n u_i^2$ and the associated inner
product by $\langle u,v\rangle_n := \tfrac{1}{n}\sum_{i=1}^n u_i v_i$.
The Dirichlet energies associated with $L$ and $\widetilde{L}$ are defined by
\[
  b_\varepsilon(u) \;=\; \frac{1}{n^2 \varepsilon^{m+2}}
    \sum_{i,j} \omega_{ij}(u_i - u_j)^2,
  \qquad
  \widetilde{b}_\varepsilon(u) \;=\; \frac{1}{n^2 \varepsilon^{m+2}}
    \sum_{i,j} \widetilde{\omega}_{ij}(u_i - u_j)^2.
\]
We denote by $\lambda_l$ and $\widetilde{\lambda}_l$ the $l$-th eigenvalue of
$L$ and of $\widetilde{L}$, indexed by the convention above. Since $L$ and
$\widetilde{L}$ are symmetric $n\times n$ matrices, $l$ ranges over
$1,\dots,n$.

Calder and Garc\'ia Trillos \cite[Theorem~2.4(2)]{calder_gt} show that, for
fixed $l$ and $p\geq1$, if $\varepsilon=C_\varepsilon(\log n/n)^{1/(m+4)}$
with $C_\varepsilon$ large enough, then for $n$ large enough, with probability
at least $1-O(n^{-p})$,
\[
  \bigl|\lambda_l-\sigma_\eta\lambda_l^{\mathcal M}\bigr|
  \;\leq\;C_l'\,\varepsilon,
\]
where $C_l'$ depends only on $l,\mathcal M,\rho$ and $\eta$. The factor $\sigma_\eta$ appears because $L$
acts on smooth functions approximately as $\sigma_\eta\Delta_\rho$
\cite[Theorem~3.3]{calder_gt}. Our results use the same bandwidth.
Theorem~\ref{thm:main} bounds the difference between the noisy and the clean
eigenvalues, and Corollary~\ref{cor:continuum-eig} combines
Theorem~\ref{thm:main} with this estimate.

We denote the relative noise scale by
\begin{equation}\label{eq:a-def}
  a \;:=\; \frac{2\delta}{\varepsilon}.
\end{equation}
The bound $\delta$ may depend on $n$. In the results of
Section~\ref{sec:results}, $\varepsilon\to0$ as $n\to\infty$, and
$\delta/\varepsilon$ is assumed to be sufficiently small.
\section{Main Results}\label{sec:results}

\subsection{Eigenvalue stability}

\begin{theorem}\label{thm:main}
Assume the setup of Section~\ref{sec:setup} and fix $l,p \geq 1$. Let
$\varepsilon = C_\varepsilon(\log n/n)^{1/(m+4)}$, where $C_\varepsilon$ is large
enough depending on $p,l,\mathcal M,\rho$ and $\eta$, and assume that
$\delta/\varepsilon$ is sufficiently small, depending only on
$l,\mathcal M,\rho$ and $\eta$. Then, for $n$ large enough,
with probability at least $1-O(n^{-p})$,
\[
  |\lambda_l - \widetilde\lambda_l|
  \;\leq\; K_l\, \frac{\delta}{\varepsilon}
  \;=\; O\!\left(\frac{\delta}{\varepsilon}\right).
\]
Here $K_l$ is the constant \eqref{eq:Kl-explicit} of Lemma~\ref{lem:aux-gap}. It
is deterministic and depends only on $l,\mathcal M,\rho$ and $\eta$.
\end{theorem}

\begin{corollary}\label{cor:continuum-eig}
Assume the setup of Section~\ref{sec:setup} and fix $l,p \geq 1$. Let
$\varepsilon = C_\varepsilon(\log n/n)^{1/(m+4)}$, where $C_\varepsilon$ is large
enough depending on $p,l,\mathcal M,\rho$ and $\eta$, and assume that
$\delta/\varepsilon$ is sufficiently small, depending only on
$l,\mathcal M,\rho$ and $\eta$. Then, for $n$ large enough,
with probability at least $1-O(n^{-p})$,
\[
  \bigl|\widetilde\lambda_l-\sigma_\eta\lambda_l^{\mathcal M}\bigr|
  \;\leq\; K_l\,\frac{\delta}{\varepsilon} \;+\; C_l'\,\varepsilon
  \;=\; O\!\left(\frac{\delta}{\varepsilon}\right)+O(\varepsilon).
\]
Here $K_l$ is the constant of Theorem~\ref{thm:main}, and $C_l'$ is the constant
from \cite[Theorem~2.4(2)]{calder_gt} quoted in Section~\ref{sec:setup},
depending only on $l,\mathcal M,\rho$ and $\eta$.
In particular the right-hand side is $O(\varepsilon)$ when
$\delta=O(\varepsilon^2)$.
\end{corollary}

\subsection{Eigenvector stability}
Define $P_l$ as the orthogonal projection of $\mathbb R^n$ onto the subspace
spanned by $l$ orthonormal eigenvectors of $L$ with eigenvalues
$\lambda_1,\dots,\lambda_l$, which we call the first $l$ eigenvectors of $L$, and
$\widetilde P_l$ in the same way for $\widetilde L$. We call $P_l$ and
$\widetilde P_l$ spectral projectors. When $\lambda_l=\lambda_{l+1}$, $P_l$
depends on the choice of eigenvectors, and likewise $\widetilde P_l$ when
$\widetilde\lambda_l=\widetilde\lambda_{l+1}$; the results below hold for every
choice. We measure the distance between projectors in the Hilbert--Schmidt
norm of operators $T$ on $(\mathbb R^n,\langle\cdot,\cdot\rangle_n)$,
\[
  \|T\|_{\mathrm{HS},n}^2\;:=\;\sum_{i=1}^n\|Te_i\|_n^2,
\]
where $e_1,\dots,e_n$ is any $\langle\cdot,\cdot\rangle_n$-orthonormal basis of
$\mathbb R^n$; the sum does not depend on the basis. The vectors $\sqrt n$
times the standard basis vectors form such a basis, so
$\|T\|_{\mathrm{HS},n}$ equals the Frobenius norm of the matrix of $T$.

\begin{theorem}\label{thm:projector-headline}
Assume the setup of Section~\ref{sec:setup} and fix $l,p \geq 1$ such that the
spectral gap $\gamma_l^{\mathcal M}$ is strictly positive. Let
$\varepsilon = C_\varepsilon(\log n/n)^{1/(m+4)}$, where $C_\varepsilon$ is large
enough depending on $p,l,\mathcal M,\rho$ and $\eta$, and assume that
$\delta/\varepsilon$ is sufficiently small, depending only on
$l,\mathcal M,\rho$ and $\eta$. Then, for $n$ large enough, with probability at
least $1-O(n^{-p})$,
\[
  \|P_l - \widetilde P_l\|_{\mathrm{HS},n}
  \;\leq\; D_l\,\sqrt{\frac{\delta}{\varepsilon}}
  \;=\; O\!\left(\sqrt{\frac{\delta}{\varepsilon}}\right),
  \qquad
  D_l:=4\sqrt{\frac{K_l\,l}{\sigma_\eta\gamma_l^{\mathcal M}}},
\]
where $K_l$ is the constant of Theorem~\ref{thm:main}.
\end{theorem}

\begin{theorem}\label{thm:eigenvector-headline}
Assume the setup of Section~\ref{sec:setup} and fix $l\ge2$, $p\geq1$. Write
\[
  \lambda_{l_-}^{\mathcal M}=\cdots=\lambda_l^{\mathcal M}
  =\cdots=\lambda_{l_+}^{\mathcal M},
  \qquad l_-\leq l\leq l_+,
\]
for the multiplicity block of $\lambda_l^{\mathcal M}$, so that $l_-$ and $l_+$
are the smallest and the largest index $k$ with
$\lambda_k^{\mathcal M}=\lambda_l^{\mathcal M}$. Hence the two spectral gaps
$\gamma^{\mathcal M}_{l_--1}$ and $\gamma^{\mathcal M}_{l_+}$ bounding the block
are strictly positive, and $l_-\geq2$ since $\lambda_1^{\mathcal M}$ is simple.
Let $\varepsilon = C_\varepsilon(\log n/n)^{1/(m+4)}$, where $C_\varepsilon$ is
large enough depending on $p,l,\mathcal M,\rho$ and $\eta$, and assume that
$\delta/\varepsilon$ is sufficiently small, depending only on
$l,\mathcal M,\rho$ and $\eta$. Then, for $n$ large enough, with probability at
least $1-O(n^{-p})$, for every
eigenvector $\widetilde u_l$ of $\widetilde L$ with eigenvalue
$\widetilde\lambda_l$ and $\|\widetilde u_l\|_n=1$ there is an eigenfunction
$f_l\in L^2(\mathcal M,\mu)$ with $\Delta_\rho f_l=\lambda_l^{\mathcal M}f_l$
and $\|f_l\|_{L^2(\mu)}=1$ such that
\[
  \|\widetilde u_l-f_l|_{X_n}\|_n
  \;=\; O\!\left(\sqrt{\frac{\delta}{\varepsilon}}+\varepsilon\right).
\]
In particular the right-hand side is $O(\varepsilon)$ when
$\delta=O(\varepsilon^3)$.
\end{theorem}

\begin{remark}[The index $l=1$]\label{rem:index-one}
Theorem~\ref{thm:eigenvector-headline} is stated for $l\geq2$ because its
proof uses the spectral gap $\gamma^{\mathcal M}_{l_--1}$, which is not defined
when $l_-=1$. The case $l=1$ is trivial. Since $\mathcal M$ is
connected, $\lambda_1^{\mathcal M}=0$ is simple and its unit eigenfunctions are the
constants $\pm1$. Under the hypotheses of Corollary~\ref{cor:continuum-eig} with
$l=2$, for $n$ large and with probability at least $1-O(n^{-p})$,
$\widetilde\lambda_1$ is simple as well,
so every unit noisy eigenvector $\widetilde u_1$ at index~$1$ is constant and
$\|\widetilde u_1-f_1|_{X_n}\|_n=0$ for the matching constant eigenfunction $f_1$.
The proof is at the end of Section~\ref{subsec:proof-thm-eigenvector}.
\end{remark}
\section{The Bracketing Construction}\label{sec:bracket}

In this section we introduce the bracketing construction, which reduces the
comparison of the clean and noisy graph Laplacians to a comparison between two
auxiliary graph Laplacians on the clean samples.
Since $\|\delta_i\|\leq\delta$ for every $i$, noise changes each pairwise
distance by at most $2\delta=a\varepsilon$, and hence the argument of $\eta$ in each
weight by at most $a$. Because $\eta$ is non-increasing, increasing or decreasing that
argument by $a$ gives lower and upper bounds for both the clean and the noisy weights.
We use these bounds to define auxiliary weights
satisfying
\begin{equation}\label{eq:weight-bracket}
  \omega_{ij}^+ \leq \omega_{ij} \leq \omega_{ij}^-,
  \qquad
  \omega_{ij}^+ \leq \widetilde{\omega}_{ij} \leq \omega_{ij}^-
\end{equation}
for $1 \leq i,j \leq n$.
Section~\ref{sec:eigenvalue} transfers this common bracket to Dirichlet
energies and eigenvalues, with the appropriate bandwidths.

Keeping the auxiliary graphs on $\{x_i\}\subset\mathcal M$ allows us to
use the estimates of \cite{calder_lip}, whose
sampling hypotheses need not hold for the noisy cloud. We express the
shifted weights through kernel profiles $\eta_a^\pm$ supported in $[0,1]$ and
bandwidths $\varepsilon^\pm=(1\mp a)\varepsilon$. The rest of this section
verifies the weight ordering and shows that, after normalization, these
kernel profiles satisfy the hypotheses of the Lipschitz estimate of
\cite{calder_lip}, stated as (CL1) and (CL2) in
Section~\ref{subsec:aux-lipschitz}, with bounds that do not depend on $a$ when
$a$ is small enough.

\begin{definition}\label{def:aux}
For a parameter $a\in[0,1)$, define the shift maps
\[
  \phi^+(t) \;:=\; a + (1-a)t,
  \qquad
  \phi^-(t) \;:=\; \bigl((1+a)t - a\bigr)_+,
\]
where $(\,\cdot\,)_+ := \max(\cdot, 0)$, and the shifted kernel profiles
$\eta^\pm_a := \eta \circ \phi^\pm$, that is,
\[
  \eta^+_a(t) \;=\; \eta\bigl(a + (1-a)t\bigr),
  \qquad
  \eta^-_a(t) \;=\; \eta\bigl(((1+a)t - a)_+\bigr).
\]
The auxiliary Laplacians use the value $a=2\delta/\varepsilon$ of
\eqref{eq:a-def}, which we assume is less than $1$. For this value, set
\[
  \varepsilon^+ \;:=\; (1-a)\varepsilon,
  \qquad
  \varepsilon^- \;:=\; (1+a)\varepsilon,
\]
and define the auxiliary weights
\[
  \omega_{ij}^+ \;=\; \eta^+_a\!\left(\tfrac{\|x_i - x_j\|}{\varepsilon^+}\right)
                 \;=\; \eta\!\left(a + \tfrac{\|x_i - x_j\|}{\varepsilon}\right),
\]
\[
  \omega_{ij}^- \;=\; \eta^-_a\!\left(\tfrac{\|x_i - x_j\|}{\varepsilon^-}\right)
                 \;=\; \eta\!\left(\bigl(\tfrac{\|x_i - x_j\|}{\varepsilon}
                                          - a\bigr)_+\right).
\]
Let $W^\pm=(\omega_{ij}^\pm)_{i,j=1}^n$ and
$D^\pm=\operatorname{diag}(\sum_{j=1}^n\omega_{ij}^\pm)_{i=1}^n$ be the
$n\times n$ weight and diagonal degree matrices of these weights,
and let $L^\pm = \tfrac{1}{n(\varepsilon^\pm)^{m+2}}(D^\pm - W^\pm)$ be
the corresponding graph Laplacians on $\{x_i\}$, with Dirichlet
energies
\[
  b^\pm(u) \;=\; \frac{1}{n^2 (\varepsilon^\pm)^{m+2}}
                 \sum_{i,j} \omega_{ij}^\pm (u_i - u_j)^2.
\]
Let $\lambda_l^\pm$ denote the $l$-th eigenvalue of $L^\pm$.
\end{definition}
\subsection{The weight bracket}
In Lemmas~\ref{lem:basic} and~\ref{lem:weights} we verify the weight ordering
\eqref{eq:weight-bracket}.

\begin{lemma}\label{lem:basic}
For all $i,j$,
\[
  \bigl|\,\|x_i-x_j\|-\|y_i-y_j\|\,\bigr|
  \;\leq\;\|\delta_i-\delta_j\|
  \;\leq\;2\delta\;=\;a\,\varepsilon.
\]
\end{lemma}

\begin{proof}
Since $y_i - y_j = (x_i - x_j) + (\delta_i - \delta_j)$, the reverse triangle
inequality gives the first bound, and
$\|\delta_i - \delta_j\| \leq \|\delta_i\| + \|\delta_j\| \leq 2\delta$ the
second; $2\delta = a\varepsilon$ by \eqref{eq:a-def}.
\end{proof}

\begin{lemma}\label{lem:weights}
For all $i,j$, $\omega_{ij}^+ \leq \omega_{ij} \leq \omega_{ij}^-$ and
$\omega_{ij}^+ \leq \widetilde{\omega}_{ij} \leq \omega_{ij}^-$.
\end{lemma}

\begin{proof}
Since $\eta$ is non-increasing on $[0,\infty)$ and $a \geq 0$,
\[
  \eta\!\left(a + \tfrac{\|x_i - x_j\|}{\varepsilon}\right)
  \;\leq\;
  \eta\!\left(\tfrac{\|x_i - x_j\|}{\varepsilon}\right)
  \;\leq\;
  \eta\!\left(\bigl(\tfrac{\|x_i - x_j\|}{\varepsilon} - a\bigr)_+\right),
\]
which is exactly $\omega_{ij}^+ \leq \omega_{ij} \leq \omega_{ij}^-$.
For the noisy weights, Lemma~\ref{lem:basic} gives
\[
  \tfrac{\|x_i - x_j\|}{\varepsilon} - a
  \;\leq\;
  \tfrac{\|y_i - y_j\|}{\varepsilon}
  \;\leq\;
  \tfrac{\|x_i - x_j\|}{\varepsilon} + a.
\]
The positive-part map $t \mapsto (t)_+ = \max(t,0)$ is non-decreasing, so
applying it to the left inequality and using
$\bigl(\tfrac{\|y_i - y_j\|}{\varepsilon}\bigr)_+ = \tfrac{\|y_i - y_j\|}{\varepsilon}$
gives
\[
  \Bigl(\tfrac{\|x_i - x_j\|}{\varepsilon} - a\Bigr)_+
  \;\leq\;
  \tfrac{\|y_i - y_j\|}{\varepsilon}
  \;\leq\;
  \tfrac{\|x_i - x_j\|}{\varepsilon} + a.
\]
Since $\eta$ is non-increasing, applying it reverses the chain and yields
$\omega_{ij}^+ \leq \widetilde{\omega}_{ij} \leq \omega_{ij}^-$.
\end{proof}

\subsection{Properties of the shifted kernel profiles}
Lemmas~\ref{lem:structural}--\ref{lem:uniform-pos} collect the properties of the
shifted kernel profiles that are used below. They hold for every value of the
parameter $a$ in a stated range, not only for $a=2\delta/\varepsilon$. We need
this because $a=2\delta/\varepsilon$ changes with $n$ and $\delta$, and
Lemmas~\ref{lem:uniform-adm} and~\ref{lem:uniform-cl} need bounds that are
uniform in $a$.
\begin{lemma}\label{lem:structural}
For every $a \in [0, 1)$, the shifted kernel profiles $\eta^+_a$ and $\eta^-_a$
satisfy:
\begin{enumerate}
  \item\label{it:nonneg} $\eta^\pm_a : [0,\infty) \to [0,\infty)$;
  \item\label{it:dec} $\eta^\pm_a$ is non-increasing on $[0,\infty)$;
  \item\label{it:supp} $\mathrm{supp}(\eta^\pm_a) \subseteq [0,1]$ and
        $\eta^\pm_a(1) = 0$;
  \item\label{it:lip} the restriction of $\eta^\pm_a$ to $[0,1]$ is
        Lipschitz, with constant at most $(1-a)\,C_\eta$ for $\eta^+_a$ and
        $(1+a)\,C_\eta$ for $\eta^-_a$, where $C_\eta$ is the Lipschitz
        constant of $\eta$ in (K3).
\end{enumerate}
\end{lemma}

\begin{proof}
The map $\phi^+$ of Definition~\ref{def:aux} is affine with slope $1-a\geq0$,
and $\phi^-$ is zero on
$[0,a/(1+a)]$ and affine with slope $1+a$ on $[a/(1+a),\infty)$.
Thus both maps are non-decreasing on $[0,\infty)$, with
$\phi^\pm(0)\in[0,a]$ and $\phi^\pm(1)=1$. Hence
$\phi^\pm([0,1]) \subseteq [0,1]$ and $\phi^\pm(t) \geq 1$ for $t \geq 1$.

\emph{(\ref{it:nonneg})} Inherited from $\eta \geq 0$.

\emph{(\ref{it:dec})} A non-decreasing map composed with a
non-increasing map is non-increasing.

\emph{(\ref{it:supp})} For $t \geq 1$, $\phi^\pm(t) \geq 1$, so
$\eta^\pm_a(t) = \eta(\phi^\pm(t)) = 0$ by (K1).

\emph{(\ref{it:lip})} On $[0,1]$, $\phi^+$ has Lipschitz constant
$1-a$, and $\phi^-$ has Lipschitz constant $1+a$ (the positive-part
map is $1$-Lipschitz, applied to an affine map of slope $1+a$).
Since $\phi^\pm([0,1]) \subseteq [0,1]$ and $\eta$ is $C_\eta$-Lipschitz on
$[0,1]$ by (K3),
\[
  \bigl|\eta^\pm_a(t) - \eta^\pm_a(t')\bigr|
  \;\leq\; C_\eta\,\bigl|\phi^\pm(t) - \phi^\pm(t')\bigr|
  \;\leq\; (1 \mp a)\,C_\eta\,|t - t'|. \qedhere
\]
\end{proof}

\begin{lemma}\label{lem:phi-bound}
For every $a \in [0,1)$ and every $t \in [0,1]$,
\[
  \bigl|\phi^\pm(t) - t\bigr| \;\leq\; a.
\]
\end{lemma}

\begin{proof}
$\phi^+(t) - t = a(1 - t) \in [0, a]$. For $\phi^-$: when $t \in
[a/(1+a), 1]$, $\phi^-(t) - t = a(t - 1) \in [-a/(1+a), 0]$; when $t
\in [0, a/(1+a))$, $\phi^-(t) - t = -t \in (-a/(1+a), 0]$.
\end{proof}

\begin{lemma}\label{lem:eta-pointwise}
The kernel profile $\eta$ is $C_\eta$-Lipschitz on all of $[0,\infty)$,
and for every $a \in [0,1)$ and every $t \in [0,\infty)$,
\[
  \bigl|\eta^\pm_a(t) - \eta(t)\bigr| \;\leq\; C_\eta\,a.
\]
\end{lemma}

\begin{proof}
For $s \in [0,1]$ and $t \geq 1$,
$\eta(t) = \eta(1) = 0$, hence
\[
  |\eta(s) - \eta(t)| \;=\; |\eta(s) - \eta(1)|
  \;\leq\; C_\eta\,|s - 1| \;\leq\; C_\eta\,|s - t|,
\]
where the last step uses $|s - t| = t - s \geq 1 - s$. Combined with
the (K3) Lipschitz bound on $[0,1]$ and the trivial bound on
$[1,\infty)$, $\eta$ is $C_\eta$-Lipschitz on $[0,\infty)$.

For $t \geq 1$, $\eta^\pm_a(t) =
\eta(t) = 0$ (Lemma~\ref{lem:structural}~\eqref{it:supp}). For $t \in
[0,1]$, by Lemma~\ref{lem:phi-bound} and the global Lipschitz bound for $\eta$,
\[
  |\eta^\pm_a(t) - \eta(t)|
  \;=\; |\eta(\phi^\pm(t)) - \eta(t)|
  \;\leq\; C_\eta\,|\phi^\pm(t) - t|
  \;\leq\; C_\eta\,a. \qedhere
\]
\end{proof}

\begin{lemma}\label{lem:moment-stab}
For every $a \in [0, 1)$,
\[
  \biggl|\int_{\mathbb{R}^m} \eta^\pm_a(|w|)\,dw \,-\, 1\biggr|
  \;\leq\; C_\eta\,|B_1^m|\,a,
  \qquad
  |\sigma_{\eta^\pm_a} - \sigma_\eta|
  \;\leq\; C_\eta\,|B_1^m|\,a,
\]
where $C_\eta$ is the Lipschitz constant of $\eta$ in (K3) and $|B_1^m|$ is the
volume of the unit ball in $\mathbb{R}^m$. In
particular, for $a \leq \sigma_\eta/(2C_\eta\,|B_1^m|)$,
\[
  \sigma_{\eta^\pm_a} \;\geq\; \tfrac{1}{2}\sigma_\eta \;>\; 0.
\]
\end{lemma}

\begin{proof}
Both $\eta$ and $\eta^\pm_a$ are supported on $\{|w| \leq 1\}$, so the
integrand $|\eta^\pm_a(|w|) - \eta(|w|)|$ vanishes for $|w| > 1$. By
Lemma~\ref{lem:eta-pointwise} and the unit normalization (K4),
\[
  \biggl|\int \eta^\pm_a(|w|)\,dw - 1\biggr|
  \;\leq\; \int_{|w|\le 1} |\eta^\pm_a(|w|) - \eta(|w|)|\,dw
  \;\leq\; C_\eta\,|B_1^m|\,a.
\]
For the second bound, using $\langle w, e_1\rangle^2 \leq 1$
on $|w| \leq 1$,
\[
  |\sigma_{\eta^\pm_a} - \sigma_\eta|
  \;\leq\; \int_{|w|\le 1} \langle w, e_1\rangle^2
            \,|\eta^\pm_a(|w|) - \eta(|w|)|\,dw
  \;\leq\; \int_{|w|\le 1} C_\eta\,a\,dw
  \;=\; C_\eta\,|B_1^m|\,a. \qedhere
\]
\end{proof}

\begin{lemma}\label{lem:uniform-pos}
Let $c_0 := \eta(0)$, $t_0 := \min(c_0/(2C_\eta),\,1/2)$, and
$a_0 := t_0/2$. For every $a \in [0, a_0)$,
\[
  \eta^\pm_a(t) \;\geq\; \tfrac{c_0}{2}
  \qquad \text{for every } t \in [0, t_0/2].
\]
In particular, $\eta^\pm_a(|w|) \geq c_0/2$ for every $w$ in the open
Euclidean ball $S_0 := B(0, t_0/2) \subset \mathbb{R}^m$, and neither the ball
nor the bound depends on $a$.
\end{lemma}

\begin{proof}
By (K3) and the definition of $t_0$, $\eta(t) \geq \eta(0) - C_\eta t = c_0
- C_\eta t \geq c_0/2$ for every $t \in [0, t_0]$.

For $t \in [0, t_0/2]$ and $a \in [0, a_0)$,
\[
  \phi^+(t) \;=\; a + (1-a)t \;\leq\; a + t \;\leq\; \tfrac{t_0}{2} + \tfrac{t_0}{2} \;=\; t_0,
\]
so $\eta^+_a(t) = \eta(\phi^+(t)) \geq c_0/2$.

For $\eta^-_a$, let $t \in [0, t_0/2]$. Since $(1+a)t-a\leq t$ for $t\leq1$,
we have $0\leq\phi^-(t)\leq t\leq t_0$, and again
$\eta^-_a(t) = \eta(\phi^-(t)) \geq c_0/2$.
\end{proof}

\subsection{Lipschitz estimates for the auxiliary Laplacians}\label{subsec:aux-lipschitz}
Section~\ref{sec:eigenvalue} needs a Lipschitz estimate for the eigenvectors of
$L^+$. Calder, Garc\'ia Trillos and Lewicka \cite[Theorem~2.3]{calder_lip} prove
such an estimate for the graph Laplacian on the clean points built, as $L$ is in
Section~\ref{sec:setup}, from a single kernel profile
$\kappa:[0,\infty)\to[0,\infty)$ in place of $\eta$. Besides hypotheses on
$\mathcal M$ and $\rho$ that those of Section~\ref{sec:setup} imply, the
estimate assumes the following.
\begin{enumerate}
  \item[(CL1)] $\kappa$ is non-increasing, supported in $[0,1]$ and Lipschitz on
        $[0,1]$, and $\int_{\mathbb R^m}\kappa(|w|)\,dw=1$.
  \item[(CL2)] There are an open set $S\subset B(0,1)\subset\mathbb R^m$ and a
        constant $c>0$ such that $\kappa(|w|)\geq c$ for every $w\in S$.
\end{enumerate}
These are \cite[(H1)(iii)]{calder_lip}, with the support condition that
\cite[Section~2]{calder_lip} imposes throughout, and \cite[(9.1)]{calder_lip}.
The lower bound (CL2) is used only in \cite[Section~9]{calder_lip}, for pairs
of points at geodesic distance at most $3\varepsilon$. There it gives
\cite[Theorem~9.1]{calder_lip}, which bounds $|f(x)-f(y)|$ for such a pair by a
factor $\theta<1$ times the supremum of $|f(x')-f(y')|$ over nearby pairs in a
larger ball, up to error terms. The proof of the global estimate
\cite[Theorem~10.1]{calder_lip} needs $\theta<1$ to absorb that supremum. The constants of
\cite[Theorem~2.3]{calder_lip}, including its smallness threshold for the
bandwidth, depend on $\kappa$.

Our kernel profiles depend on $a$, and $a=2\delta/\varepsilon$ changes with $n$
and $\delta$, so we need these constants to be the same for all of them. For a
family of profiles this requires, among other things, that $\theta$ stay below
a fixed number smaller than $1$, which is why Lemma~\ref{lem:uniform-pos} bounds
the shifted profiles from below on one ball, by one constant.
Lemma~\ref{lem:uniform-adm} normalizes the shifted kernel profiles and shows that
the normalized profiles satisfy (CL1), and (CL2) with $S=S_0$, with bounds that
do not depend on $a$. Lemma~\ref{lem:uniform-cl} is the estimate of
\cite[Theorem~2.3]{calder_lip} for the normalized auxiliary Laplacians
$\bar L^\pm$ of Definition~\ref{def:aux-normalized}. What it adds to that
theorem is that the constants depend only on $\mathcal M$, $\rho$ and $\eta$,
which Appendix~\ref{app:cl-audit} verifies.

\begin{lemma}\label{lem:uniform-adm}
Let
\[
  a_\star \;:=\; \min\Bigl(\tfrac{t_0}{2},\, \tfrac{1}{2C_\eta\,|B_1^m|},\,
                            \tfrac{\sigma_\eta}{2C_\eta\,|B_1^m|},\,
                            \tfrac12\Bigr),
\]
with $t_0$ as in Lemma~\ref{lem:uniform-pos}, $C_\eta$ the Lipschitz constant
of (K3) and $|B_1^m|$ the volume of the unit ball in $\mathbb R^m$. It depends
only on $m$ and $\eta$. For every $a\in[0,a_\star)$ and for both signs, the
quantities
\[
  I_a^\pm := \int_{\mathbb{R}^m}\eta_a^\pm(|w|)\,dw,
  \qquad
  \bar\eta_a^\pm := \frac{1}{I_a^\pm}\eta_a^\pm
\]
are well-defined, with $\tfrac12\le I_a^\pm\le\tfrac32$, and the normalized kernel
profiles $\bar\eta_a^\pm$ satisfy the following bounds, none of which depends
on $a$:
\begin{align*}
  \mathrm{supp}(\bar\eta_a^\pm) &\subseteq [0,1],
  &\int_{\mathbb{R}^m}\bar\eta_a^\pm(|w|)\,dw &= 1, \\
  \mathrm{Lip}_{[0,1]}(\bar\eta_a^\pm) &\leq 3C_\eta,
  &\|\bar\eta_a^\pm\|_\infty &\leq 2\eta(0), \\
  \frac{\sigma_\eta}{3} \leq \sigma_{\bar\eta_a^\pm} &\leq 3\sigma_\eta,
  &\bar\eta_a^\pm(|w|) &\geq \frac{c_0}{3}
       \quad\text{on } S_0.
\end{align*}
Here $\mathrm{Lip}_{[0,1]}$ is the smallest Lipschitz constant of the
restriction to $[0,1]$, $\|\cdot\|_\infty$ is the supremum, and $c_0$ and the
ball $S_0=B(0,t_0/2)$ are as in Lemma~\ref{lem:uniform-pos}. In particular,
$\bar\eta_a^\pm$ satisfies (CL1), and (CL2) with $S=S_0$ and $c=c_0/3$.
\end{lemma}

\begin{proof}
Let $a<a_\star$. Since $a<t_0/2=a_0$, $a\le1/(2C_\eta|B_1^m|)$ and
$a\le\sigma_\eta/(2C_\eta|B_1^m|)$, Lemmas~\ref{lem:moment-stab}
and~\ref{lem:uniform-pos} give
\[
  \tfrac12\le I_a^\pm\le \tfrac32,
  \qquad
  \tfrac12\sigma_\eta\le \sigma_{\eta_a^\pm}\le \tfrac32\sigma_\eta,
  \qquad
  \eta_a^\pm(|w|)\ge c_0/2\quad\text{on }S_0.
\]
In particular $I_a^\pm>0$. The normalized profile $\bar\eta_a^\pm$ is the
shifted profile $\eta_a^\pm$ divided by the constant $I_a^\pm$. So it has unit
mass, and it keeps the support, non-negativity and monotonicity of
$\eta_a^\pm$ from Lemma~\ref{lem:structural}. Its Lipschitz constant on
$[0,1]$, its supremum and its lower bound on $S_0$ are those of $\eta_a^\pm$
multiplied by $(I_a^\pm)^{-1}\in[\tfrac23,2]$. By
Lemma~\ref{lem:structural}~\eqref{it:lip} and $a<\tfrac12$,
$\mathrm{Lip}_{[0,1]}(\eta_a^\pm)\le(1+a)C_\eta\le\tfrac32C_\eta$, and
$\|\eta_a^\pm\|_\infty=\eta(\phi^\pm(0))\le\eta(0)$ because $\eta$ is
non-increasing. Hence $\mathrm{Lip}_{[0,1]}(\bar\eta_a^\pm)\le3C_\eta$,
$\|\bar\eta_a^\pm\|_\infty\le2\eta(0)$, and
$\bar\eta_a^\pm(|w|)\ge\tfrac23\cdot\tfrac{c_0}{2}=\tfrac{c_0}{3}$ on $S_0$.
The second moment is divided by $I_a^\pm$ as well:
\[
  \sigma_{\bar\eta_a^\pm}
  = \frac{\sigma_{\eta_a^\pm}}{I_a^\pm}
  \in \left[\frac{\sigma_\eta/2}{3/2},\frac{3\sigma_\eta/2}{1/2}\right]
  = [\sigma_\eta/3,3\sigma_\eta].
\]
Finally, these properties give (CL1), and the lower bound on the open ball
$S_0\subset B(0,1)$ gives (CL2).
\end{proof}

\begin{remark}[Standing assumption]\label{rem:standing}
From now on, $a$ is again the relative noise scale $2\delta/\varepsilon$ of
\eqref{eq:a-def}, and we assume $a<a_\star$. In the results of
Section~\ref{sec:results}, this is part of the assumption that
$\delta/\varepsilon$ is sufficiently small.
\end{remark}

\begin{definition}\label{def:aux-normalized}
Denote by $\bar L^\pm:=(I_a^\pm)^{-1}L^\pm$ the \emph{normalized auxiliary
Laplacians}. They are the graph Laplacians on $\{x_i\}$ built from the normalized
profiles $\bar\eta_a^\pm$ at bandwidths $\varepsilon^\pm$, since $\bar\eta_a^\pm
=(I_a^\pm)^{-1}\eta_a^\pm$ rescales every auxiliary weight by the same factor.
The normalization is that of the mass of the kernel profile, as in (K4). Like
$L$, the operators $\bar L^\pm$ are unnormalized graph Laplacians in the sense
of Section~\ref{sec:setup}.
\end{definition}

\begin{lemma}
\label{lem:uniform-cl}
There exist constants $\mathcal C_{\mathrm{CL}}^{\mathrm{unif}}$,
$c_\star^{\mathrm{unif}}$, $C_{\mathbb P}$, $c_{\mathbb P}>0$, depending only on
$\mathcal M$, $\rho$ and the kernel profile $\eta$, and in particular not on
$a$, with the following property. Let $s\in\{+,-\}$ be one of the two signs and
$\Lambda>0$ a cutoff, and assume
$\varepsilon^s\le c_\star^{\mathrm{unif}}/(\Lambda+1)$. There is an event
$\mathcal E^s_{n,1}(\Lambda)$, depending only on $x_1,\ldots,x_n$, of
probability at least
\[
  1-C_{\mathbb P}(\varepsilon^s)^{-6m}
     \exp\!\bigl(-c_{\mathbb P}n(\varepsilon^s)^{m+4}\bigr)
   -2n\exp\!\bigl(-c_{\mathbb P}n(\Lambda+1)^{-m}\bigr)
\]
such that, given $\mathcal E^s_{n,1}(\Lambda)$, every nonzero $f$ on $X_n$ with
$\|\bar L^s f\|_\infty<\Lambda\|f\|_\infty$ satisfies
\[
  |f(x_i)-f(x_j)|
  \le \mathcal C_{\mathrm{CL}}^{\mathrm{unif}}(\Lambda+1)^{m+1}\,\|f\|_{1,n}
      \bigl(d_{\mathcal M}(x_i,x_j)+\varepsilon^s\bigr)
  \qquad\text{for all }i,j,
\]
where $\|f\|_\infty:=\max_i|f(x_i)|$ and $\|f\|_{1,n}:=\frac1n\sum_i|f(x_i)|$.
\end{lemma}

\begin{proof}
Fix the sign $s$. By Definition~\ref{def:aux-normalized}, $\bar L^s$ is the graph
Laplacian on the clean points built from the profile $\bar\eta_a^s$ at the
bandwidth $\varepsilon^s$, scaled as $L$ is in Section~\ref{sec:setup}. The
clean points are i.i.d.\ samples from $\mu$, and by Lemma~\ref{lem:uniform-adm}
the profile $\bar\eta_a^s$ satisfies (CL1), and (CL2) with $S=S_0$ and
$c=c_0/3$. So \cite[Theorem~2.3]{calder_lip}, with $\bar\eta_a^s$ and
$\varepsilon^s$ in place of $\eta$ and $\varepsilon$, gives an event and an
estimate of the asserted form. Its hypothesis on $f$ is $\lambda_f<\Lambda$,
where $\lambda_f:=\|\bar L^sf\|_\infty/\|f\|_\infty$, which is the condition in
the statement. The theorem has two smallness conditions on the bandwidth:
$\varepsilon^s\le c/(\Lambda+1)$ for a constant $c$, and $\varepsilon^s$ below a
threshold.

It remains to choose the constants independently of $a$. The constants of
\cite[Theorem~2.3]{calder_lip}, including $c$ and the threshold, depend on
$\mathcal M$ and $\rho$, and on the profile only through four quantities: an
upper bound on its supremum, an upper bound on its Lipschitz constant on
$[0,1]$, a lower bound on its second moment, and a positive lower bound on a
fixed ball. Proposition~\ref{prop:cl-uniform} in Appendix~\ref{app:cl-audit}
verifies this. By Lemma~\ref{lem:uniform-adm},
\[
  \|\bar\eta_a^s\|_\infty\le 2\eta(0),\qquad
  \mathrm{Lip}_{[0,1]}(\bar\eta_a^s)\le 3C_\eta,\qquad
  \sigma_{\bar\eta_a^s}\ge\tfrac13\sigma_\eta,\qquad
  \bar\eta_a^s(|w|)\ge\tfrac{c_0}{3}\ \text{for}\ w\in S_0,
\]
and these bounds are the same for every $a<a_\star$ and both signs. The source
uses the lower bound through \cite[Proposition~9.2]{calder_lip}, which is
stated for an annulus, so Appendix~\ref{app:cl-audit} applies it on the annulus
$B(0,t_0/2)\setminus\bar B(0,t_0/4)$, which lies in $S_0$ and inherits the
bound. Proposition~\ref{prop:cl-uniform}, applied to the family of
Lemma~\ref{lem:our-family-admissible}, therefore gives a single prefactor
$\mathcal C_{\mathrm{CL}}^{\mathrm{unif}}$, a single pair
$C_{\mathbb P},c_{\mathbb P}$ of probability constants, and a single constant
$c_\star^{\mathrm{unif}}$ below both $c$ and the threshold, valid for every
$a<a_\star$ and both signs. Since $\Lambda>0$, the assumption
$\varepsilon^s\le c_\star^{\mathrm{unif}}/(\Lambda+1)$ implies both smallness
conditions.
\end{proof}

\section{Eigenvalue Stability}\label{sec:eigenvalue}
In this section we prove Theorem~\ref{thm:main} and Corollary~\ref{cor:continuum-eig}.
By Lemma~\ref{lem:dirichlet} and Corollary~\ref{cor:eig-bracket}, for every point cloud both
$\lambda_l$ and $\widetilde\lambda_l$ lie between $(1-a)^{m+2}\lambda_l^+$ and
$(1+a)^{m+2}\lambda_l^-$. Proposition~\ref{prop:aux-eig} shows that, given an event
$\mathcal E_n$ of high probability, the difference between the upper and the lower bound,
and hence $|\lambda_l-\widetilde\lambda_l|$, is $O(\delta/\varepsilon)$.
The proof of Proposition~\ref{prop:aux-eig} applies the Courant--Fischer principle
\cite[Theorem~4.2.6, p.~236]{horn_johnson} for $\lambda_l^-$ to the span of the first $l$
eigenvectors of $L^+$. For $u$ in this span, Lemma~\ref{lem:aux-gap} bounds
${(1+a)^{m+2}b^-(u)-(1-a)^{m+2}b^+(u)}$ by a constant multiple of $a\,\|u\|_n^2$, using the
Lipschitz estimate of \cite{calder_lip} for these eigenvectors (Proposition~\ref{prop:hpe}).
This estimate applies because $L^+$ is built on the clean points (Section~\ref{sec:bracket}).
Section~\ref{subsec:hpe} defines $\mathcal E_n$, bounds its probability and proves the
estimate given $\mathcal E_n$, and Section~\ref{subsec:proof-thm-main} completes the proofs.

\subsection{Bracketing of Dirichlet energies and eigenvalues}
\label{subsec:energy-bracket}
Lemma~\ref{lem:dirichlet} transfers the weight bracket of Lemma~\ref{lem:weights} to the
Dirichlet energies, and Corollary~\ref{cor:eig-bracket} transfers it to the eigenvalues by
the Courant--Fischer principle \cite[Theorem~4.2.6, p.~236]{horn_johnson}.

\begin{lemma}\label{lem:dirichlet}
For every $u \in \mathbb{R}^n$,
\[
  (1-a)^{m+2}\, b^+(u) \;\leq\; b_\varepsilon(u) \;\leq\;
  (1+a)^{m+2}\, b^-(u),
\]
and the same inequalities hold with $b_\varepsilon$ replaced by
$\widetilde{b}_\varepsilon$.
\end{lemma}

\begin{proof}
Using $(\varepsilon^\pm)^{m+2} = (1 \mp a)^{m+2}\varepsilon^{m+2}$ and
Lemma~\ref{lem:weights},
\[
  (1-a)^{m+2} b^+(u)
  \;=\; \frac{1}{n^2 \varepsilon^{m+2}}\sum_{i,j} \omega_{ij}^+ (u_i - u_j)^2
  \;\leq\; \frac{1}{n^2 \varepsilon^{m+2}}\sum_{i,j} \omega_{ij}(u_i - u_j)^2
  \;=\; b_\varepsilon(u),
\]
and the upper bound is analogous. The same argument with
$\widetilde{\omega}_{ij}$ in place of $\omega_{ij}$ proves the inequalities
for $\widetilde b_\varepsilon$.
\end{proof}

\begin{corollary}\label{cor:eig-bracket}
For each $l$,
\[
  (1-a)^{m+2}\, \lambda_l^+ \;\leq\; \lambda_l \;\leq\;
  (1+a)^{m+2}\, \lambda_l^-,
\]
and the same inequalities hold for $\widetilde{\lambda}_l$. Consequently
\begin{equation}\label{eq:gap-bound}
  |\lambda_l - \widetilde{\lambda}_l|
  \;\leq\;
  (1+a)^{m+2}\lambda_l^- \,-\, (1-a)^{m+2}\lambda_l^+.
\end{equation}
\end{corollary}

\begin{proof}
By the Courant--Fischer principle \cite[Theorem~4.2.6, p.~236]{horn_johnson}
and the identity $b_\varepsilon(u)=2\langle Lu,u\rangle_n$,
\[
  \lambda_l \;=\; \frac12\min_{V \in \mathcal{G}_l}\; \max_{u \in V\setminus\{0\}}\;
                  \frac{b_\varepsilon(u)}{\|u\|_n^2},
\]
where $\mathcal{G}_l$ is the family of $l$-dimensional subspaces. Apply
Lemma~\ref{lem:dirichlet} inside the min-max. The same argument applies
to $\widetilde{\lambda}_l$. Both $\lambda_l$ and $\widetilde{\lambda}_l$ lie
inside the interval $[(1-a)^{m+2}\lambda_l^+, (1+a)^{m+2}\lambda_l^-]$,
yielding \eqref{eq:gap-bound}.
\end{proof}
\subsection{The high-probability event}\label{subsec:hpe}

Corollary~\ref{cor:eig-bracket} holds for every point cloud. To bound the
right-hand side of \eqref{eq:gap-bound}, Section~\ref{subsec:proof-thm-main}
needs two properties of the clean points that hold only with high
probability: a bound on the number of sample points in each ball
$B(x_i,\varepsilon^-)$, and hence on the number of pairs of sample points at
distance at most $\varepsilon^-$, and a Lipschitz estimate for the first $l$
orthonormal eigenvectors $u$ of $L^+$, which bounds $|u(x_i)-u(x_j)|$ by a
constant multiple of $d_{\mathcal M}(x_i,x_j)+\varepsilon^+$.
Proposition~\ref{prop:aux-lip} below gives the Lipschitz estimate only for
eigenvectors whose eigenvalues lie below a cutoff fixed before the clean
points are drawn. A third property, the estimate of
\cite[Theorem~2.4(2)]{calder_gt} for the $l$-th eigenvalue of the clean
Laplacian $L$, keeps $\lambda_l^+$ below a deterministic constant, and the
cutoff $\Lambda_l$, defined after Proposition~\ref{prop:aux-lip}, is chosen
above twice this constant. Definition~\ref{def:hpe} introduces the corresponding event
$\mathcal E_n$ of the clean points, and Proposition~\ref{prop:hpe} bounds its
probability and proves the Lipschitz estimate given $\mathcal E_n$.

\begin{proposition}\label{prop:aux-lip}
Fix a cutoff $\Lambda>0$ and assume
$\varepsilon^+\le c_\star^{\mathrm{unif}}/(\Lambda+1)$, where
$c_\star^{\mathrm{unif}}$ and the constants
$\mathcal C_{\mathrm{CL}}^{\mathrm{unif}},C_{\mathbb P},c_{\mathbb P}$ below are
those of Lemma~\ref{lem:uniform-cl}. Let
$\mathcal E^+_{n,1}(\Lambda)\subseteq\mathcal M^n$ be the event of
Lemma~\ref{lem:uniform-cl} for the sign $s=+$. The event
$\mathcal E^+_{n,1}(\Lambda)$ depends only on the clean points
$x_1,\ldots,x_n$, and
\begin{equation}\label{eq:lip-event-prob}
  \mathbb P_n\bigl(\mathcal E^+_{n,1}(\Lambda)^c\bigr)
  \;\le\;
  C_{\mathbb P}(\varepsilon^+)^{-6m}\exp\!\bigl(-c_{\mathbb P}n(\varepsilon^+)^{m+4}\bigr)
  +2n\exp\!\bigl(-c_{\mathbb P}n(\Lambda+1)^{-m}\bigr).
\end{equation}
Given $\mathcal E^+_{n,1}(\Lambda)$, the following holds for every eigenvector
$u$ of $L^+$ with $\|u\|_n=1$ whose eigenvalue $\lambda$ satisfies
$(I_a^+)^{-1}\lambda<\Lambda$:
\[
  |u(x_i)-u(x_j)|
  \le \mathcal C_{\mathrm{CL}}^{\mathrm{unif}}(\Lambda+1)^{m+1}
      \bigl(d_{\mathcal M}(x_i,x_j)+\varepsilon^+\bigr)
  \qquad\text{for all }i,j.
\]
\end{proposition}

\begin{proof}
The probability bound is that of Lemma~\ref{lem:uniform-cl}. Let $u$ be as in
the statement. Since $\bar L^+=(I_a^+)^{-1}L^+$, $u$ is an eigenvector of
$\bar L^+$ with eigenvalue $\bar\lambda=(I_a^+)^{-1}\lambda$, so
$\|\bar L^+u\|_\infty=\bar\lambda\|u\|_\infty<\Lambda\|u\|_\infty$.
Lemma~\ref{lem:uniform-cl} therefore applies to $u$, and the bound follows from
$\|u\|_{1,n}\le\|u\|_n=1$.
\end{proof}

From now on we fix $l\geq1$ and take the bandwidth
\begin{equation}\label{eq:scaling}
  \varepsilon \;=\; C_\varepsilon(\log n/n)^{1/(m+4)},
\end{equation}
where $C_\varepsilon>0$ is a constant, fixed independently of $n$. Set
\[
  B_l \;:=\; 2^{m+2}\bigl(\sigma_\eta\,\lambda_l^{\mathcal M} + C_l'\bigr),
  \qquad
  \Lambda_l \;:=\; 2B_l + 1,
\]
where $\lambda_l^{\mathcal M}$ is the $l$-th eigenvalue of $\Delta_\rho$ and
$C_l'$ is the constant of \cite[Theorem~2.4(2)]{calder_gt}. Both are
deterministic, depending only on $l,\mathcal M,\rho$ and the kernel profile
$\eta$. We assume that
${\varepsilon^+\le c_\star^{\mathrm{unif}}/(\Lambda_l+1)}$, the bandwidth
condition of Proposition~\ref{prop:aux-lip} at $\Lambda=\Lambda_l$, and that
$\varepsilon\le1$. Since $\varepsilon^+\le\varepsilon$, the first holds
whenever $\varepsilon\le c_\star^{\mathrm{unif}}/(2B_l+2)$, so under
\eqref{eq:scaling} both hold for all $n$ large enough. Finally, $C_1$ denotes
the constant, depending only on $\mathcal M$ and $\rho$, that the proof of
Proposition~\ref{prop:hpe} obtains from \cite[Corollary~3.3]{calder_lip}.

\begin{definition}[High-probability event]\label{def:hpe}
Let $\mathcal{E}_n \subseteq \mathcal{M}^n$ be the intersection
$\mathcal E_{n,1}\cap\mathcal E_{n,2}\cap\mathcal E_{n,3}$ of the following
three events.
\begin{enumerate}
  \item\label{it:lip-event}\emph{(Lipschitz bound for the eigenvectors
        of $L^+$.)} $\mathcal E_{n,1}:=\mathcal E^+_{n,1}(\Lambda_l)$, the
        event of Proposition~\ref{prop:aux-lip} at $\Lambda=\Lambda_l$.
  \item\label{it:degree-event}\emph{(Degree upper bound at radius
        $\varepsilon^-$.)} $\mathcal E_{n,2}$ is the event that, for every
        $i \in \{1, \ldots, n\}$,
        \[
          \#\{j : \|x_i - x_j\| \leq \varepsilon^-\}
          \;\leq\; C_1\, n\,(\varepsilon^-)^m.
        \]
  \item\label{it:lambda-plus-event}\emph{(Eigenvalue convergence for the
        clean Laplacian $L$.)} $\mathcal E_{n,3}$ is the event that
        $|\lambda_l-\sigma_\eta\lambda_l^{\mathcal M}|\leq C_l'\varepsilon$.
\end{enumerate}
\end{definition}

The constant $C_l'$ and the event
$\mathcal E_{n,3}$ come from \cite[Theorem~2.4(2)]{calder_gt}. Appendix~\ref{app:cgt-hypotheses} states
the hypotheses of that theorem in our notation, chooses the parameters
$\theta$ and $\delta_e$ of \cite{calder_gt} so that they hold for $n$ large
enough (Lemma~\ref{lem:cgt-parameters}), and shows that, for each $p\geq1$,
once $C_\varepsilon$ is large enough,
\begin{equation}\label{eq:cgt-cited}
  \mathbb P_n\bigl(\mathcal E_{n,3}\bigr)
  \;\geq\;1-O(n^{-p}),
  \qquad
  \bigl|\lambda_l-\sigma_\eta\lambda_l^{\mathcal M}\bigr|
  \;\leq\;C_l'\varepsilon
  \quad\text{given }\mathcal E_{n,3}.
\end{equation}

\begin{proposition}\label{prop:hpe}
The event $\mathcal E_n$ of Definition~\ref{def:hpe} has the following
properties.
\begin{enumerate}
  \item\label{it:hpe-prob} Let $p\geq1$. If $C_\varepsilon$ is large enough,
        depending on $p,l,\mathcal M,\rho$ and $\eta$, then
        $\mathbb P_n(\mathcal E_n)\geq1-O(n^{-p})$ for $n$ large enough.
  \item\label{it:hpe-lambda-plus} Given $\mathcal E_n$, $\lambda_l^+\leq B_l$.
  \item\label{it:hpe-lip} Let $u^1,\ldots,u^l$ be $L^2(X_n)$-orthonormal
        eigenvectors of $L^+$ corresponding to its first $l$ eigenvalues.
        Given $\mathcal E_n$, for every $k\le l$ and all $i,j$,
        \[
          |u^k(x_i)-u^k(x_j)|
          \le C_{\mathrm{aux}}(B_l+1)^{m+1}
              \bigl(d_{\mathcal M}(x_i,x_j)+\varepsilon^+\bigr),
        \]
        where
        \begin{equation}\label{eq:caux-def}
          C_{\mathrm{aux}} := 2^{m+1}\,\mathcal C_{\mathrm{CL}}^{\mathrm{unif}}
        \end{equation}
        depends only on $\mathcal M,\rho$ and $\eta$, through the constant of
        Lemma~\ref{lem:uniform-cl}.
\end{enumerate}
\end{proposition}

\begin{proof}
In this proof the constants $C, c > 0$ may differ from
display to display, but each depends only on $\mathcal{M},\rho$ and
the kernel profile $\eta$.

\emph{The probability bound.}
For $\mathcal E_{n,1}$, by \eqref{eq:lip-event-prob} with
$\Lambda=\Lambda_l$, the complement has probability at most
\begin{equation}\label{eq:p1}
  P_1 \;:=\;
    C\,(\varepsilon^+)^{-6m}
      \exp\!\bigl(-c\,n\,(\varepsilon^+)^{m+4}\bigr)
    \;+\;
    2n\,\exp\!\bigl(-c\,n\,(\Lambda_l+1)^{-m}\bigr),
\end{equation}
and these constants do not depend on $a$.

For $\mathcal E_{n,2}$, by
\cite[Corollary~3.3]{calder_lip}, the complement has probability
at most
\begin{equation}\label{eq:p2}
  P_2 \;:=\;
  2n\,\exp\!\bigl(-c\,n\,(\varepsilon^-)^m\bigr).
\end{equation}
\cite[Corollary~3.3]{calder_lip} bounds the number of sample points in a
ball $B(x,\varepsilon^-)$ whose centre $x$ does not depend on the sample points.
Here the centre is itself a sample point, so for each $i$ we condition on
$x_i=x$. The other $n-1$ points are still independent samples from $\mu$,
and the corollary applies to the number of them in $B(x,\varepsilon^-)$.
For $n(\varepsilon^-)^m$ large enough, the term $j=i$ and the change from
$n-1$ to $n$ points are absorbed into $C_1$ and $c$. Integrating over $x$
and taking a union bound over $i$ gives \eqref{eq:p2}.

For $\mathcal E_{n,3}$, set
\[
  P_3:=\mathbb P_n\bigl(\mathcal E_{n,3}^c\bigr).
\]
By a union bound over the three events,
\begin{equation}\label{eq:hpe-prob}
\begin{aligned}
  \mathbb{P}_n(\mathcal{E}_n)
  \;\geq\; 1 - P_1 - P_2 - P_3
  \;=\; 1
  &\;-\; C\,(\varepsilon^+)^{-6m}
         \exp\!\bigl(-c\,n\,(\varepsilon^+)^{m+4}\bigr)
   \;-\; 2n\,\exp\!\bigl(-c\,n\,(\Lambda_l+1)^{-m}\bigr) \\
  &\;-\; 2n\,\exp\!\bigl(-c\,n\,(\varepsilon^-)^m\bigr) \\
  &\;-\; P_3.
\end{aligned}
\end{equation}
By \eqref{eq:cgt-cited} (Appendix~\ref{app:cgt-hypotheses}), $P_3=O(n^{-p})$.
Since $a<a_\star\leq1/2$ (Remark~\ref{rem:standing}), we have
$\tfrac12\varepsilon\leq\varepsilon^\pm\leq\tfrac32\varepsilon$. After increasing $C_\varepsilon$ if
necessary, the remaining terms in \eqref{eq:hpe-prob} are likewise
$O(n^{-p})$. Hence
$\mathbb{P}_n(\mathcal{E}_n) \geq 1-O(n^{-p})$.

\emph{The bound on $\lambda_l^+$.}
By Corollary~\ref{cor:eig-bracket} applied to the clean Laplacian,
$(1-a)^{m+2}\,\lambda_l^+\le\lambda_l$, so
$\lambda_l^+\le(1-a)^{-(m+2)}\,\lambda_l$. Given $\mathcal{E}_n$,
item~\ref{it:lambda-plus-event} of Definition~\ref{def:hpe} gives
$\lambda_l\le\sigma_\eta\,\lambda_l^{\mathcal M}+C_l'$, since $\varepsilon\le1$.
Since $a<a_\star\le1/2$, we have $(1-a)^{-(m+2)}\le2^{m+2}$, whence
$\lambda_l^+\le2^{m+2}\bigl(\sigma_\eta\,\lambda_l^{\mathcal M}+C_l'\bigr)=B_l$.

\emph{The Lipschitz estimate.}
Item~\ref{it:lip-event} of Definition~\ref{def:hpe} is the event
$\mathcal E_{n,1}=\mathcal E^+_{n,1}(\Lambda_l)$ of Proposition~\ref{prop:aux-lip},
with $\Lambda_l=2B_l+1$.
For $k\le l$, the bound $\lambda_l^+\le B_l$ proved above
and $I_a^+\ge\tfrac12$ (Lemma~\ref{lem:uniform-adm}) give
\[
  (I_a^+)^{-1}\lambda_k^+\le 2\lambda_k^+\le 2\lambda_l^+\le 2B_l
  <2B_l+1=\Lambda_l,
\]
so the cutoff condition of Proposition~\ref{prop:aux-lip} is met. Its
bandwidth condition $\varepsilon^+\le c_\star^{\mathrm{unif}}/(\Lambda_l+1)$
is assumed before Definition~\ref{def:hpe}. Applying
the proposition at $\Lambda=\Lambda_l$ and using $\Lambda_l+1=2(B_l+1)$ with
$C_{\mathrm{aux}}=2^{m+1}\mathcal C_{\mathrm{CL}}^{\mathrm{unif}}$ gives the
stated bound.
\end{proof}
\subsection[Proofs of the eigenvalue results]{Proofs of Theorem~\ref{thm:main} and Corollary~\ref{cor:continuum-eig}}\label{subsec:proof-thm-main}

It remains to bound the right-hand side of \eqref{eq:gap-bound}, which
Proposition~\ref{prop:aux-eig} does given $\mathcal E_n$. This section ends
with the proofs of Theorem~\ref{thm:main} and
Corollary~\ref{cor:continuum-eig}. We compare Euclidean and geodesic
distances by \cite[(3.2)]{calder_lip}: there exists $R > 0$, depending only
on $\mathcal{M}$, such that for all $x,y \in \mathcal{M}$ with
$\|x-y\| \leq R/2$,
\[
  \|x-y\| \;\leq\; d_{\mathcal{M}}(x,y)
                    \;\leq\; \|x-y\| + \tfrac{8}{R^2}\|x-y\|^3.
\]
From here on we also assume that $\varepsilon\leq R/3$; under the scaling
\eqref{eq:scaling}, this holds for large $n$.

\begin{lemma}\label{lem:dist-comp}
Let $a<1/2$ and $\varepsilon\leq R/3$. If $x,y\in\mathcal M$ satisfy
$\|x-y\| \leq \varepsilon^- = (1+a)\varepsilon$, then
$d_{\mathcal{M}}(x,y) \leq 3\,\varepsilon^- \leq 5\varepsilon$.
\end{lemma}

\begin{proof}
Since $a<1/2$ and $\varepsilon\leq R/3$, we have
$\varepsilon^-<\tfrac32\varepsilon\leq R/2$, so the hypothesis
$\|x-y\|\leq \varepsilon^-$ puts the pair in the range of
\cite[(3.2)]{calder_lip}. Moreover,
\[
  \frac{8}{R^2}\|x-y\|^3
  = \frac{8}{R^2}\|x-y\|^2\,\|x-y\|
  \leq \frac{8}{R^2}\left(\frac R2\right)^2\varepsilon^-
  = 2\varepsilon^-.
\]
Therefore $d_{\mathcal M}(x,y)\leq 3\varepsilon^-$. Since $a<1/2$,
$\varepsilon^-=(1+a)\varepsilon\leq 3\varepsilon/2$, hence
$d_{\mathcal M}(x,y)\leq 9\varepsilon/2\leq 5\varepsilon$.
\end{proof}

\begin{lemma}\label{lem:aux-gap}
Let $u^1,\ldots,u^l$ be
$L^2(X_n)$-orthonormal eigenvectors of $L^+$ corresponding to its first
$l$ eigenvalues, and set
\[
  S^+ := \mathrm{Span}\{u^1,\ldots,u^l\}.
\]
Then, given $\mathcal{E}_n$, for every $u\in S^+$,
\[
  (1+a)^{m+2}\, b^-(u) \;-\; (1-a)^{m+2}\, b^+(u)
  \;\leq\; K_l\, a\, \|u\|_n^2,
\]
where one may take
\begin{equation}\label{eq:Kl-explicit}
  K_l := 2(3/2)^mC_1C_\eta\,l\,\bigl(6C_{\mathrm{aux}}(B_l+1)^{m+1}\bigr)^2,
\end{equation}
with $C_\eta$ the Lipschitz constant of $\eta$ in (K3), $C_{\mathrm{aux}}$ as in
\eqref{eq:caux-def}, and $B_l$ and $C_1$ as in Section~\ref{subsec:hpe}.
\end{lemma}

\begin{proof}
Write $u = \sum_{k=1}^l \beta_k u^k$, so that $\|u\|_n^2 = \sum_k
\beta_k^2$.

By item~\ref{it:hpe-lip} of Proposition~\ref{prop:hpe} and $\varepsilon^+\leq\varepsilon$,
given $\mathcal{E}_n$ we have, for every $k\le l$ and all $i,j$,
\[
  |u^k(x_i) - u^k(x_j)| \;\leq\; C_{\mathrm{aux}}\,(B_l + 1)^{m+1}
                                  \bigl(d_{\mathcal{M}}(x_i, x_j) + \varepsilon\bigr).
\]
For pairs $(i,j)$ with $\|x_i - x_j\| \leq \varepsilon^-$, Lemma~\ref{lem:dist-comp}
gives $d_{\mathcal{M}}(x_i, x_j) \leq 5\varepsilon$. Substituting,
\begin{equation}\label{eq:eig-close}
  |u^k(x_i) - u^k(x_j)| \;\leq\; 6 C_{\mathrm{aux}}\,(B_l + 1)^{m+1}\,\varepsilon
  \;=:\; M_l\,\varepsilon.
\end{equation}

For such a pair, Cauchy--Schwarz and \eqref{eq:eig-close} give
\[
  |u_i - u_j|
  \;\leq\; \sum_{k=1}^l |\beta_k|\,|u^k(x_i) - u^k(x_j)|
  \;\leq\; M_l\,\varepsilon\,\sum_{k=1}^l |\beta_k|
  \;\leq\; M_l\,\varepsilon\,\sqrt{l}\,\|u\|_n,
\]
hence
\begin{equation}\label{eq:lin-comb}
  (u_i - u_j)^2 \;\leq\; M_l^2\,\varepsilon^2\,l\,\|u\|_n^2.
\end{equation}

Next we compare the weights. The function $\eta$ is
$C_\eta$-Lipschitz on $[0,\infty)$ (Lemma~\ref{lem:eta-pointwise}), so for
every $i,j$, with $r_{ij} := \|x_i - x_j\|/\varepsilon$,
\[
  \omega_{ij}^- - \omega_{ij}^+
  \;=\; \eta\bigl((r_{ij} - a)_+\bigr) - \eta(a + r_{ij})
  \;\leq\; C_\eta\,\bigl|(r_{ij} - a)_+ - (a + r_{ij})\bigr|
  \;\leq\; 2 C_\eta a.
\]
Moreover $\omega_{ij}^- = 0$ whenever $\|x_i - x_j\| > \varepsilon^-$.

Finally, the number of pairs $(i,j)$
with $\|x_i - x_j\| \leq \varepsilon^-$ is at most $n \cdot C_1 n
(\varepsilon^-)^m \leq C_1\,(1+a)^m\,n^2\,\varepsilon^m \leq (3/2)^m C_1\,
n^2\,\varepsilon^m$ given $\mathcal{E}_n$ (item~\ref{it:degree-event} of
Definition~\ref{def:hpe}).
Combining this with the weight bound and \eqref{eq:lin-comb},
\begin{align*}
  (1+a)^{m+2} b^-(u) - (1-a)^{m+2} b^+(u)
  &\;=\; \frac{1}{n^2 \varepsilon^{m+2}}
         \sum_{\substack{i,j \\ \|x_i-x_j\|\le\varepsilon^-}}
         (\omega_{ij}^- - \omega_{ij}^+)(u_i - u_j)^2 \\
  &\;\leq\; \frac{1}{n^2 \varepsilon^{m+2}}
            \cdot (3/2)^m C_1 n^2\varepsilon^m \cdot 2C_\eta a \cdot M_l^2\varepsilon^2\,l\,\|u\|_n^2 \\
  &\;=\; \underbrace{2(3/2)^m C_1 C_\eta\,M_l^2\,l}_{=: K_l}\, a\, \|u\|_n^2. \qedhere
\end{align*}
\end{proof}

\begin{proposition}\label{prop:aux-eig}
Given $\mathcal{E}_n$,
\[
  (1+a)^{m+2}\lambda_l^- \;\leq\; (1-a)^{m+2}\lambda_l^+
                            \,+\, \frac12 K_l\, a,
\]
with $K_l$ as in Lemma~\ref{lem:aux-gap}.
\end{proposition}

\begin{proof}
The subspace $S^+$ of Lemma~\ref{lem:aux-gap} is $l$-dimensional, so the
Courant--Fischer principle, applied as in the proof of
Corollary~\ref{cor:eig-bracket}, gives
\[
  (1+a)^{m+2}\lambda_l^-
  \;=\; \frac12\min_{V\in \mathcal{G}_l} \max_{u \in V\setminus\{0\}}
        \frac{(1+a)^{m+2} b^-(u)}{\|u\|_n^2}
  \;\leq\; \frac12\max_{u \in S^+\setminus\{0\}}
        \frac{(1+a)^{m+2} b^-(u)}{\|u\|_n^2},
\]
where, as there, $\mathcal G_l$ is the family of $l$-dimensional subspaces.
By Lemma~\ref{lem:aux-gap},
\[
  \frac{(1+a)^{m+2} b^-(u)}{\|u\|_n^2}
  \;\leq\; \frac{(1-a)^{m+2} b^+(u)}{\|u\|_n^2} + K_l\,a.
\]
The maximum of the first term over $u \in S^+\setminus\{0\}$ is exactly
$2(1-a)^{m+2}\lambda_l^+$, since $b^+(u)=2\langle L^+u,u\rangle_n$ and $S^+$
is spanned by orthonormal eigenvectors of $L^+$ with eigenvalues
$\lambda_1^+,\dots,\lambda_l^+$. Multiplying the displayed inequality by the
outside factor $1/2$ gives the claim.
\end{proof}

\begin{proof}[Proof of Theorem~\ref{thm:main}]
By Proposition~\ref{prop:hpe}, $\mathbb P_n(\mathcal E_n)\geq1-O(n^{-p})$.
Given $\mathcal E_n$, combine \eqref{eq:gap-bound} with
Proposition~\ref{prop:aux-eig}:
\[
  |\lambda_l - \widetilde{\lambda}_l|
  \;\leq\; (1+a)^{m+2}\lambda_l^- - (1-a)^{m+2}\lambda_l^+
  \;\leq\; \frac12 K_l\, a
  \;\leq\; \frac{K_l\, \delta}{\varepsilon}. \qedhere
\]
\end{proof}

\begin{proof}[Proof of Corollary~\ref{cor:continuum-eig}]
By Proposition~\ref{prop:hpe}, $\mathbb P_n(\mathcal E_n)\geq1-O(n^{-p})$.
Given $\mathcal E_n$, both estimates below hold: the first by the proof of
Theorem~\ref{thm:main}, and the second by \eqref{eq:cgt-cited}, since
$\mathcal E_n$ is contained in the Calder--Garc\'ia Trillos event
$\mathcal E_{n,3}$ of item~\ref{it:lambda-plus-event} of
Definition~\ref{def:hpe}. The triangle inequality gives
\[
  \bigl|\widetilde\lambda_l-\sigma_\eta\lambda_l^{\mathcal M}\bigr|
  \;\leq\;
  \underbrace{\bigl|\widetilde\lambda_l-\lambda_l\bigr|}
             _{\leq\, K_l\, \delta/\varepsilon}
  \;+\;
  \underbrace{\bigl|\lambda_l-\sigma_\eta\lambda_l^{\mathcal M}\bigr|}
             _{\leq\, C_l'\,\varepsilon}
  \;\leq\; K_l\,\frac{\delta}{\varepsilon}+C_l'\,\varepsilon. \qedhere
\]
\end{proof}

\section{Eigenvector Stability}\label{sec:projector}

In this section we prove Theorems~\ref{thm:projector-headline}
and~\ref{thm:eigenvector-headline}, the analogues of Theorem~\ref{thm:main} and
Corollary~\ref{cor:continuum-eig} for eigenvectors.
Theorem~\ref{thm:projector-headline} bounds the distance between the spectral
projectors $P_l$ and $\widetilde P_l$ of Section~\ref{sec:results}, rather than
between single eigenvectors. When
$\lambda_l^{\mathcal M}$ is a multiple
eigenvalue, $L$ has several eigenvalues close to $\lambda_l$,
and a perturbation larger than the differences between them can mix their
eigenvectors.
So an eigenvector of $L$ for $\lambda_l$ need not be close to any
eigenvector of $\widetilde L$. The projector $P_l$ does not have this
problem: Theorem~\ref{thm:projector-headline} needs only
$\gamma_l^{\mathcal M}=\lambda_{l+1}^{\mathcal M}-\lambda_l^{\mathcal M}>0$, and
$\lambda_l^{\mathcal M}$ itself may be a multiple eigenvalue.

The proof of Theorem~\ref{thm:projector-headline} bounds the distances from
$P_l$ and from $\widetilde P_l$ to $P_l^+$, the orthogonal projection onto the
subspace spanned by the first $l$ eigenvectors of $L^+$.
Lemma~\ref{lem:operator-bracket} restates the energy bracket of
Lemma~\ref{lem:dirichlet} as the inequalities
$A\preceq L\preceq B$ and $A\preceq\widetilde L\preceq B$, where
$A=(1-a)^{m+2}L^+$ and $B=(1+a)^{m+2}L^-$ are the comparison operators and
$S\preceq T$ means that $T-S$ is positive semidefinite.
Lemma~\ref{lem:operator-bracket} also bounds $\langle(B-A)u,u\rangle_n$ by
$\tfrac12K_l\,a\,\|u\|_n^2$ for every $u$ in the range of $P_l^+$, given
$\mathcal E_n$, using Lemma~\ref{lem:aux-gap}.
Lemma~\ref{lem:fcpb} and Proposition~\ref{prop:aux-projector} show that $P_l$ and
$\widetilde P_l$ are close to $P_l^+$ if $\tfrac12K_l\,a$ is small compared with the
spectral gap $\Gamma_l^+$ between the $l$-th and $(l+1)$-st eigenvalues of $A$,
defined in \eqref{eq:gap-defs}.
Under \eqref{eq:gap-smallness}, Lemma~\ref{lem:gap-transfer} bounds $\Gamma_l^+$ from below by
$\tfrac12\sigma_\eta\gamma_l^{\mathcal M}$, given the events of
Definition~\ref{def:hpe} at the indices $l$ and $l+1$,
and the triangle inequality then gives Theorem~\ref{thm:projector-headline}.

In Section~\ref{subsec:proof-thm-eigenvector} we apply the bound of
Theorem~\ref{thm:projector-headline} on both sides of the multiplicity block
$l_-,\dots,l_+$ of $\lambda_l^{\mathcal M}$, at the indices $l_--1$ and
$l_+$, and combine it with the
proof of \cite[Theorem~2.7(2)]{calder_gt} to prove
Theorem~\ref{thm:eigenvector-headline}. The claim of
Remark~\ref{rem:index-one} is proved at the end of that section.

\subsection{Notation}

Throughout Section~\ref{sec:projector}, $1\le l<n$ is fixed and $\mathbb R^n$
carries the inner product $\langle\cdot,\cdot\rangle_n$ of $L^2(X_n)$. For an
operator $T$ on $\mathbb R^n$, $T^\ast$ denotes its adjoint; since
$\langle\cdot,\cdot\rangle_n$ is $1/n$ times the Euclidean inner product,
$T^\ast$ is the transpose of $T$. We write $\operatorname{Tr}(T)$ for the
trace of $T$, the sum of the diagonal entries of its matrix in any basis. The
trace does not depend on the inner product: in a
$\langle\cdot,\cdot\rangle_n$-orthonormal basis $e_1,\dots,e_n$, the $i$-th
diagonal entry of the matrix of $T$ is $\langle Te_i,e_i\rangle_n$, so
\begin{equation}\label{eq:trace-basis}
  \operatorname{Tr}(T) \;=\; \sum_{i=1}^n \langle Te_i,e_i\rangle_n.
\end{equation}
The Hilbert--Schmidt norm of Section~\ref{sec:results} satisfies
$\|T\|_{\mathrm{HS},n}^2=\operatorname{Tr}(T^\ast T)$. We define the
operator norm $\|T\|_{\mathrm{op},n}:=\max_{\|u\|_n=1}\|Tu\|_n$, which is the
spectral norm of the matrix of $T$. The operator norm is at most
$\|T\|_{\mathrm{HS},n}$: a unit vector $u$ can be completed to a
$\langle\cdot,\cdot\rangle_n$-orthonormal basis, and then $\|Tu\|_n^2$ is one
of the terms in the sum that defines $\|T\|_{\mathrm{HS},n}^2$ in
Section~\ref{sec:results}. For a subspace $V$ and $v\in\mathbb R^n$ we set
$\operatorname{dist}_n(v,V):=\min_{w\in V}\|v-w\|_n$. For symmetric $S,T$ we
write $S\preceq T$ if $\langle Su,u\rangle_n\leq\langle Tu,u\rangle_n$ for
every $u\in\mathbb R^n$, that is, if $T-S$ is positive semidefinite.

We express the bracket of Lemma~\ref{lem:dirichlet} and
Corollary~\ref{cor:eig-bracket} in terms of the \emph{comparison operators}
\begin{equation}\label{eq:AB-def}
  A\;:=\;(1-a)^{m+2}L^+,
  \qquad
  B\;:=\;(1+a)^{m+2}L^-.
\end{equation}
The eigenvalues of $A$ are $\alpha_k=(1-a)^{m+2}\lambda_k^+$, and we define the
spectral gap of $A$ between its $l$-th and $(l+1)$-st eigenvalues by
\begin{equation}\label{eq:gap-defs}
  \Gamma_l^+ \;:=\; \alpha_{l+1}-\alpha_l.
\end{equation}

Let $P_l^+$ be the orthogonal projection onto the subspace spanned by the
first $l$ eigenvectors of $L^+$, that is, by $l$ orthonormal eigenvectors of
$L^+$ with eigenvalues $\lambda_1^+,\dots,\lambda_l^+$. Since
$(1-a)^{m+2}>0$, these are also eigenvectors of $A$, with eigenvalues
$\alpha_1,\dots,\alpha_l$. If
$\Gamma_l^+>0$, that is, $\lambda_l^+<\lambda_{l+1}^+$, the range of $P_l^+$
is the span of all eigenvectors of $L^+$ with eigenvalue at most
$\lambda_l^+$: this span contains the range, and it has dimension $l$
because only $\lambda_1^+,\dots,\lambda_l^+$ are at most $\lambda_l^+$. So
$P_l^+$ does not depend on the choice of eigenvectors.
The same argument shows that $P_l$ does not depend on the choice when
$\lambda_l<\lambda_{l+1}$, and $\widetilde P_l$ when
$\widetilde\lambda_l<\widetilde\lambda_{l+1}$. If $\Gamma_l^+=0$, then
$P_l^+$ depends on the choice of eigenvectors, and we fix one.

The size of $B-A$ on the range of $P_l^+$ is measured by
\begin{equation}\label{eq:tau-def}
  \tau_l \;:=\;
  \max_{0\neq u\in\operatorname{Range}(P_l^+)}
    \frac{\langle(B-A)u,u\rangle_n}{\|u\|_n^2}.
\end{equation}
For a given point cloud, $\tau_l$ can be computed from $L^+$ and $L^-$, and
Section~\ref{sec:numerics} measures it. The bound $\tau_l\leq\tfrac12K_l\,a$ of
Lemma~\ref{lem:operator-bracket} below cannot be evaluated, because $K_l$
contains the prefactor $\mathcal C_{\mathrm{CL}}^{\mathrm{unif}}$ of
Lemma~\ref{lem:uniform-cl}, which is not available in effective form.

\subsection{Comparing spectral projectors}
Lemma~\ref{lem:operator-bracket} proves $A\preceq L\preceq B$,
$A\preceq\widetilde L\preceq B$ and, given $\mathcal E_n$, the bound
$\tau_l\leq\tfrac12K_l\,a$.
Lemma~\ref{lem:fcpb} shows that if $A$ has a positive spectral gap between its
$l$-th and $(l+1)$-st eigenvalues and $B-A$ is small on the span of the first
$l$ eigenvectors of $A$, then, for every symmetric $M$ with
$A\preceq M\preceq B$, the orthogonal projection onto the subspace spanned by
the first $l$ eigenvectors of $M$ is close to that of $A$ in the
Hilbert--Schmidt norm. Lemma~\ref{lem:fcpb} holds for arbitrary symmetric
$A\preceq B$ on $(\mathbb R^n,\langle\cdot,\cdot\rangle_n)$ and every $l$, not
only for the comparison operators \eqref{eq:AB-def} and the index fixed in
Section~\ref{sec:projector}. The proof of Lemma~\ref{lem:fcpb} uses
Lemma~\ref{lem:hs-identity}, which expresses the Hilbert--Schmidt distance
between two orthogonal projections of the same rank as a trace.

\begin{lemma}\label{lem:operator-bracket}
The comparison operators satisfy
\begin{equation}\label{eq:op-bracket}
  A\;\preceq\;L\;\preceq\;B,
  \qquad
  A\;\preceq\;\widetilde L\;\preceq\;B.
\end{equation}
Moreover, given $\mathcal E_n$, $\tau_l\leq\tfrac12K_l\,a$, with $K_l$ as in
\eqref{eq:Kl-explicit}.
\end{lemma}

\begin{proof}
Since $b_\varepsilon(u)=2\langle Lu,u\rangle_n$,
$\widetilde b_\varepsilon(u)=2\langle\widetilde Lu,u\rangle_n$ and
$b^\pm(u)=2\langle L^\pm u,u\rangle_n$, the definitions \eqref{eq:AB-def}
give
\[
  (1-a)^{m+2}b^+(u)=2\langle A u,u\rangle_n,
  \qquad
  (1+a)^{m+2}b^-(u)=2\langle B u,u\rangle_n.
\]
Hence Lemma~\ref{lem:dirichlet} gives \eqref{eq:op-bracket}, and
$2\langle(B-A)u,u\rangle_n=(1+a)^{m+2}b^-(u)-(1-a)^{m+2}b^+(u)$. For the
bound on $\tau_l$, Lemma~\ref{lem:aux-gap} with
$S^+=\operatorname{Range}(P_l^+)$ therefore gives, given $\mathcal E_n$,
\[
  2\langle(B-A)u,u\rangle_n\leq K_l\,a\,\|u\|_n^2
  \qquad\text{for every }u\in\operatorname{Range}(P_l^+).
\]
Dividing by $2\|u\|_n^2$ for $u\neq0$ and taking the maximum over $u$ gives
$\tau_l\leq\tfrac12K_l\,a$.
\end{proof}

\begin{lemma}\label{lem:hs-identity}
Let $P,Q$ be $\langle\cdot,\cdot\rangle_n$-orthogonal projections of the
same rank $l$. Then
\begin{equation}\label{eq:hs-identity}
  \|P-Q\|_{\mathrm{HS},n}^2
  \;=\;
  2\,\operatorname{Tr}\!\bigl((I-P)Q\bigr)
  \;=\;
  2\bigl(l-\operatorname{Tr}(PQ)\bigr),
\end{equation}
and $\|P-Q\|_{\mathrm{HS},n}^2\leq2l$.
\end{lemma}

\begin{proof}
Since $P,Q$ are self-adjoint idempotents of trace $l$,
\[
  \|P-Q\|_{\mathrm{HS},n}^2
  = \operatorname{Tr}\!\bigl((P-Q)^2\bigr)
  = \operatorname{Tr}(P)+\operatorname{Tr}(Q)-2\operatorname{Tr}(PQ)
  = 2\bigl(l-\operatorname{Tr}(PQ)\bigr),
\]
and $l-\operatorname{Tr}(PQ)=\operatorname{Tr}\bigl((I-P)Q\bigr)$ because
$\operatorname{Tr}(Q)=l$. By cyclicity of the trace and $P^2=P$, $Q^2=Q$,
we have
$\operatorname{Tr}(PQ)=\operatorname{Tr}\bigl((QP)^\ast QP\bigr)
=\|QP\|_{\mathrm{HS},n}^2\geq0$, which gives the bound.
\end{proof}

\begin{lemma}\label{lem:fcpb}
Let $A\preceq B$ be symmetric operators on
$(\mathbb R^n,\langle\cdot,\cdot\rangle_n)$, let
$\alpha_1\leq\cdots\leq\alpha_n$ be the eigenvalues of $A$, and let
$\tau\geq0$. For any integer $l$ with $1\le l<n$, let $P$ be the orthogonal
projection onto the subspace spanned by $l$ orthonormal eigenvectors of $A$
with eigenvalues $\alpha_1,\ldots,\alpha_l$, and assume
\[
  \Gamma:=\alpha_{l+1}-\alpha_l>0
  \qquad\text{and}\qquad
  \langle (B-A)u,u\rangle_n\leq\tau\|u\|_n^2
  \quad\text{for all }u\in\operatorname{Range}(P).
\]
Then for every symmetric $M$ with $A\preceq M\preceq B$, with eigenvalues
$\mu_1\leq\cdots\leq\mu_n$, and every orthogonal projection $Q$ onto the
subspace spanned by $l$ orthonormal eigenvectors of $M$ with eigenvalues
$\mu_1,\ldots,\mu_l$,
\begin{equation}\label{eq:fcpb}
  \|Q - P\|_{\mathrm{HS},n}^2
  \;\leq\; \frac{2\,l\,\tau}{\Gamma}.
\end{equation}
\end{lemma}

\begin{proof}
Let $u^1,\ldots,u^n$ be an orthonormal eigenbasis of $A$ with
$Au^k=\alpha_ku^k$ and $\operatorname{Range}(P)
=\operatorname{Span}\{u^1,\ldots,u^l\}$, put
$r_k:=\langle Qu^k,u^k\rangle_n=\|Qu^k\|_n^2\in[0,1]$, and set
\[
  s\;:=\;\operatorname{Tr}\!\bigl((I-P)Q\bigr)
   \;=\;\sum_{k>l}r_k
   \;=\;\sum_{k\leq l}(1-r_k).
\]
The first equality follows by evaluating the trace in the basis
$\{u^k\}$, and the last because $\sum_k r_k=\operatorname{Tr}(Q)=l$. We bound
$\operatorname{Tr}(AQ)-\operatorname{Tr}(AP)$ from below by $\Gamma s$ and
from above by $l\tau$.

\emph{Lower bound.} Expanding both traces in the basis $\{u^k\}$ and
using $\alpha_k\geq\alpha_{l+1}$ for $k>l$, $\alpha_k\leq\alpha_l$ for
$k\leq l$, and $r_k,1-r_k\geq0$,
\[
  \operatorname{Tr}(AQ)-\operatorname{Tr}(AP)
  =\sum_{k>l}\alpha_kr_k-\sum_{k\leq l}\alpha_k(1-r_k)
  \;\geq\;\alpha_{l+1}s-\alpha_l s
  \;=\;\Gamma s.
\]

\emph{Upper bound.} Let $S\preceq T$ be symmetric and $R$ an orthogonal
projection of rank $l$. Evaluating both traces by \eqref{eq:trace-basis} in
an orthonormal basis whose first $l$ vectors $v_1,\ldots,v_l$ span
$\operatorname{Range}(R)$ gives
\[
  \operatorname{Tr}(SR)=\sum_{i=1}^l\langle Sv_i,v_i\rangle_n
  \leq\sum_{i=1}^l\langle Tv_i,v_i\rangle_n=\operatorname{Tr}(TR).
\]
By the Ky Fan minimum principle
\cite[Corollary~4.3.39, p.~248]{horn_johnson},
$\operatorname{Tr}(MQ)$ is the minimum of $\operatorname{Tr}(MR)$
over rank-$l$ orthogonal projections $R$, so it is at most
$\operatorname{Tr}(MP)$. Combining this with the inequality above for
$R=Q$ and for $R=P$ gives
\[
  \operatorname{Tr}(AQ)
  \;\overset{A\preceq M}{\leq}\;\operatorname{Tr}(MQ)
  \;\overset{\text{Ky Fan}}{\leq}\;\operatorname{Tr}(MP)
  \;\overset{M\preceq B}{\leq}\;\operatorname{Tr}(BP).
\]
Subtracting $\operatorname{Tr}(AP)$ and evaluating the resulting trace
in the basis $\{u^k\}$,
\[
  \operatorname{Tr}(AQ)-\operatorname{Tr}(AP)
  \;\leq\;\operatorname{Tr}\!\bigl((B-A)P\bigr)
  \;=\;\sum_{k=1}^l\langle (B-A)u^k,u^k\rangle_n
  \;\leq\; l\tau.
\]

The two bounds give $\Gamma s\leq l\tau$, and by \eqref{eq:hs-identity},
$\|Q-P\|_{\mathrm{HS},n}^2=2s\leq 2l\tau/\Gamma$.
\end{proof}

\subsection[Proof of the projector theorem]{Proof of Theorem~\ref{thm:projector-headline}}\label{subsec:proof-thm-projector}
Proposition~\ref{prop:aux-projector} applies Lemma~\ref{lem:fcpb} to the
comparison operators, with $M=L$ and with $M=\widetilde L$.

\begin{proposition}\label{prop:aux-projector}
For every point cloud with $\Gamma_l^+>0$ and every choice of $P_l$ and
$\widetilde P_l$,
\begin{equation}\label{eq:aux-projector-bound}
  \max\bigl\{\|P_l-P_l^+\|_{\mathrm{HS},n},\,
             \|\widetilde P_l-P_l^+\|_{\mathrm{HS},n}\bigr\}
  \;\leq\; \sqrt{\frac{2l\tau_l}{\Gamma_l^+}}.
\end{equation}
\end{proposition}

\begin{proof}
By Lemma~\ref{lem:operator-bracket}, $A\preceq L\preceq B$ and
$A\preceq\widetilde L\preceq B$. In particular $A\preceq B$, so
$\tau_l\geq0$, and by \eqref{eq:tau-def},
$\langle(B-A)u,u\rangle_n\leq\tau_l\|u\|_n^2$ for every
$u\in\operatorname{Range}(P_l^+)$. Since $\Gamma_l^+>0$,
Lemma~\ref{lem:fcpb} applies to the comparison operators \eqref{eq:AB-def}
with $P=P_l^+$, $\Gamma=\Gamma_l^+$ and $\tau=\tau_l$. With $M=L$ and
$Q=P_l$, and with $M=\widetilde L$ and $Q=\widetilde P_l$,
Lemma~\ref{lem:fcpb} gives, for every choice of $P_l$ and $\widetilde P_l$,
\begin{equation}\label{eq:aux-projector-squared}
  \|P_l-P_l^+\|_{\mathrm{HS},n}^2\;\leq\;\frac{2l\tau_l}{\Gamma_l^+},
  \qquad
  \|\widetilde P_l-P_l^+\|_{\mathrm{HS},n}^2\;\leq\;\frac{2l\tau_l}{\Gamma_l^+}.
\end{equation}
Taking square roots in \eqref{eq:aux-projector-squared} proves
\eqref{eq:aux-projector-bound}.
\end{proof}

Next, Lemma~\ref{lem:gap-transfer} gives a lower bound on $\Gamma_l^+$. The
proof of Lemma~\ref{lem:gap-transfer} uses the results of
Section~\ref{sec:eigenvalue} at the two indices $l$ and $l+1$, so it needs the
event of Definition~\ref{def:hpe} for both indices at once. Since the event
$\mathcal E_n$ of Definition~\ref{def:hpe} depends on the index $l$, we write
$\mathcal E_n^{(k)}$ for the event of Definition~\ref{def:hpe} with $l$
replaced by $k$, so that $\mathcal E_n=\mathcal E_n^{(l)}$. The event
$\mathcal E_n^{(k)}$ is defined once the bandwidth assumptions made before
Definition~\ref{def:hpe} hold with $l$ replaced by $k$, which under
\eqref{eq:scaling} is the case for $n$ large enough. Set
\begin{equation}\label{eq:projector-event}
  \mathcal E_n^{\mathrm{proj}}
  \;:=\; \mathcal E_n^{(l)}\cap\mathcal E_n^{(l+1)}.
\end{equation}
Given $\mathcal E_n^{\mathrm{proj}}$, the results of
Section~\ref{sec:eigenvalue} hold at both indices.

\begin{lemma}\label{lem:gap-transfer}
Assume that $\gamma_l^{\mathcal M}>0$ and
\begin{equation}\label{eq:gap-smallness}
  (C_l'+C_{l+1}')\varepsilon
  +\tfrac12 K_{l+1}\,a
  \;\leq\; \tfrac12\sigma_\eta\gamma_l^{\mathcal M}.
\end{equation}
Then, given $\mathcal E_n^{\mathrm{proj}}$, the spectral gap $\Gamma_l^+$ of
$A$ satisfies
\[
  \Gamma_l^+
  \;\geq\; \tfrac12\sigma_\eta\gamma_l^{\mathcal M}\;>\;0.
\]
\end{lemma}

\begin{proof}
The lower bound of Corollary~\ref{cor:eig-bracket} reads
$\alpha_l\leq\lambda_l$. The upper bound of Corollary~\ref{cor:eig-bracket} at
index $l+1$, followed by Proposition~\ref{prop:aux-eig} at index $l+1$ (which
holds given $\mathcal E_n^{(l+1)}$), gives
\[
  \lambda_{l+1}
  \;\leq\; (1+a)^{m+2}\lambda_{l+1}^-
  \;\leq\; \alpha_{l+1}+\tfrac12K_{l+1}\,a.
\]
Moreover, \eqref{eq:cgt-cited} at index $l+1$, given $\mathcal E_n^{(l+1)}$,
and at index $l$, given $\mathcal E_n^{(l)}$, gives
$\lambda_{l+1}-\lambda_l\geq\sigma_\eta\gamma_l^{\mathcal M}
 -{(C_l'+C_{l+1}')\varepsilon}$. Combining these bounds with
\eqref{eq:gap-smallness},
\[
  \Gamma_l^+=\alpha_{l+1}-\alpha_l
  \;\geq\;\lambda_{l+1}-\lambda_l-\tfrac12K_{l+1}\,a
  \;\geq\;\sigma_\eta\gamma_l^{\mathcal M}
          -(C_l'+C_{l+1}')\varepsilon-\tfrac12K_{l+1}\,a
  \;\geq\;\tfrac12\sigma_\eta\gamma_l^{\mathcal M}. \qedhere
\]
\end{proof}

\begin{proof}[Proof of Theorem~\ref{thm:projector-headline}]
Since $\gamma_l^{\mathcal M}>0$, condition~\eqref{eq:gap-smallness} holds
for $n$ large enough and $\delta/\varepsilon$ small enough.
Lemma~\ref{lem:gap-transfer} then gives
$\Gamma_l^+\geq\tfrac12\sigma_\eta\gamma_l^{\mathcal M}>0$ given
$\mathcal E_n^{\mathrm{proj}}$. The event $\mathcal E_n^{\mathrm{proj}}$ is
contained in $\mathcal E_n$, so given $\mathcal E_n^{\mathrm{proj}}$,
Proposition~\ref{prop:aux-projector}, the triangle
inequality, the bound $\tau_l\leq\tfrac12K_l\,a$ of
Lemma~\ref{lem:operator-bracket} and $a=2\delta/\varepsilon$ yield
\begin{equation}\label{eq:projector-pair-bound}
\begin{aligned}
  \|P_l-\widetilde P_l\|_{\mathrm{HS},n}
  &\leq \|P_l-P_l^+\|_{\mathrm{HS},n}
         +\|\widetilde P_l-P_l^+\|_{\mathrm{HS},n} \\
  &\leq 2\sqrt{\frac{2l\tau_l}{\Gamma_l^+}}
   \leq 2\sqrt{\frac{2K_l\,l\,a}{\sigma_\eta\gamma_l^{\mathcal M}}}
   =D_l\sqrt{\frac{\delta}{\varepsilon}}.
\end{aligned}
\end{equation}
Proposition~\ref{prop:hpe} at the indices $l$ and $l+1$ and a union bound give
$\mathbb P_n(\mathcal E_n^{\mathrm{proj}})\geq1-O(n^{-p})$ once $C_\varepsilon$
is large enough. This proves the theorem.
\end{proof}

Since $\|T\|_{\mathrm{op},n}\leq\|T\|_{\mathrm{HS},n}$, the bound of
Theorem~\ref{thm:projector-headline} also holds in the operator norm.

\subsection[Proof of the eigenvector theorem]{Proof of Theorem~\ref{thm:eigenvector-headline}}\label{subsec:proof-thm-eigenvector}
Throughout Section~\ref{subsec:proof-thm-eigenvector}, $l\geq2$ and $p\geq1$
are as in
Theorem~\ref{thm:eigenvector-headline}, and $l_-,\dots,l_+$ is the
multiplicity block of $\lambda:=\lambda_l^{\mathcal M}$. Since
$\lambda_1^{\mathcal M}=0$ is simple and $l\geq2$, $\lambda>0$. For $L$ and
$\widetilde L$, choose orthonormal eigenbases ordered by increasing
eigenvalue and let $P_k$ and $\widetilde P_k$ be the orthogonal projections
onto the spans of their first $k$ vectors, so that
$\operatorname{Range}(P_j)\subseteq\operatorname{Range}(P_k)$ and
$\operatorname{Range}(\widetilde P_j)\subseteq\operatorname{Range}(\widetilde P_k)$
for $j\leq k$. Define
\[
  P_{l_-:l_+}:=P_{l_+}-P_{l_--1},
  \qquad
  \widetilde P_{l_-:l_+}:=\widetilde P_{l_+}-\widetilde P_{l_--1}.
\]
These are the orthogonal projections onto the subspaces spanned by the basis
vectors with indices $l_-,\dots,l_+$. We define $S_n:=\operatorname{Range}(P_{l_-:l_+})$,
the span of the eigenvectors of $L$ with these indices;
$E_\lambda:=\ker(\Delta_\rho-\lambda I)$, the eigenspace of $\Delta_\rho$
for $\lambda$, of dimension $r:=l_+-l_-+1$; and
${J_l:=\{l_--1,\,l_-,\,l_+,\,l_++1\}}$, the set of indices at which we use the
events $\mathcal E_n^{(j)}$.

Corollary~\ref{cor:cluster-vector} shows that every unit vector in
$\operatorname{Range}(\widetilde P_{l_-:l_+})$ is close to $S_n$, and
Lemma~\ref{lem:cgt-cluster-correspondence} shows that every unit vector in
$S_n$ is close to a unit eigenfunction in $E_\lambda$. The statement of
\cite[Theorem~2.7(2)]{calder_gt} gives the second property only for
eigenvectors of $L$. The proof of Theorem~\ref{thm:eigenvector-headline}
needs this property for a vector $v\in S_n$ that is a linear combination of
eigenvectors of $L$ for $\lambda_{l_-},\dots,\lambda_{l_+}$. These
eigenvalues are close to each other but in general distinct, so $v$ need not
be an eigenvector of $L$. Lemma~\ref{lem:cgt-cluster-correspondence}
therefore uses the proof of \cite[Theorem~2.7(2)]{calder_gt} instead of its statement, and
assumes the hypotheses of that proof, (A1)--(A4) and (S2$_N$) of
Appendix~\ref{app:cgt-hypotheses} with $N=l_++1$.

\begin{corollary}\label{cor:cluster-vector}
Suppose that \eqref{eq:gap-smallness} holds at the indices $l_--1$ and
$l_+$. Then, given $\bigcap_{j\in J_l}\mathcal E_n^{(j)}$,
\begin{equation}\label{eq:cluster}
  \|P_{l_-:l_+}-\widetilde P_{l_-:l_+}\|_{\mathrm{HS},n}
  \;\leq\; C_{l_-,l_+}\sqrt{\frac{\delta}{\varepsilon}},
  \qquad C_{l_-,l_+}:=D_{l_+}+D_{l_--1}.
\end{equation}
Consequently every unit
$\widetilde u\in\operatorname{Range}(\widetilde P_{l_-:l_+})$ satisfies
$\operatorname{dist}_n(\widetilde u,S_n)
 \leq C_{l_-,l_+}\sqrt{\delta/\varepsilon}$.
\end{corollary}

\begin{proof}
The spectral gaps $\gamma^{\mathcal M}_{l_--1}$ and $\gamma^{\mathcal M}_{l_+}$
are strictly positive, since $l_-$ and $l_+$
are the ends of the multiplicity block, and $\bigcap_{j\in J_l}\mathcal E_n^{(j)}$
is contained in
$\mathcal E_n^{(k)}\cap\mathcal E_n^{(k+1)}$ for $k=l_+$ and for $k=l_--1$.
Hence \eqref{eq:projector-pair-bound} holds with $l$ replaced by $l_+$ and
by $l_--1$, and \eqref{eq:cluster} follows from the triangle inequality
\[
  \|P_{l_-:l_+}-\widetilde P_{l_-:l_+}\|_{\mathrm{HS},n}
  \leq\|P_{l_+}-\widetilde P_{l_+}\|_{\mathrm{HS},n}
      +\|P_{l_--1}-\widetilde P_{l_--1}\|_{\mathrm{HS},n}.
\]
For the bound on $\operatorname{dist}_n(\widetilde u,S_n)$, the identity
$\widetilde P_{l_-:l_+}\widetilde u=\widetilde u$ gives
\[
\begin{aligned}
  \operatorname{dist}_n(\widetilde u,S_n)
  &= \|(I-P_{l_-:l_+})\widetilde u\|_n
   = \|(\widetilde P_{l_-:l_+}-P_{l_-:l_+})\widetilde u\|_n\\
  &\leq \|P_{l_-:l_+}-\widetilde P_{l_-:l_+}\|_{\mathrm{op},n}
   \leq \|P_{l_-:l_+}-\widetilde P_{l_-:l_+}\|_{\mathrm{HS},n}.
\end{aligned}
\qedhere
\]
\end{proof}

\begin{lemma}
\label{lem:cgt-cluster-correspondence}
Let $N:=l_++1$, and assume that $\varepsilon$, $\theta$ and $\delta_e$ satisfy
(A1)--(A4) and (S2$_N$) of Appendix~\ref{app:cgt-hypotheses}.
Then there exist an event $\mathcal V_{n,\lambda}$ and constants
$C_\lambda,c_\lambda>0$, depending only on $\lambda$, $\mathcal M$, $\rho$
and $\eta$, such that,
given $\mathcal V_{n,\lambda}$, for every $v\in S_n$ with $\|v\|_n=1$ there is an
$f\in E_\lambda$ with
\begin{equation}\label{eq:cgt-cluster-vector}
  \|f\|_{L^2(\mu)}=1,
  \qquad
  \|v-f|_{X_n}\|_n
  \;\leq\;\sqrt r\,C_\lambda\varepsilon.
\end{equation}
The failure probability of $\mathcal V_{n,\lambda}$ is at most
$C_\lambda n\exp(-c_\lambda n\theta^2\delta_e^m)
 +C_\lambda n\exp(-c_\lambda n\varepsilon^{m+4})$, as in
\cite[Theorem~2.7(2)]{calder_gt}.
\end{lemma}

\begin{proof}
We use the proof of \cite[Theorem~2.7(2)]{calder_gt}, which treats an
eigenvalue $\lambda>0$ of $\Delta_\rho$, as here. In its notation the block
of $\lambda$ has the indices $i+1,\dots,i+k$, so that $i=l_--1$ and $k=r$.
That proof assumes (S2$_N$) for some $N\geq i+k$ and uses the resulting
bound \cite[(4.13)]{calder_gt} also at the index $i+k+1=l_++1$, which is why we take
$N=l_++1$. Under (A1)--(A4) and (S2$_N$), \cite[Theorem~2.4(1)]{calder_gt}
gives, with probability at least $1-Cn\exp(-cn\theta^2\delta_e^m)$,
\[
  |\lambda_j-\sigma_\eta\lambda_j^{\mathcal M}|
  \;\leq\;\tfrac12\sigma_\eta g_j,
  \qquad j=1,\dots,N.
\]
This is (4.13) of that proof, with the factor $\sigma_\eta$ discussed in
Appendix~\ref{app:cgt-hypotheses}. Since $g_{l_--1}\leq\gamma^{\mathcal M}_{l_--1}$
and $g_{l_++1}\leq\gamma^{\mathcal M}_{l_+}$, (4.13) places $\lambda_{l_--1}$ and
$\lambda_{l_++1}$, and hence every eigenvalue of $L$ whose index lies
outside the block, at distance at least $\tfrac12\sigma_\eta g_l$ from
$\sigma_\eta\lambda$. Here $g_l$, defined in Appendix~\ref{app:cgt-hypotheses},
equals $\min\{\gamma^{\mathcal M}_{l_--1},\gamma^{\mathcal M}_{l_+}\}$ because
$l_-\leq l\leq l_+$.

Let $f_1,\ldots,f_r$ be an $L^2(\mu)$-orthonormal basis of $E_\lambda$.
Using this separation, \cite[(4.14)]{calder_gt} gives, for each $j$,
\[
  \operatorname{dist}_n\bigl(f_j|_{X_n},S_n\bigr)
  \;\leq\;\frac{2}{\sigma_\eta g_l}\,
  \bigl\|Lf_j|_{X_n}-\sigma_\eta(\Delta_\rho f_j)|_{X_n}\bigr\|_n.
\]
Pointwise consistency \cite[Theorem~3.3]{calder_gt}, the approximation of
$\sigma_\eta\Delta_\rho f$ by $Lf|_{X_n}$ described in
Section~\ref{sec:setup}, applies to the $f_j$, which are $C^3$
\cite[Section~2.3]{calder_gt}. It bounds the norms on the right-hand side by a
constant multiple of $\varepsilon$, with probability at least
$1-2rn\exp(-cn\varepsilon^{m+4})$. The proof of
\cite[Theorem~2.7(2)]{calder_gt} then combines this bound with (4.14) and
concludes that, with probability at least
$1-Cn\exp(-cn\theta^2\delta_e^m)-2rn\exp(-cn\varepsilon^{m+4})$, there is an
$\langle\cdot,\cdot\rangle_n$-orthonormal basis $v_1,\ldots,v_r$ of $S_n$
with
\begin{equation}\label{eq:cgt-basis-correspondence}
  \|v_j-f_j|_{X_n}\|_n
  \;\leq\;C_\lambda\varepsilon,
  \qquad j=1,\ldots,r.
\end{equation}
We take $\mathcal V_{n,\lambda}$ to be the event that such a basis exists.
This event is defined without $\theta$ and $\delta_e$, and the probability
bound just quoted gives the failure probability in the statement. The last
sentence of that proof uses this basis to derive the statement of
\cite[Theorem~2.7(2)]{calder_gt}, which concerns eigenvectors of $L$. We use
the basis itself, which gives the bound for every unit vector of $S_n$.

Now take any unit $v\in S_n$. There are coefficients
$\beta_1,\ldots,\beta_r\in\mathbb R$ such that
\[
  v=\sum_{j=1}^r \beta_jv_j,
  \qquad \sum_{j=1}^r \beta_j^2=1.
\]
Define $f:=\sum_{j=1}^r \beta_jf_j$. Since all $f_j$ belong to
$E_\lambda$ and are orthonormal in $L^2(\mu)$,
$f\in E_\lambda$ and $\|f\|_{L^2(\mu)}=1$. Moreover,
\[
  \|v-f|_{X_n}\|_n
  \leq \sum_{j=1}^r |\beta_j|\,
                    \|v_j-f_j|_{X_n}\|_n
  \leq C_\lambda\varepsilon
        \sum_{j=1}^r|\beta_j|
  \leq \sqrt r\,C_\lambda\varepsilon,
\]
where the last inequality is Cauchy--Schwarz. This proves
\eqref{eq:cgt-cluster-vector}.
\end{proof}

\begin{proof}[Proof of Theorem~\ref{thm:eigenvector-headline}]
Since the spectral gaps $\gamma^{\mathcal M}_{l_--1}$ and
$\gamma^{\mathcal M}_{l_+}$ are strictly positive,
\eqref{eq:gap-smallness} holds at these two indices for $n$ large enough and
$\delta/\varepsilon$ small enough. Let $N:=l_++1$, and take
$\theta:=\theta_N$ and $\delta_e:=\theta_N\varepsilon$ as in
Lemma~\ref{lem:cgt-parameters} of Appendix~\ref{app:cgt-hypotheses}. For $n$
large enough, (A1)--(A4) and
(S2$_N$) hold, which are the hypotheses of
Lemma~\ref{lem:cgt-cluster-correspondence}. With
$\mathcal V_{n,\lambda}$ the event of
Lemma~\ref{lem:cgt-cluster-correspondence}, set
\[
  \mathcal H_{n,l}
  :=\mathcal V_{n,\lambda}\cap\bigcap_{j\in J_l}\mathcal E_n^{(j)}.
\]

Given $\mathcal H_{n,l}$, let $\widetilde u_l$ be an eigenvector of
$\widetilde L$ with eigenvalue $\widetilde\lambda_l$ and
$\|\widetilde u_l\|_n=1$. We show that $\widetilde u_l$ is close to a unit
vector $v\in S_n$, and that $v$ is close to a unit eigenfunction of
$\Delta_\rho$ for $\lambda$. Finally we bound $\mathbb P_n(\mathcal H_{n,l})$.

The eigenbasis of $\widetilde L$ chosen at the start of
Section~\ref{subsec:proof-thm-eigenvector} was arbitrary. Since $\widetilde u_l$ is a unit eigenvector for the $l$-th
eigenvalue $\widetilde\lambda_l$, we may take it as the $l$-th vector of this
basis. Then $\widetilde u_l$ lies in
$\operatorname{Range}(\widetilde P_{l_-:l_+})$, because $l_-\leq l\leq l_+$.
The hypothesis of Corollary~\ref{cor:cluster-vector}, \eqref{eq:gap-smallness}
at the indices $l_--1$ and $l_+$, holds by the first paragraph of this proof,
and $\mathcal H_{n,l}$ is contained in the event
$\bigcap_{j\in J_l}\mathcal E_n^{(j)}$ of that corollary. The corollary gives
\[
  t:=\operatorname{dist}_n(\widetilde u_l,S_n)
  \;\leq\;C_{l_-,l_+}\sqrt{\frac{\delta}{\varepsilon}}.
\]
For $\delta/\varepsilon$ small enough, $t<1$, so
$z:=P_{l_-:l_+}\widetilde u_l$ is nonzero. Set $v:=z/\|z\|_n$, a unit
vector in $S_n$. The orthogonal decomposition
$\widetilde u_l=z+(I-P_{l_-:l_+})\widetilde u_l$ gives
$\|z\|_n^2=1-t^2$ and $\langle\widetilde u_l,v\rangle_n=\|z\|_n$. Since
$\|z\|_n\leq1$, we have $\|z\|_n\geq\|z\|_n^2$, so
\[
  \|\widetilde u_l-v\|_n^2
  =2-2\langle\widetilde u_l,v\rangle_n
  =2-2\|z\|_n
  \leq2-2\|z\|_n^2
  =2t^2,
\]
and therefore
\begin{equation}\label{eq:normalized-clean-cluster-vector}
  \|\widetilde u_l-v\|_n
  \;\leq\;\sqrt2\,C_{l_-,l_+}
             \sqrt{\frac{\delta}{\varepsilon}}.
\end{equation}

By Lemma~\ref{lem:cgt-cluster-correspondence} there is an
$f_l\in E_\lambda$, that is, with
$\Delta_\rho f_l=\lambda_l^{\mathcal M}f_l$, such that
$\|f_l\|_{L^2(\mu)}=1$ and
\[
  \|v-f_l|_{X_n}\|_n
  \;\leq\;\sqrt r\,C_\lambda\varepsilon.
\]
Combining this bound with \eqref{eq:normalized-clean-cluster-vector} by the
triangle inequality gives
\begin{equation}\label{eq:noisy-continuum-eigenvector}
  \|\widetilde u_l-f_l|_{X_n}\|_n
  \;\leq\;
  \sqrt2\,C_{l_-,l_+}\sqrt{\frac{\delta}{\varepsilon}}
  +\sqrt r\,C_\lambda\varepsilon.
\end{equation}
The constants $C_{l_-,l_+}$, $C_\lambda$ and $r$ do not depend on $n$ or on
the displacements, so \eqref{eq:noisy-continuum-eigenvector} is the asserted
$O\bigl(\sqrt{\delta/\varepsilon}+\varepsilon\bigr)$ estimate. If
$\delta=O(\varepsilon^3)$, then $\sqrt{\delta/\varepsilon}=O(\varepsilon)$
and the right-hand side is $O(\varepsilon)$.

It remains to bound $\mathbb P_n(\mathcal H_{n,l})$. By
Lemma~\ref{lem:cgt-cluster-correspondence}, the failure probability of
$\mathcal V_{n,\lambda}$ is at most
\[
  C_\lambda n\exp\bigl(-c_\lambda n\theta^2\delta_e^m\bigr)
  +C_\lambda n\exp\bigl(-c_\lambda n\varepsilon^{m+4}\bigr).
\]
By Lemma~\ref{lem:cgt-parameters} with $c=c_\lambda$, the first term is at most
$C_\lambda n^{-p}$ for $n$ large, and the second is at most
$C_\lambda n^{-p}$ once $c_\lambda C_\varepsilon^{m+4}\geq p+1$. These two
bounds, Proposition~\ref{prop:hpe} at each index $j\in J_l$ and a union bound
give $\mathbb P_n(\mathcal H_{n,l})\geq1-O(n^{-p})$.
\end{proof}

\begin{proof}[Proof of Remark~\ref{rem:index-one}]
Integration by parts gives
\[
  \langle\Delta_\rho f,f\rangle_{L^2(\mu)}
  =\tfrac12\int_{\mathcal M}\rho^2|\nabla f|^2\,d\operatorname{Vol}_{\mathcal M}.
\]
Since $\mathcal M$ is connected, the eigenfunctions for
$\lambda_1^{\mathcal M}=0$ are the constant functions, so
$\lambda_1^{\mathcal M}$ is simple, and since $\mu$ is a probability measure,
its unit eigenfunctions are the constants $\pm1$. The noisy graph Laplacian
has $\widetilde\lambda_1=0$, with the constant vectors as eigenvectors. Given
$\mathcal E_n^{(2)}$, which has probability at least $1-O(n^{-p})$ by
Proposition~\ref{prop:hpe}, the proof of Corollary~\ref{cor:continuum-eig} at
index $2$ gives
$|\widetilde\lambda_2-\sigma_\eta\lambda_2^{\mathcal M}|
 \leq K_2\,\delta/\varepsilon+C_2'\varepsilon$. Since
$\lambda_2^{\mathcal M}>0$, this gives $\widetilde\lambda_2>0$ for $n$ large
enough and $\delta/\varepsilon$ small enough, so $\widetilde\lambda_1$ is
simple. Every
eigenvector $\widetilde u_1$ for $\widetilde\lambda_1$ with
$\|\widetilde u_1\|_n=1$ is therefore a constant vector $(c,\dots,c)$ with
$c^2=\|\widetilde u_1\|_n^2=1$. The constant function $f_1\equiv c$ is then a
unit eigenfunction of $\Delta_\rho$ for $\lambda_1^{\mathcal M}$ whose
restriction is
$f_1|_{X_n}=\widetilde u_1$, so
\[
  \|\widetilde u_1-f_1|_{X_n}\|_n=0. \qedhere
\]
\end{proof}

\section{Numerical Experiments}\label{sec:numerics}

This section reports numerical experiments on the bounds from which
Theorems~\ref{thm:main} and~\ref{thm:projector-headline} are proved. Let
\[
  W_l(a) \;:=\; (1+a)^{m+2}\lambda^-_l-(1-a)^{m+2}\lambda^+_l,
\]
the difference between the $l$-th eigenvalues of the comparison operators $B$ and $A$ of
\eqref{eq:AB-def}. Recall from \eqref{eq:gap-defs} that $\Gamma^+_l$ is the difference
between the $(l+1)$-st and the $l$-th eigenvalue of $A$, and from \eqref{eq:tau-def} that
$\tau_l$ is the largest value of $\langle(B-A)u,u\rangle_n/\|u\|_n^2$ over nonzero $u$ in
the range of $P^+_l$, the span of the first $l$ eigenvectors of $L^+$. The proofs of both
theorems pass through the chains
\begin{equation}\label{eq:num-chain}
  |\lambda_l-\widetilde\lambda_l|\;\le\;\underbrace{W_l(a)}_{\text{computable}}
    \;\le\;\tfrac12K_l\,a,
  \qquad
  \|P_l-\widetilde P_l\|_{\mathrm{HS},n}
    \;\le\;\underbrace{2\sqrt{2l\tau_l/\Gamma^+_l}}_{\text{computable}}
    \;\le\;D_l\sqrt{\delta/\varepsilon}.
\end{equation}
We call $|\lambda_l-\widetilde\lambda_l|$ the \emph{eigenvalue error},
$\|P_l-\widetilde P_l\|_{\mathrm{HS},n}$ the \emph{projector distance}, and
$2\sqrt{2l\tau_l/\Gamma^+_l}$ the \emph{computable projector bound}. For a given point
cloud both middle terms can be computed from $L^+$ and $L^-$, and we call them the
\emph{computable bounds}.

The experiments measure the computable bounds and stop there. The constant $K_l$ does
not depend on $\varepsilon$ or $n$, but it contains the prefactor
$\mathcal C^{\mathrm{unif}}_{\mathrm{CL}}$ of Lemma~\ref{lem:uniform-cl}, which is not
available in effective form, so the right-hand sides of \eqref{eq:num-chain} cannot be
evaluated. The first inequality in each chain holds for every point cloud (for the
projectors, every point cloud with $\Gamma^+_l>0$). These first inequalities are
Corollary~\ref{cor:eig-bracket}, and Proposition~\ref{prop:aux-projector} with the
triangle inequality, so checking them tests the implementation. The second inequalities,
Proposition~\ref{prop:aux-eig} and \eqref{eq:projector-pair-bound}, hold for $a<a_\star$
(Remark~\ref{rem:standing}), given the event $\mathcal E_n$ for the eigenvalues and, under
\eqref{eq:gap-smallness}, given the event $\mathcal E_n^{\mathrm{proj}}$ of
\eqref{eq:projector-event} for the projectors.

\begin{figure}[htbp]
  \centering
  \includegraphics[width=\textwidth]{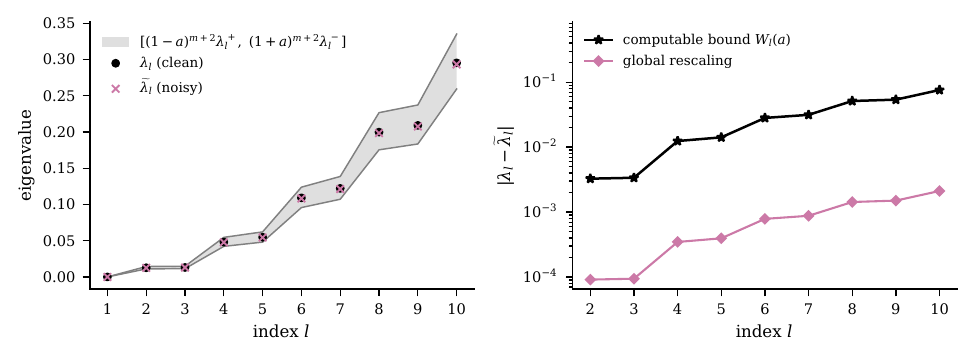}
  \caption{The bracket on one point cloud ($S^1$, $n=\bracketN$, $\varepsilon=\bracketEps$,
  $a=\bracketA$, global rescaling). Left: at
  every index, $\lambda_l$ and $\widetilde\lambda_l$ lie between
  $(1-a)^{m+2}\lambda^+_l$ and $(1+a)^{m+2}\lambda^-_l$. The eigenvalues come in near
  pairs because every nonzero eigenvalue of the circle has multiplicity two. Right: the
  computable bound $W_l(a)$ and the eigenvalue error it bounds. Index $1$ is omitted
  there, since $\lambda_1=\widetilde\lambda_1=0$ up to rounding.}
  \label{fig:bracket}
\end{figure}

\begin{figure}[htbp]
  \centering
  \includegraphics[width=0.49\textwidth]{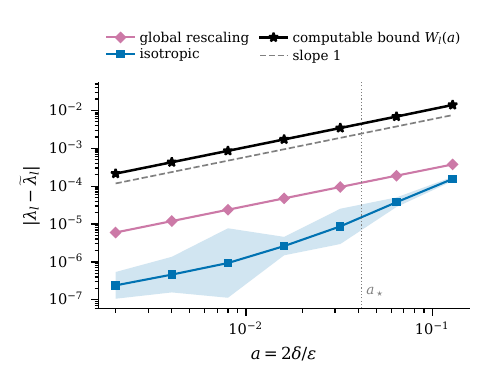}\hfill
  \includegraphics[width=0.49\textwidth]{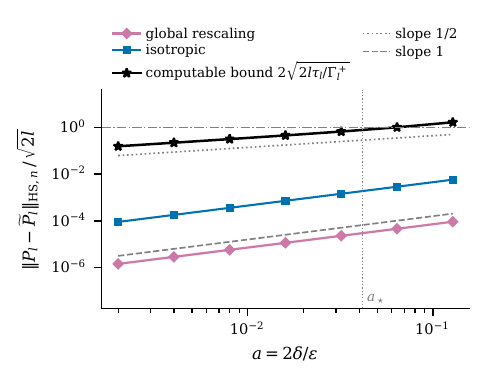}
  \caption{Dependence on $a=2\delta/\varepsilon$ ($\projManifold$, $n=\projN$,
  $\varepsilon=\projEps$, $l=\lFocus$,
  \seedsRates\ seeds). Left: the eigenvalue error under global rescaling and isotropic
  noise, and the computable bound $W_l(a)$. Right: the projector distance and the
  computable bound $2\sqrt{2l\tau_l/\Gamma^+_l}$, both divided by $\sqrt{2l}$, so the
  dash-dotted line at one is the trivial bound. The vertical line is at $a_\star$, and
  the sloping grey lines have the slopes given in the legends. Shaded bands show the
  scatter between seeds; where no band is visible it is narrower than the curve.}
  \label{fig:rates}
\end{figure}

\begin{figure}[!htbp]
  \centering
  \includegraphics[width=\textwidth]{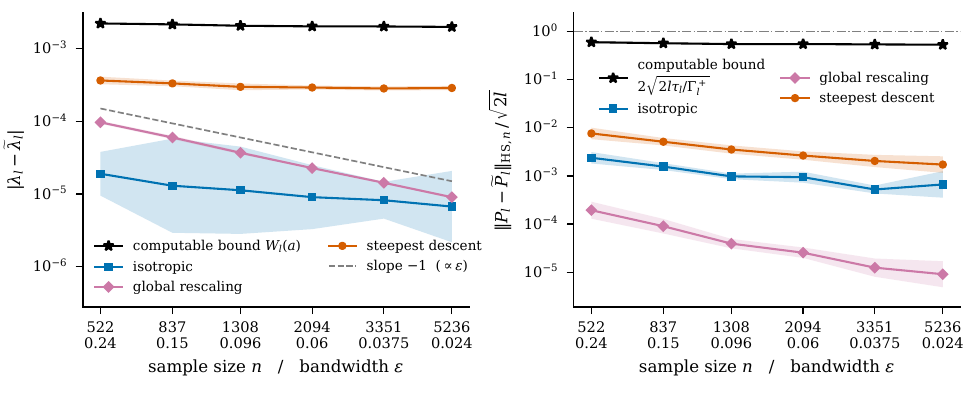}
  \caption{Fixed $a=\refineA$ and expected degree $\refineDegree$ as $\varepsilon$
  decreases and $n$ grows ($S^1$, $l=\refineL$,
  \refineSeeds\ seeds). The tick labels on the
  horizontal axis give $n$ and $\varepsilon$. Left: the eigenvalue errors under the three
  kinds of noise, and the computable bound $W_l(a)$. The dashed grey line has slope $-1$
  against $n$, that of a quantity proportional to $\varepsilon$. Right: the projector
  distance under the same noise and the computable
  bound $2\sqrt{2l\tau_l/\Gamma^+_l}$, both divided by $\sqrt{2l}$, so the dash-dotted
  line at one is the trivial bound.}
  \label{fig:fixed-a-refine}
\end{figure}

\begin{figure}[!htbp]
  \centering
  \includegraphics[width=\textwidth]{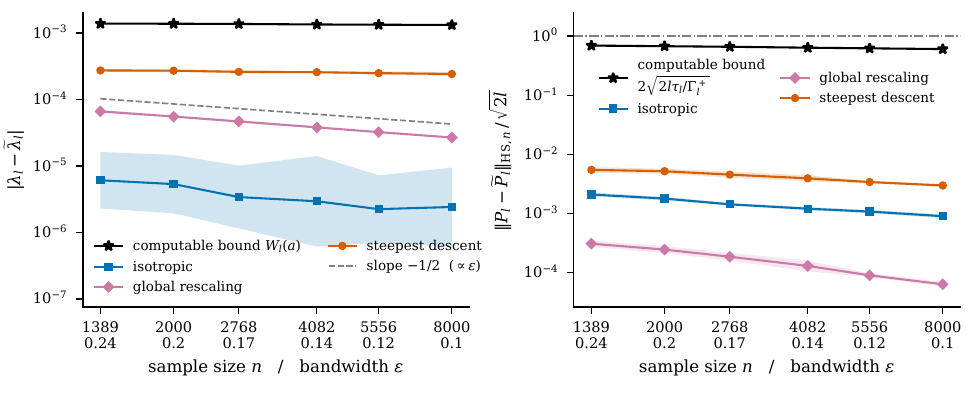}
  \caption{The experiment of Figure~\ref{fig:fixed-a-refine} on the unit sphere: fixed
  $a=\sTwoA$ and expected degree $\sTwoDegree$ as $\varepsilon$ decreases and $n$ grows
  ($S^2$, $l=\sTwoL$, \sTwoSeeds\ seeds).
  The panels, kinds of noise and units are those of Figure~\ref{fig:fixed-a-refine}. Since
  $n\varepsilon^2$ is constant along the sweep, a quantity proportional to $\varepsilon$
  has slope $-1/2$ against $n$, the slope of the dashed grey line.}
  \label{fig:fixed-a-sphere}
\end{figure}

There are three experiments, and the first inequalities of \eqref{eq:num-chain} hold up
to rounding error for every point cloud in all of them. The first experiment
(Figures~\ref{fig:bracket} and~\ref{fig:rates}) varies $a$ for fixed clean points and
bandwidth: $W_l(a)$ and $\tau_l$ grow linearly in $a$, and the computable projector
bound grows roughly like $\sqrt a$. Varying $a$ alone cannot show how the computable bounds
depend on $\varepsilon$, which needs $a$ fixed while $\varepsilon$ decreases. The other
two experiments (Figures~\ref{fig:fixed-a-refine} and~\ref{fig:fixed-a-sphere}) do this
on $S^1$ and $S^2$. Both computable bounds stay nearly constant while $\varepsilon$
decreases from $\refineEpsHi$ to $\refineEpsLo$ on $S^1$ and from $\sTwoEpsHi$ to
$\sTwoEpsLo$ on $S^2$. Under noise built from the clean points to make the eigenvalue
error large, this error is between $\refineWorstPctLo\%$ and $\refineWorstPctHi\%$ of
$W_l(a)$ on $S^1$ and between $\sTwoWorstPctLo\%$ and $\sTwoWorstPctHi\%$ on $S^2$, so
it stays of order $\delta/\varepsilon$ as $\varepsilon$ decreases. In these two
experiments $n$ grows so that every point keeps the same expected number of neighbours.
Then $\varepsilon$ is of order $n^{-1/m}$ and decreases much faster than the bandwidth
$\varepsilon=C_\varepsilon(\log n/n)^{1/(m+4)}$ of the theorems, along which a tenfold
decrease of $\varepsilon$ would need $n$ to grow more than $10^{m+4}$-fold. The faster
decrease is legitimate because the first inequalities of \eqref{eq:num-chain} hold at
every bandwidth, so the computable bounds still bound the eigenvalue error and the
projector distance.

\paragraph{Setup.} The experiments of Figures~\ref{fig:bracket}--\ref{fig:fixed-a-refine}
sample i.i.d.\ points from the uniform distribution on the circle $S^1\subset\mathbb R^2$,
and that of Figure~\ref{fig:fixed-a-sphere} from the uniform distribution on the unit
sphere $S^2\subset\mathbb R^3$. Both manifolds satisfy the hypotheses of
Section~\ref{sec:setup}. Other manifolds and non-constant densities are not tested. The
kernel profile is $\eta(t)=(1-t)_+$, normalized to unit mass on $B(0,1)\subset\mathbb
R^m$, which satisfies (K1)--(K5). The index $l$ is the end of a multiplicity block, so
that $\gamma^{\mathcal M}_l>0$ as Theorem~\ref{thm:projector-headline} requires. Noise
is applied in the ambient space, $\mathbb R^2$ or $\mathbb R^3$, and the noisy points are
not projected back to the manifold. The \emph{expected degree} is the expected number of
other sample points within distance $\varepsilon$ of a given sample point, where the
distance is that of the ambient space.

The theorems hold for any displacements $\delta_1,\dots,\delta_n$ with
$\|\delta_i\|\le\delta$, including noise that depends on the clean points
(Lemma~\ref{lem:denoise-uniform}). We use three kinds of
noise, each with $\max_i\|\delta_i\|=\delta$. The \emph{global rescaling}
$y_i=\bar x+(1+\delta/R)(x_i-\bar x)$, with $\bar x$ the centroid of the clean points and
$R=\max_i\|x_i-\bar x\|$, multiplies every pairwise distance by $1+\delta/R$.
\emph{Isotropic} noise moves each point by $\delta$ in an independent, uniformly
distributed direction. \emph{Steepest-descent} noise, used only in
Figures~\ref{fig:fixed-a-refine} and~\ref{fig:fixed-a-sphere} and defined before them,
is built from the clean points to make the eigenvalue error large.

Each experiment is run on several independent point clouds, one for each seed. Plotted
values are means over the seeds, and a shaded band shows the seed-to-seed scatter. Since
each plotted quantity is positive and each axis logarithmic, the band runs from
$\mu e^{-s}$ to $\mu e^{s}$ about the mean $\mu$, where $s$ is the standard deviation of
the logarithms of the seed values. The band displays the observed scatter and is not a
confidence interval. Every exponent quoted below is a least-squares fit in logarithmic
coordinates over the plotted grid, or over its part with $a<a_\star$ where stated. The
sloping grey lines in the figures are reference lines, with the slopes given in the
legends.

\paragraph{The bracket.} Figure~\ref{fig:bracket} shows
Corollary~\ref{cor:eig-bracket} on one point cloud, at the first ten indices. In the
right panel, $W_l(a)$ and the eigenvalue error both grow with $l$, and the error is
nearly the same fraction of $W_l(a)$ at every index. The bracket uses only $a<1$ and the
monotonicity of $\eta$, so it also holds at the values of $a$ above $a_\star$ in
Figure~\ref{fig:rates}.

\paragraph{Dependence on $a$.}
Figure~\ref{fig:rates} varies $a$ at $n=\projN$, $\varepsilon=\projEps$ and
$l=\lFocus$. Each point cloud is the same at every value of $a$, and so are the
directions of the isotropic noise, so only $\delta$ changes. In the left panel,
$W_l(a)$ is parallel to the reference line of slope one: its fitted exponent in $a$ is
$\expEigWall$. The observed eigenvalue error is at most about $\errWallMaxPct\%$ of
$W_l(a)$ under both kinds of noise. The theorems do not predict how that error depends
on $a$, and under isotropic noise it varies strongly between seeds.

The right panel shows the projectors. The computable projector bound is built from
$\tau_l$, which is not plotted. In this experiment the fitted exponent of $\tau_l$ in $a$
is $\expTau$. The computable projector bound has exponent $\expCert$ ($\expCertIn$
over $a<a_\star$) and stays close to the reference line of slope $1/2$, as the square
root in \eqref{eq:num-chain} suggests. It crosses the trivial bound
$\sqrt{2l}$ of Lemma~\ref{lem:hs-identity} between $a=\aCrossLo$ and $a=\aCrossHi$, so
at this sample size and bandwidth it is informative only for smaller $a$.

The observed projector distances have exponent $\expProjScaling$ under global rescaling
and $\expProjSphere$ under isotropic noise. They are measured at one bandwidth, whereas
Theorem~\ref{thm:projector-headline} bounds the distance uniformly in $n$ and
$\varepsilon$, so they do not show that its square root can be dropped.

\paragraph{Steepest-descent noise.} Let $u$ be an eigenvector of $L$ with eigenvalue
$\lambda_l$ and $\|u\|_n=1$, and set
\[
  \delta_i\;=\;\delta\,\frac{v_i}{\|v_i\|},
  \qquad
  v_i\;=\;\sum_{j\neq i:\,\|x_i-x_j\|<\varepsilon}(u_i-u_j)^2\,e_{ij},
  \qquad
  e_{ij}\;=\;\frac{x_i-x_j}{\|x_i-x_j\|},
\]
with $\delta_i=0$ if $v_i=0$. So every nonzero $\delta_i$ has length $\delta$ and points
in the direction of the sum of the vectors $e_{ij}$, weighted by $(u_i-u_j)^2$. Each
$e_{ij}$ points away from the neighbour $x_j$. Moving $x_i$ away from $x_j$ lowers the
weight $\omega_{ij}$, and so lowers $\langle Lu,u\rangle_n$, in which $\omega_{ij}$
appears with the factor $(u_i-u_j)^2$.

This noise is meant to make the eigenvalue error $|\lambda_l-\widetilde\lambda_l|$
large. Since $\delta$ is small compared with $\varepsilon$
($a=2\delta/\varepsilon$ is $\refineA$ in Figure~\ref{fig:fixed-a-refine} and $\sTwoA$ in
Figure~\ref{fig:fixed-a-sphere}), we look at how $\langle Lu,u\rangle_n$, which equals
$\lambda_l$, changes when every clean point starts to move in a chosen direction. Let
$h_1,\dots,h_n\in\mathbb R^d$ with $\|h_i\|\le1$, and let $L(s)$ be the graph Laplacian
of the points $x_i+sh_i$, so that $L(0)=L$. Write $\eta(t)=c\,(1-t)$ on $[0,1]$, with
$c>0$ the normalizing constant. If the $x_i$ are distinct and no $\|x_i-x_j\|$ equals
$\varepsilon$, which holds almost surely here, the chain rule and $e_{ji}=-e_{ij}$ give
\[
  \frac{d}{ds}\langle L(s)u,u\rangle_n\Big|_{s=0}
  \;=\;-\frac{c}{n^2\varepsilon^{m+3}}\sum_{i=1}^n\langle h_i,v_i\rangle.
\]
By the Cauchy--Schwarz inequality, this derivative is smallest for
$h_i=\delta_i/\delta$. So the noise moves the points in the directions in which
$\langle L(s)u,u\rangle_n$ decreases fastest at $s=0$, which gives the noise its name.

\paragraph{Fixed $a$ on $S^1$.}
Figure~\ref{fig:fixed-a-refine} holds $a=\refineA=\refineAFrac\,a_\star$ fixed while
$\varepsilon$ decreases from $\refineEpsHi$ to $\refineEpsLo$, so that
$\delta=a\varepsilon/2$ decreases from $\refineDeltaHi$ to $\refineDeltaLo$. The sample
size grows from $\refineNLo$ to $\refineNHi$, which keeps the expected degree at
$\refineDegree$. Then $\varepsilon$ is of order $n^{-1/m}$, which does not satisfy the
condition $\varepsilon\gg(\log n/n)^{1/m}$ under which Calder and Garc\'ia Trillos prove
that the eigenvalues of $L$ converge to those of $\sigma_\eta\Delta_\rho$
\cite[the discussion after Theorem~2.4]{calder_gt}. The probability bound of
Proposition~\ref{prop:hpe} does not apply either, so only the first inequalities of
\eqref{eq:num-chain} are assured.

In the left panel, $W_l(a)$ stays nearly constant: fitted against $n$, it scales as
$n^{\refineWallExp}$. This is consistent with Proposition~\ref{prop:aux-eig}, which bounds $W_l(a)$ by $\tfrac12K_l\,a$ given
$\mathcal E_n$, a bound that does not depend on $\varepsilon$. Under steepest-descent
noise the eigenvalue error, averaged over the seeds, is between $\refineWorstPctLo\%$ and
$\refineWorstPctHi\%$ of $W_l(a)$ at every bandwidth, so for this noise the error stays of
order $\delta/\varepsilon$ as $\varepsilon$ decreases. Under global rescaling and
isotropic noise the error is a decreasing fraction of $W_l(a)$. The global-rescaling
error scales as $n^{\refineScalingExp}$, parallel to the reference line of slope $-1$,
that is, in proportion to $\varepsilon$. So it is of order $\delta$ rather than
$\delta/\varepsilon$. This is expected: the rescaling changes every weight by at most
$C_\eta\delta/R$, a factor $\varepsilon/(2R)$ smaller than the bound $C_\eta a$ that holds
for any noise. The isotropic error also decreases, but it varies too much between seeds
to fit a power law.

Since every choice of displacements with $\|\delta_i\|\le\delta$ has eigenvalue error at
most $W_l(a)$, the largest error over all such choices lies between the error under
steepest-descent noise and $W_l(a)$; these experiments do not locate it. The expected
degree is the same at every bandwidth, so the experiment does not show that the error
under steepest-descent noise stays of order $\delta/\varepsilon$ along the bandwidth of
the theorems, where the expected degree grows with $n$.

The right panel of Figure~\ref{fig:fixed-a-refine} shows the projectors on the same
point clouds. The computable projector bound also stays nearly constant: it lies
between $\refineCertLo$ and $\refineCertHi$ times the trivial bound $\sqrt{2l}$ at every
bandwidth and, fitted against $n$, scales as $n^{\refineCertExp}$. This is consistent
with the bound $D_l\sqrt{\delta/\varepsilon}=D_l\sqrt{a/2}$ of
\eqref{eq:num-chain}, which holds given $\mathcal E_n^{\mathrm{proj}}$ and depends on $a$
alone. The observed distances decrease overall under all three kinds of noise. The
distance under steepest-descent noise is the largest of the three, and the computable
bound exceeds it by a factor that grows from $\refineProjGapLo$ to $\refineProjGapHi$.
Since steepest-descent noise is not chosen to move the projector, the experiment does
not show how large the projector distance can be at fixed $a$.

\paragraph{Fixed $a$ on $S^2$.}
Figure~\ref{fig:fixed-a-sphere} repeats this experiment on the unit sphere $S^2$, which,
unlike $S^1$, has dimension $m=2$ and is not intrinsically flat. Here $l=\sTwoL$, the end
of the first nonzero block, and $a=\sTwoA=\sTwoAFrac\,a_\star$, with $a_\star$ taken at
$m=2$. To hold the expected degree fixed, $n$ must now grow like $\varepsilon^{-2}$, so a
tenfold decrease of $\varepsilon$ would need a hundredfold increase of $n$, beyond the
reach of the dense eigensolvers used here. The range is therefore shorter: $\varepsilon$
decreases from $\sTwoEpsHi$ to $\sTwoEpsLo$ and $\delta$ from $\sTwoDeltaHi$ to
$\sTwoDeltaLo$, while $n$ grows from $\sTwoNLo$ to $\sTwoNHi$ at expected degree
$\sTwoDegree$. Over this shorter range, the fitted exponents are less informative than
on $S^1$.

The results agree with those on $S^1$. Both computable bounds stay nearly constant:
fitted against $n$, $W_l(a)$ scales as $n^{\sTwoWallExp}$, and the computable projector
bound lies between $\sTwoCertLo$ and $\sTwoCertHi$ times the trivial bound and scales as
$n^{\sTwoCertExp}$. Under steepest-descent noise the eigenvalue error, averaged over the
seeds, is between $\sTwoWorstPctLo\%$ and $\sTwoWorstPctHi\%$ of $W_l(a)$ at every
bandwidth. The global-rescaling error scales as $n^{\sTwoScalingExp}$; since
$\varepsilon$ is proportional to $n^{-1/2}$ here, it is again proportional to
$\varepsilon$. The isotropic error decreases overall, with a large seed-to-seed scatter.
The observed projector distances decrease. The distance under steepest-descent noise is
again the largest, and the computable bound exceeds it by a factor that grows from
$\sTwoProjGapLo$ to $\sTwoProjGapHi$.

\begin{remark}[Scope of these experiments]\label{rem:numerics-scope}
The experiments do not test the probability estimates for $\mathcal E_n$ and
$\mathcal E_n^{\mathrm{proj}}$, or the uniformity of the constants of
Lemma~\ref{lem:uniform-cl}, which Appendix~\ref{app:cl-audit} verifies, and no
experiment shows that a bound is sharp.
\end{remark}

\FloatBarrier
\appendix

\section{Uniformity of the Eigenvector Lipschitz Estimate over a Kernel Profile Family}
\label{app:cl-audit}

The constants in the eigenvector Lipschitz estimate of Calder, Garc\'ia
Trillos and Lewicka \cite[Theorem~2.3]{calder_lip} are allowed to depend on
the kernel profile. To apply that result to the family
$\{\bar\eta_a^\pm : a<a_\star\}$, this appendix verifies that the constants,
probability constants and smallness thresholds of that theorem can be chosen
uniformly over an admissible kernel profile family
(Definition~\ref{def:admissible}).

This appendix is not a new proof of that theorem.
It follows the argument for \cite[Theorem~2.3]{calder_lip} in the source's own
order and records, for each step, which quantities of the profile its constants
depend on and which smallness thresholds it imposes, including the implicit one
behind each ``$\varepsilon\ll1$''. The claims below can only be checked
against the source, and the appendix is meant to be read with
\cite{calder_lip} open. For $m=1$ two steps of the source need a small
change, made where they occur, in Sections~\ref{app:coupling}
and~\ref{app:contraction}.

Throughout, the manifold $\mathcal M$, its induced Riemannian metric $g$ and
the density $\rho$ are fixed. We denote by $\kappa$ a generic
member of an admissible family $\mathcal K$ (Definition~\ref{def:admissible}),
and use $A_{\varepsilon,\kappa}$, $\bar A_{\varepsilon,\kappa}$,
$\Delta_{\varepsilon,\kappa}$, $d_{\varepsilon,\kappa}$,
$d_{\varepsilon,\kappa,X_n}$, $\sigma_\kappa$ for the objects of \cite{calder_lip}
built with $\kappa$ in place of $\eta$. A result is called \emph{uniform} when
its constants and its implicit smallness thresholds can be chosen once for all
$\kappa\in\mathcal K$. Constants are named after their place in the source:
$C_{3.6}$ is the constant of \cite[Lemma~3.6]{calder_lip}, $c_{3.4(ii)}$ the
exponential constant of \cite[Proposition~3.4(ii)]{calder_lip}, and $C_{8.8}$ the
constant of display (8.8) in \cite{calder_lip}. An unnamed $C$ or $c$ is a
generic constant of the source, and what it depends on is stated where it
appears.

\subsection{Admissible families}

\begin{definition}[Admissible kernel profile family]\label{def:admissible}
Let $0<s_1<s_2<1$, and let
$S=B(0,s_2)\setminus\bar B(0,s_1)$ be the corresponding annulus. (The source
writes $\delta_1,\delta_2$ for $s_1,s_2$; we do not, because $\delta_i$ denotes
the displacements.)
A family $\mathcal K$ of profiles $\kappa:[0,\infty)\to[0,\infty)$ is
\emph{admissible} with data $(K_\infty, K_L, \sigma_*, S, c_S)$ if every
$\kappa\in\mathcal K$ satisfies:
\begin{enumerate}
  \item[(U1)] \emph{(structure)} $\kappa\ge0$, $\mathrm{supp}(\kappa)\subseteq[0,1]$,
        $\kappa$ is non-increasing, and $\int_{B(0,1)}\kappa(|w|)\,dw=1$;
  \item[(U2)] \emph{(size)} $\|\kappa\|_\infty\le K_\infty$;
  \item[(U3)] \emph{(regularity)} $\mathrm{Lip}_{[0,1]}(\kappa)\le K_L$;
  \item[(U4)] \emph{(second moment)} $\sigma_\kappa\ge\sigma_*>0$;
  \item[(U5)] \emph{(positivity)} $\kappa(|w|)\ge c_S>0$ for every
        $w\in B(0,s_2)$.
\end{enumerate}
\end{definition}

For every $\kappa\in\mathcal K$, condition (U1), together with (U3), gives
(CL1) of Section~\ref{subsec:aux-lipschitz}, and (U5) gives (CL2), which is
\cite[(9.1)]{calder_lip}, with the annulus $S$ and $c=c_S$. An annulus is the form of set for which
\cite[Proposition~9.2]{calder_lip} is stated. The
lower bound on the rest of the ball $B(0,s_2)$ is used in
Section~\ref{app:contraction}. The kernel $w\mapsto\kappa(|w|)$ is radial for
every profile, so radial symmetry is not listed as a condition.

\begin{lemma}\label{lem:our-family-admissible}
The family $\{\bar\eta_a^\pm : a\in[0,a_\star)\}$ is admissible with
\[
  K_\infty=2\eta(0),\quad K_L=3C_\eta,\quad \sigma_*=\tfrac13\sigma_\eta,\quad
  S=B(0,t_0/2)\setminus\bar B(0,t_0/4),\quad c_S=\tfrac{c_0}{3}.
\]
\end{lemma}

\begin{proof}
By Lemma~\ref{lem:structural}, $\eta_a^\pm$ is non-negative and
non-increasing, and so is $\bar\eta_a^\pm$, which is $\eta_a^\pm$ divided by
the positive constant $I_a^\pm$. The rest of (U1) and (U2)--(U5) are the
bounds of Lemma~\ref{lem:uniform-adm}, with $B(0,s_2)=B(0,t_0/2)$.
\end{proof}

\begin{proposition}\label{prop:cl-uniform}
Let $\mathcal K$ be admissible with data $(K_\infty, K_L, \sigma_*, S, c_S)$.
Then the constants, the probability constants and the smallness thresholds,
stated and implicit, of \cite[Theorem~2.3]{calder_lip} can be chosen to
depend only on
$\mathcal M$, $\rho$ and the data. In particular they can be chosen the same
for every $\kappa\in\mathcal K$.
\end{proposition}

Sections~\ref{app:concentration}--\ref{app:global} prove
Proposition~\ref{prop:cl-uniform}.
Lemma~\ref{lem:uniform-cl} is Proposition~\ref{prop:cl-uniform} applied to the
family of Lemma~\ref{lem:our-family-admissible}.

\subsection{Concentration and counting}\label{app:concentration}

\cite[Lemma~3.2]{calder_lip} is a Bernstein estimate for
$\sum_i\psi(x_i)\mathbf 1_{\{|x_i-x|\le\varepsilon\}}$, with $\psi$ a bounded
Borel function, valid for $\varepsilon^2\le t\le1$. Its constant depends only
on $\mathcal M$. The density enters through an explicit factor
$\|\rho\|_{C^0}$, and $\psi$ enters only through a factor $\|\psi\|_\infty$ in
the size of the deviation. In particular, the exponent of the probability
bound does not involve $\psi$. We use this below, where the lemma is applied
to functions $\psi$ built from $\kappa$. The requirement $\varepsilon\ll1$ of
Lemma~3.2 is geometric: it comes from the inclusions between Euclidean and
geodesic balls and from volume estimates for small balls.
\cite[Corollary~3.3]{calder_lip} applies the lemma with
$\psi\equiv1$ at one fixed value of $t$ depending only on $\mathcal M$ and
$\rho$, and \cite[Propositions~3.4 and~3.5]{calder_lip} involve no kernel. So
all four results are uniform.

\subsection{Degree control}\label{app:degree}

\cite[Lemma~3.6]{calder_lip} gives
\[
  \mathbb P\Bigl(\sup_{x\in\mathcal M}|d_{\varepsilon,\kappa,X_n}(x)-d_{\varepsilon,\kappa}(x)|\le C_{3.6}t\Bigr)
  \ \ge\ 1-C(t\varepsilon)^{-m}\exp(-ct^2n\varepsilon^m),
\]
with constants depending on $\mathcal M$, $\rho$, $\|\kappa\|_\infty$ and
$\mathrm{Lip}_{[0,1]}(\kappa)$, which (U2) and (U3) bound. The sup-norm enters
when Lemma~3.2 is applied to $\psi=\kappa(|\cdot-x|/\varepsilon)$. The
Lipschitz constant enters in \cite[(3.11)]{calder_lip}, which bounds the
difference between $\kappa(|y-x|/\varepsilon)$ and
$\kappa(d_{\mathcal M}(x,y)/\varepsilon)$ and also uses monotonicity. Both
enter again when the estimate is extended from a net to all $x\in\mathcal M$.
The requirement $\varepsilon\ll1$ of Lemma~3.6 comes from Lemma~3.2,
Corollary~3.3, Proposition~3.4 and the bound \cite[(3.2)]{calder_lip} between
Euclidean and geodesic distance, so its threshold depends only on $\mathcal M$
and $\rho$.

\cite[Corollary~3.7]{calder_lip} adds the expansion
$d_{\varepsilon,\kappa}=\rho+O(\varepsilon^2)$ of \cite[(3.9)]{calder_lip},
with a constant that does not depend on $\kappa$. Normalization gives the
leading term $\rho$, the first-order term vanishes by radial symmetry, and
the second-order term is at most
$C\varepsilon^2\int\kappa(|w|)\,dw=C\varepsilon^2$. Hence, given the event of
Lemma~3.6, $\sup_{\mathcal M}|d_{\varepsilon,\kappa,X_n}-\rho|\le C_{3.7}t$.
Since $C_{3.7}$ is uniform, the requirement
$t\ll1$ of the source has a uniform threshold, below which
$d_{\varepsilon,\kappa,X_n}\in[\tfrac12\rho_{\min},2\rho_{\max}]$.
\cite[Remark~3.8]{calder_lip}, which bounds the number of sample points in a
ball and is used in \cite[Theorem~4.2]{calder_lip}, concerns only the
indicator profile and does not involve $\kappa$.

\subsection{Concentration of the double convolution}\label{app:double-kernel}

\cite[Lemma~4.1]{calder_lip} bounds
$\sup_{x,y}|g_{n,\varepsilon}(x,y)-g_\varepsilon(x,y)|$ by $C_{4.1}t\varepsilon^{-m}$,
where $g_{n,\varepsilon}$ and $g_\varepsilon$ are the empirical and the
continuum double convolution of the kernel, defined in its proof.
With $\kappa_\varepsilon(s):=\varepsilon^{-m}\kappa(s/\varepsilon)$, applying
\cite[Lemma~3.2]{calder_lip} to
$\psi=\kappa_\varepsilon(|\cdot-x|)\kappa_\varepsilon(|\cdot-y|)/(n\rho)$ produces
the factor $\|\kappa\|_\infty^2$, and passing from a $t\varepsilon$-net to
arbitrary points produces $\mathrm{Lip}_{[0,1]}(\kappa)\|\kappa\|_\infty$ and
$\|\kappa\|_\infty^2$. These factors appear only in the constant $C_{4.1}$.
The requirement $\varepsilon\ll1$ of the lemma is that of Lemma~3.2,
Corollary~3.3 and Proposition~3.4, so its threshold depends only on
$\mathcal M$ and $\rho$.

The exponent in Lemma~3.2 does not involve $\psi$. Hence
$c_{4.1}=\min\{c_{3.2},c_{3.3},c_{3.4(ii)}\}$ does not depend on $\kappa$,
while (U2) and (U3) bound the remaining constants of Lemma~4.1.
Lemma~4.1 uses neither the normalization in (U1), nor (U4), nor (U5).

\subsection{Interpolation and the discrete-to-nonlocal transfer}\label{app:interpolation}

The estimate of \cite[Theorem~4.2]{calder_lip} holds given the intersection of
the three events of \cite[(4.5)]{calder_lip}: the event that the bound of
Remark~3.8 holds for the number of sample points in each ball of radius
$2\varepsilon$, the event of Lemma~3.6, and the event of Lemma~4.1. The
profile enters the constant of Theorem~4.2 through $C_{3.6}$ and $C_{4.1}$
and, in the proof itself, through $\|\kappa\|_\infty$,
$\mathrm{Lip}_{[0,1]}(\kappa)$ and the normalization. No other quantity of the
profile appears. Theorem~4.2 requires $\varepsilon\ll1$ for Remark~3.8,
Lemmas~3.6 and~4.1 and \cite[(3.2)]{calder_lip}, and $t\ll1$ to bound the
degree from below as in Section~\ref{app:degree}. Both thresholds are uniform.

\cite[Theorem~2.9]{calder_lip} is then Theorem~4.2 with $t=\varepsilon^2$,
combined with the two algebraic identities \cite[(4.10)]{calder_lip}. With
$I_{\varepsilon,X_n}$ the interpolation operator of the source, these express
$f-I_{\varepsilon,X_n}f$ and
$I_{\varepsilon,X_n}f-A_{\varepsilon,\kappa}(I_{\varepsilon,X_n}f)$ through the
graph Laplacian and the nonlocal Laplacian. The event of Theorem~2.9, which we
denote by $E_{2.9}$, and its probability constants are those of Theorem~4.2.
The identities add to the constant of Theorem~2.9 only the degree bounds of
Section~\ref{app:degree}. Substituting $t=\varepsilon^2$ in the probability
bound of Theorem~4.2 gives $(t\varepsilon)^{-2m}=\varepsilon^{-6m}$ and
$t^2n\varepsilon^m=n\varepsilon^{m+4}$, so
\[
  \mathbb P(E_{2.9})\ \ge\ 1-C^{\mathrm{prob}}_{2.9}\,\varepsilon^{-6m}
  \exp\bigl(-c^{\mathrm{prob}}_{2.9}\,n\varepsilon^{m+4}\bigr).
\]
Conditions (U1)--(U3) bound the constants of Theorem~2.9. The theorem uses
neither (U4) nor (U5).

\subsection{Geometry and averaging operators}\label{app:averaging}

\cite[Lemma~5.1]{calder_lip} and \cite[Proposition~5.2]{calder_lip} involve
no kernel, so their constants and thresholds are the same for every $\kappa$.
\cite[Lemma~6.3]{calder_lip} bounds the difference between
$A_{\varepsilon,\kappa}$ and $\bar A_{\varepsilon,\kappa}$,
\[
  A_{\varepsilon,\kappa}f(x)-\bar A_{\varepsilon,\kappa}f(x)
  =O(\varepsilon^2)\,\|f\|_{L^\infty(B_{\mathcal M}(x,\varepsilon))},
\]
with a constant depending only on $(\mathcal M,g)$ and $\rho$. Its proof uses
the normalization, radial symmetry, non-negativity and support of $\kappa$,
and the second moment $\sigma_\kappa$, which these bound by $1$. Consequently
\[
  \frac{\|\bar A_{\varepsilon,\kappa}f-f\|_\infty}{\varepsilon^2}
  \ \le\ \frac{\|A_{\varepsilon,\kappa}f-f\|_\infty}{\varepsilon^2}+C\|f\|_\infty,
\]
with one constant for every $\kappa\in\mathcal K$. Moreover, the identity
$\Delta_{\varepsilon,\kappa}f=-d_{\varepsilon,\kappa}(A_{\varepsilon,\kappa}f-f)/\varepsilon^2$
of \cite[Remark~6.2(ii)]{calder_lip} and the lower bound $d_{\varepsilon,\kappa}\ge\tfrac12\rho_{\min}$, which
follows from $d_{\varepsilon,\kappa}=\rho+O(\varepsilon^2)$ for $\varepsilon$
below a uniform threshold, give
$\|A_{\varepsilon,\kappa}f-f\|_\infty/\varepsilon^2
\le2\rho_{\min}^{-1}\|\Delta_{\varepsilon,\kappa}f\|_\infty$. The source uses
these two steps to deduce \cite[Theorem~2.11]{calder_lip} from
\cite[Theorem~10.1]{calder_lip}, which is stated with
$\bar A_{\varepsilon,\kappa}$. They add no dependence on the profile.

\subsection{The coupled random walk}\label{app:coupling}

\cite[Lemma~7.1]{calder_lip} is an identity obtained by a change of
variables, which uses only radial symmetry and the normalization.
\cite[Proposition~7.2]{calder_lip} identifies the one-step transition
operator of the walk, before the walk is stopped, with
$\bar A_{\varepsilon,\kappa}$. Since $\mathrm{supp}(\kappa)\subseteq[0,1]$,
every jump has length at most $\varepsilon$. The construction also
requires $\varepsilon\|\nabla\log\rho\|_\infty<\tfrac12$, which does not
involve $\kappa$.

For $m=1$ the field of isometries $Q$ of \cite[Section~7.1]{calder_lip} needs
a small change. The source builds it from a rotation $Q_0\in SO(m)$ with
$Q_0(-e_1)=e_1$, and $SO(1)$ contains none. Lemma~7.1 requires only a linear
isometry, so we take $Q_0=-1$. Then $Q(x)$ is multiplication by the sign of
$\langle\nabla\log\rho(x),e_1(x)\rangle$, and by $+1$ where this vanishes,
with $e_1(x)$ the parallel unit frame of the source. So
$|\nabla\log\rho(x)|\,Q(x)e_1(x)=\nabla\log\rho(x)$, which is
\cite[(7.1)]{calder_lip}, for either choice of $e_1(x)$, and in particular in
the frame of \cite[(8.2)]{calder_lip}.

Let $L_\rho$ be the Lipschitz constant of $\nabla\log\rho$, with tangent
vectors compared by parallel transport, and let $r$ be the radius of
\cite[Theorem~8.1]{calder_lip}. In the case of \cite[(8.11)]{calder_lip},
where the drifts at the positions of both walks exceed $\sigma_\kappa/4$ in
norm, the source bounds the distance between their unit directions, compared
by parallel transport, by $8L_\rho r/\sigma_\kappa$. This is less than $2$
when $r\le c\sigma_*$, a condition imposed below. For $m=1$ two unit vectors
at one point are equal or opposite, so the two directions are equal. Since
$Q(x)$ sends $t\,e_1(x)$ to $t$ times the unit direction of the drift at $x$,
for either choice of $e_1(x)$, the difference of the two isometries that
enters the drift term is zero. The other case
uses only that each $Q(x)$ is an isometry. Hence \cite[(8.12)]{calder_lip}
holds as in the source, with no new quantity of the profile.

Let $(X_j)_{j\ge0}$ be the walk started at $x_0$ (the source writes the time
index as $n$), and let $\tau$ be a stopping time that is at most the exit
time of the walk from $B_{\mathcal M}(x_0,r)$. \cite[Lemma~7.3]{calder_lip} gives
$\mathbb E[d_{\mathcal M}(X_\tau,x_0)^2]\le\tfrac32 m\sigma_\kappa\varepsilon^2\mathbb E[\tau]$.
The lemma holds when
$C_{7.3}(\varepsilon+r)\le\tfrac12 m\sigma_*$, which is the same condition
for every $\kappa\in\mathcal K$.

In the statement of \cite[Theorem~8.1]{calder_lip} the profile appears only
in the coefficient $4r/\sigma_\kappa$, which is at most $4r/\sigma_*$, and in
the operator $\bar A_{\varepsilon,\kappa}$. Besides that of Lemma~7.3, the
hypothesis $\varepsilon\ll r\ll1$ of Theorem~8.1 contains three requirements
that involve the profile.
The lower bound
$4\varepsilon^2\sigma_\kappa+O(\varepsilon^2r+\varepsilon^3)\ge3\varepsilon^2\sigma_\kappa$
of \cite[(8.8)]{calder_lip} holds when $C_{8.8}(r+\varepsilon)\le\sigma_*$.
In the case where both drifts exceed $\sigma_\kappa/4$ in norm,
\cite[(8.11)]{calder_lip} requires $8L_\rho r/\sigma_\kappa$ to be small,
which holds when $r\le c\sigma_*$. The supermartingale inequality before
\cite[(8.14)]{calder_lip} requires
$\tfrac12\sigma_\kappa+C(r+\varepsilon)\le\sigma_\kappa/r$, which holds when
$r\le1$ and $C(r+\varepsilon)\le\tfrac12\sigma_*$. These three requirements
and that of Lemma~7.3 hold when $r\le c\sigma_*$ and $\varepsilon\le c'r$,
with constants $c,c'$ depending only on $\mathcal M$ and $\rho$, as do
$C_{7.3}$, $C_{8.8}$ and $C$. The remaining requirements do not involve the
profile. We fix a radius $r_0\le c\sigma_*$ that satisfies them. Theorem~8.1 then holds with $r=r_0$ for every
$\kappa\in\mathcal K$, for $\varepsilon$ below one threshold.

\subsection{The contraction estimate}\label{app:contraction}

Condition (U5) is used only in this section.
\cite[Proposition~9.2]{calder_lip} is stated for an annulus on which the
kernel is bounded below by a positive constant. We apply it with the annulus
$S$ and the constant $c_S$ of Definition~\ref{def:admissible}. With
$H_N:=\{w:\langle w,e_1\rangle\in(\tfrac1{4N},\tfrac1{2N})\}$ and
$\alpha_\kappa:=\int_{H_N}\kappa(|w|)\,dw$, Proposition~9.2 asserts that there is an
integer $N>1$ such that (i) $\alpha_\kappa\in(0,1)$ and (ii)
$|S\cap\varphi_{\varepsilon,x_0,y_0}(S)|>\tfrac12|S|$ whenever
$d_{\mathcal M}(x_0,y_0)\le\varepsilon/N$ and $\varepsilon\ll1$. Here
$\varphi_{\varepsilon,x_0,y_0}$ is the map of \cite[(9.2)]{calder_lip}.

Property (ii) and its proof involve only the set $S$, the map
$\varphi_{\varepsilon,x_0,y_0}$ and the expansion \cite[(9.3)]{calder_lip},
whose constants depend only on $(\mathcal M,g)$. The kernel does not appear
in them. So (ii) holds for every $N$ above a threshold that depends only on
$S$, and its threshold for $\varepsilon$ depends only on $S$ and
$(\mathcal M,g)$. The same is true of the estimate
$|S\mathbin{\triangle}\varphi_\varepsilon\operatorname{Refl}\varphi_\varepsilon(S)|
=O(\varepsilon^2)$, which the proof of \cite[Theorem~9.1]{calder_lip} obtains
by the same argument ($\varphi_\varepsilon$ and $\operatorname{Refl}$ are as in
the source).

For property (i) the source takes $N>1/s_2$. Then $\tfrac1{2N}<s_2$, so the
open set $B(0,s_2)\cap H_N$ contains $t\,e_1$ for every
$t\in(\tfrac1{4N},\tfrac1{2N})$. Since $\kappa(|w|)\ge c_S$ on $B(0,s_2)$ by
(U5), for every $\kappa\in\mathcal K$
\begin{equation}\label{eq:alpha-star}
  \alpha_\kappa\ \ge\ \alpha_*\ :=\ c_S\bigl|B(0,s_2)\cap H_N\bigr|\ >\ 0.
\end{equation}
For $m=1$ the ball is needed here: $S=(-s_2,-s_1)\cup(s_1,s_2)$ and
$H_N=(\tfrac1{4N},\tfrac1{2N})$ are disjoint once $N\ge1/(2s_1)$. We fix $N$,
depending only on $S$, such that (i) and (ii) hold.

The proof of \cite[Theorem~9.1]{calder_lip} uses at most $L=4(3N-1)$ steps,
each of which shortens the distance between the two points by at least
$\varepsilon/(4N)$, and its final estimate holds with $\bar N=60N$. The
contraction factor is
$\theta_\kappa=1-\tfrac{c_S|S|}{2}\alpha_\kappa^{L}$, so by
\eqref{eq:alpha-star}
\[
  \theta_\kappa\ \le\ \theta_*\ :=\ 1-\tfrac{c_S|S|}{2}\alpha_*^{L}\ <\ 1,
\]
and the theorem holds with the contraction factor $\theta_*$ for every
$\kappa\in\mathcal K$. The proof bounds the coefficients of
\cite[(9.11)]{calder_lip} by $24N$ in its final estimate, so no bound on
$(1-\alpha_\kappa)^{-1}$ enters. The remaining constant $C_{9.1}$ and the
thresholds for $\varepsilon$ come from the geometric expansions, from
property (ii), and from requiring the $O(\varepsilon^2)$ error of each step
to be at most $\varepsilon/(4N)$. They depend only on $(\mathcal M,g)$, $S$,
$c_S$ and $N$. Theorem~9.1 involves no sample points.

\subsection{Global regularity and the eigenvector estimate}\label{app:global}

The proof of \cite[Theorem~10.1]{calder_lip} combines Theorems~8.1 and~9.1
into \cite[(10.1)]{calder_lip}, valid for $d_{\mathcal M}(x_0,y_0)<r$, where
$r$ is the radius of Theorem~8.1. We take $r=r_0$, the radius of
Section~\ref{app:coupling}. The finitely many smallness conditions of
Theorems~8.1 and~9.1, together with the requirement $5\varepsilon<r_0$, hold
for $\varepsilon$ below one uniform threshold. The proof takes the supremum
of (10.1) over pairs of points at distance at most $5\varepsilon$ and absorbs
the term with the factor $\theta_\kappa$ into the left-hand side. This uses
$1/(1-\theta_\kappa)\le1/(1-\theta_*)$, and the coefficient
$4r_0/\sigma_\kappa$ of Theorem~8.1 is at most $4r_0/\sigma_*$. The general
case follows by compactness of $\mathcal M$, as in the source, with constants
depending only on $\mathcal M$ and $r_0$. \cite[Theorem~2.11]{calder_lip} then follows by the two
steps of Section~\ref{app:averaging}. Theorems~10.1 and~2.11 involve no
sample points.

The estimate of \cite[Theorem~2.1]{calder_lip} holds given the intersection of
$E_{2.9}$ with the event of Lemma~3.6 at $t=\varepsilon^2$. By a union bound
and $(t\varepsilon)^{-m}=\varepsilon^{-3m}\le\varepsilon^{-6m}$, the probability of this intersection
satisfies the bound displayed for $E_{2.9}$ in
Section~\ref{app:interpolation}, with uniform constants
$C^{\mathrm{prob}}_{2.1}$ and $c^{\mathrm{prob}}_{2.1}$ in place of
$C^{\mathrm{prob}}_{2.9}$ and $c^{\mathrm{prob}}_{2.9}$. The Lipschitz
constant $C_{2.1}$ of the theorem is
a combination of the constants of Theorems~2.9 and~2.11 and of the lower
bound on the degree. Each of these is uniform by the preceding sections.

Let $\Lambda$ be the cutoff of Theorem~2.3. The estimate of
\cite[Theorem~2.3]{calder_lip} holds given the intersection of the event of
Theorem~2.1 with the event of \cite[(11.2)]{calder_lip} that every ball
$B_{\mathcal M}(x_i,r)$ contains at least $C'nr^m$ sample points. This $r$ is
not the radius $r_0$. The proof needs $r+\varepsilon\le1/(2C_{2.1}(\Lambda+1))$
and takes $r$ of order $1/(\Lambda+1)$. We take $r=c_r/(\Lambda+1)$ and
require $\varepsilon\le c_\varepsilon/(\Lambda+1)$, with $c_r$ and
$c_\varepsilon$ chosen so that this inequality holds and $r$ is below the
threshold of Corollary~3.3. Since $C_{2.1}$ is uniform, so are $c_r$ and
$c_\varepsilon$. The source applies Corollary~3.3 at radius $r/2$ around each
sample point, so the event of (11.2) fails with probability at most
$2n\exp(-c\,n(\Lambda+1)^{-m})$, where $c$ depends only on $\mathcal M$,
$\rho$ and $c_r$. Through $c_r$, the constant $c$ depends on the data of the
family, and it is the same for every $\kappa\in\mathcal K$.

Altogether, for every $\kappa\in\mathcal K$ and every $\varepsilon\le
c_\varepsilon/(\Lambda+1)$ below the uniform thresholds, the estimate of
\cite[Theorem~2.3]{calder_lip} holds with one prefactor and with probability
at least
\[
  1-C^{\mathrm{prob}}_{2.1}\,\varepsilon^{-6m}
     \exp\bigl(-c^{\mathrm{prob}}_{2.1}\,n\varepsilon^{m+4}\bigr)
   -2n\exp\bigl(-c\,n\,(\Lambda+1)^{-m}\bigr).
\]
This is the bound of Lemma~\ref{lem:uniform-cl}, with
$C_{\mathbb P}=C^{\mathrm{prob}}_{2.1}$ and $c_{\mathbb P}$ the smaller of
$c^{\mathrm{prob}}_{2.1}$ and $c$. Every constant and threshold above has been chosen in terms of $\mathcal M$,
$\rho$ and the data of $\mathcal K$ alone. This proves
Proposition~\ref{prop:cl-uniform}.
\section[Denoising with the Projection Algorithm of Aizenbud and Sober]{Denoising with the Projection Algorithm\\ of Aizenbud and Sober}
\label{app:denoising}

This appendix applies our results to the model of Aizenbud and Sober
\cite{aizenbud}, in which the observations are uniform in a tube of fixed
radius $\sigma$ around $\mathcal M$. We take the nearest points of the
observations on $\mathcal M$ as the clean points. The displacements are then bounded by
$\delta=\sigma$, which does not decrease with $n$, so $\delta/\varepsilon$ is
not small and our results do not apply to the observations directly. We
therefore first denoise the observations with the projection algorithm of
\cite{aizenbud}, which estimates the nearest point on $\mathcal M$ of a given
query point from a dataset drawn from the distribution of the observations,
and converges to it as the dataset grows (Theorem~\ref{thm:as-quoted}).
We apply it at
each observation, always with the same additional dataset, which is
independent of the observations. If, for each $n$, the additional dataset is
large enough, then with probability at least $1-O(n^{-p})$ the graph on the
denoised points has the rates that \cite{calder_gt} gives for the clean graph,
in the number $n$ of graph vertices.

We change some of the notation of \cite{aizenbud} to avoid clashes with
ours. Their ambient and intrinsic dimensions $D$ and $d$, number $n$ of data
points, query point $r$, failure probability $\delta$, reach $\tau$ and
exponent $r_1$ are here $d$, $m$, $N$, $y$, $q$, $\tau_{\mathcal M}$ and
$r_k$.

\subsection{The tube model and the clean points}\label{subsec:tube-model}

Let $\mathcal M\subset\mathbb R^d$ satisfy the hypotheses of
Section~\ref{sec:setup}: compact, connected, orientable, smooth, of dimension
$m$, without boundary. Its \emph{reach} $\tau_{\mathcal M}$ is the largest
$\tau$ such that every point of $\mathbb R^d$ at distance less than $\tau$ from
$\mathcal M$ has a unique nearest point on $\mathcal M$
\cite[Definition~4.1]{federer}; see also \cite[Definition~1]{aizenbud}. For a
compact smooth submanifold it is positive
\cite[p.~432 and Theorem~4.12]{federer}. For $0<\sigma<\tau_{\mathcal M}$ let
\[
  \mathcal M_\sigma:=\bigl\{y\in\mathbb R^d:\operatorname{dist}(y,\mathcal M)<\sigma\bigr\}
\]
be the open tube of radius $\sigma$ around $\mathcal M$, and let $\pi$ denote
the nearest-point projection onto $\mathcal M$, which is defined on
$\mathcal M_\sigma$ because $\sigma<\tau_{\mathcal M}$. In the model of
\cite[Section~2.1]{aizenbud}, the observations are uniform in
$\mathcal M_\sigma$, and we take their nearest points on $\mathcal M$ as the
clean points:
\begin{equation}\label{eq:tube-model}
  y_1,\ldots,y_n\ \text{ sampled i.i.d.\ from }\operatorname{Unif}(\mathcal M_\sigma),
  \qquad
  x_i:=\pi(y_i),
  \qquad
  \delta_i:=y_i-x_i,
\end{equation}
where $\operatorname{Unif}(\mathcal M_\sigma)$ is Lebesgue measure on the tube
normalized to a probability. Only the $y_i$ are observed. Throughout this
appendix $\mathcal M$ and $\sigma$ are fixed, and only the numbers of
observations grow.

For $x\in\mathcal M$, let $N_x\mathcal M$ be the normal space of $\mathcal M$
at $x$. For $v\in N_x\mathcal M$, the shape operator $S_v(x)$ is the
symmetric linear map of the tangent space $T_x\mathcal M$ with
$\langle S_v(x)u,w\rangle=\langle\mathrm{II}_x(u,w),v\rangle$, where
$\mathrm{II}$ is the second fundamental form. It measures how $\mathcal M$
curves in the direction $v$, and it enters the density of the clean points
through the change of variables $y=x+v$.

\begin{lemma}\label{lem:tube-density}
Let $0<\sigma<\tau_{\mathcal M}$ and let $y$ have distribution
$\operatorname{Unif}(\mathcal M_\sigma)$. Then $x=\pi(y)$ has distribution
$\mu_\sigma=\rho_\sigma\,d\operatorname{Vol}_{\mathcal M}$, with density
\begin{equation}\label{eq:tube-density}
  \rho_\sigma(x)
  =\frac{1}{|\mathcal M_\sigma|}
    \int_{\{v\in N_x\mathcal M:\ \|v\|<\sigma\}}
      \det\bigl(I-S_v(x)\bigr)\,dv,
\end{equation}
where $|\mathcal M_\sigma|$ is the Lebesgue measure of the tube and $dv$ is
Lebesgue measure on $N_x\mathcal M$. The density $\rho_\sigma$ is smooth and
strictly positive on $\mathcal M$, hence bounded above and below by positive
constants.
\end{lemma}

\begin{proof}
Since $\sigma<\tau_{\mathcal M}$, the map $\Psi(x,v)=x+v$, defined for
$x\in\mathcal M$ and $v\in N_x\mathcal M$ with $\|v\|<\sigma$, parametrizes
$\mathcal M_\sigma$, and $\pi(\Psi(x,v))=x$
\cite[Theorem~4.8(12)]{federer}. By the Weingarten formula
\cite[(3.3)]{chen_submanifolds}, the derivative of a normal vector field with
value $v$ at $x$ in a tangent direction $u$ has tangential part $-S_v(x)u$.
Hence the Jacobian determinant of $\Psi$ with respect to
$d\operatorname{Vol}_{\mathcal M}(x)\,dv$ is $\det(I-S_v(x))$, as in the first
step of the proof of Weyl's tube formula
\cite[Section~2]{burtscher_heckman}. It follows from
\cite[Theorem~4.8(7)]{federer} that
$\|S_v(x)\|_{\mathrm{op}}\leq\|v\|/\tau_{\mathcal M}$, so
\[
  \det(I-S_v(x))\geq(1-\sigma/\tau_{\mathcal M})^m>0.
\]
The change of variables $y=\Psi(x,v)$ gives \eqref{eq:tube-density}. In a
local orthonormal frame of the normal bundle that depends smoothly on $x$,
\eqref{eq:tube-density} integrates a function that is smooth in $(x,v)$ over
a fixed ball. Hence $\rho_\sigma$ is smooth, and it is positive because the
integrand is. Since $\mathcal M$ is compact, $\rho_\sigma$ is bounded above
and below by positive constants.
\end{proof}

\begin{lemma}\label{lem:tube-instance}
Under \eqref{eq:tube-model}, the clean points $x_1,\ldots,x_n$ are i.i.d.\
with distribution $\mu_\sigma$, and $\|\delta_i\|<\sigma$ for every $i$.
Thus the hypotheses of Section~\ref{sec:setup} on the manifold, the sampling
density and the displacements hold with $\rho=\rho_\sigma$ and
$\delta=\sigma$.
\end{lemma}

\begin{proof}
Each clean point $x_i=\pi(y_i)$ is the image of $y_i$ alone under the same
measurable map $\pi$. Since the observations are i.i.d., so are the clean
points, with common distribution $\mu_\sigma$ by Lemma~\ref{lem:tube-density}. That
lemma also gives the smoothness of $\rho_\sigma$ and its bounds. Finally,
$\|\delta_i\|=\operatorname{dist}(y_i,\mathcal M)<\sigma$ because
$y_i\in\mathcal M_\sigma$.
\end{proof}

The continuum operator is $\Delta_{\rho_\sigma}$, defined as in
Section~\ref{sec:setup} with $\rho=\rho_\sigma$. Although the observations
are uniform in the tube, $\rho_\sigma$ need not be constant when $m\geq2$.
For a surface in $\mathbb R^3$, for example, \eqref{eq:tube-density} gives
$\rho_\sigma(x)=\bigl(2\sigma+\tfrac23\sigma^3K(x)\bigr)/|\mathcal M_\sigma|$,
where $K$ is the Gauss curvature of $\mathcal M$. For $m=1$, $\rho_\sigma$ is
constant: the tangent space is one-dimensional, so $\det(I-S_v(x))-1$ is
linear in $v$ and integrates to zero over the ball.

\subsection{The projection algorithm and the residual}\label{subsec:as-algorithm}

For $k\geq2$, $q\in(0,1)$, a query point $y\in\mathcal M_\sigma$ and a
dataset $R$ of $N$ points of $\mathcal M_\sigma$, let $\mathcal A_{k,q}(y;R)$
denote the output of Algorithms~1 and~2 of \cite{aizenbud}, in the form
analysed in their theorems. It is an estimate of $\pi(y)$, obtained by
fitting polynomials of degree at most $k-1$ to the points of $R$ near $y$,
over an estimate of the tangent space of $\mathcal M$ at $\pi(y)$ that is
refined at each iteration. As in their analysis, the fits use the robust
regression and derivative estimates of \cite[(12)]{aizenbud}, and the number
of iterations is that of \cite[Lemma~27]{aizenbud}, which depends on $N$ and
$q$. The algorithm also needs $m$, $\sigma$ and a suitable lower bound on
the reach.

We take $\mathcal A_{k,q}(y;R)$ to be a single point of $\mathbb R^d$ for
every $y$ and $R$, and we assume that it is a measurable function of the
pair $(y,R)$. The proof of Proposition~\ref{prop:as-residual} uses this to
integrate over a random query.

\begin{theorem}[{\cite[Theorems~1 and~3]{aizenbud}}]\label{thm:as-quoted}
Assume that $\mathcal M\in C^k$ is a compact $m$-dimensional submanifold of
$\mathbb R^d$ without boundary, that
$\tau_{\mathcal M}/\sigma>C_\tau\sqrt{d\log d}$ for the constant $C_\tau$ of
the cited theorem, and that $R$ consists of $N$ data points sampled i.i.d.\
from $\operatorname{Unif}(\mathcal M_\sigma)$. Let $y\in\mathcal M_\sigma$.
There are a constant $C_2$ and, for each $q\in(0,1)$, a threshold
$N_0(y,q)$ such that for every $N>N_0(y,q)$, with probability at least $1-q$,
\begin{equation}\label{eq:as-quoted}
  \bigl\|\mathcal A_{k,q}(y;R)-\pi(y)\bigr\|
  \;\leq\;C_2\ln(1/q)\left(\frac{N}{(\ln\ln N)^2}\right)^{-r_k},
  \qquad
  r_k:=\frac{k-1}{2k+m}.
\end{equation}
\end{theorem}

Aizenbud and Sober \cite[Theorem~1]{aizenbud} also bound the distance from
$\mathcal A_{k,q}(y;R)$ to $\mathcal M$, with the larger exponent
$k/(2k+m)$. We use the bound \eqref{eq:as-quoted} on the distance to
$\pi(y)$, because the results of Section~\ref{sec:results} need each
denoised point to be close to its own clean point. Since $\mathcal M$ is smooth, it
is $C^k$ for every $k$, so $k\geq2$ is a free parameter of the algorithm.
Enlarging $N_0(y,q)$ if necessary, we assume that $N_0(y,q)\geq3$, so that
$\ln\ln N>0$ in \eqref{eq:as-quoted}.

In Theorem~\ref{thm:as-quoted} the query point $y$ is fixed, and the
probability is over the dataset $R$. To apply the theorem at the
observations $y_i$, which are random, we run the algorithm on an additional
dataset, independent of the observations.

\begin{definition}[Denoiser and denoised Laplacian]
\label{def:as-denoiser}
Under \eqref{eq:tube-model}, fix $k\geq2$, $q\in(0,1)$ and an integer
$N\geq1$. Let the \emph{additional dataset} $R_N$ consist of $N$ data
points sampled i.i.d.\ from $\operatorname{Unif}(\mathcal M_\sigma)$,
independently of $y_1,\ldots,y_n$, and set
\begin{equation}\label{eq:as-denoiser}
  \hat x_i:=\mathcal A_{k,q}(y_i;R_N),
  \qquad i=1,\ldots,n.
\end{equation}
All $n$ queries use the same additional dataset. The \emph{residual}
\begin{equation}\label{eq:residual}
  \widehat\delta:=\max_{1\leq i\leq n}\|\hat x_i-x_i\|
\end{equation}
is the largest distance from a denoised point to its clean point. Let
$\widehat L$ be the graph Laplacian of Section~\ref{sec:setup} built on
$\hat x_1,\ldots,\hat x_n$, let $\hat\lambda_1\leq\cdots\leq\hat\lambda_n$ be
its eigenvalues, and define $\widehat P_l$ from $\widehat L$ as $P_l$ is
defined from $L$ in Section~\ref{sec:results}.
\end{definition}

\begin{proposition}\label{prop:as-residual}
Assume \eqref{eq:tube-model} with
$\tau_{\mathcal M}/\sigma>C_\tau\sqrt{d\log d}$, and let $k\geq2$ and
$p\geq1$. For $n\geq2$ set $q_n:=n^{-(p+1)}$, and run the denoiser of
Definition~\ref{def:as-denoiser} with $q=q_n$. For every $n\geq2$ and $t>0$
there is an integer $N_1(n,t)$ such that, for every size
$N\geq N_1(n,t)$ of the additional dataset,
\[
  \mathbb P\bigl(\widehat\delta>t\bigr)\;\leq\;2n^{-p},
\]
where the probability is over the observations $y_1,\ldots,y_n$ and the
additional dataset $R_N$.
\end{proposition}

\begin{proof}
Fix $n\geq2$ and $t>0$, and write $\nu:=\operatorname{Unif}(\mathcal M_\sigma)$,
so that each $y_i$ has distribution $\nu$ and $R_N$ has distribution
$\nu^{\otimes N}$, the product measure on $\mathcal M_\sigma^N$. We bound the
probability that one distance $\|\hat x_i-x_i\|$ exceeds $t$, first with the
random query $y_i$ replaced by a fixed query, using
Theorem~\ref{thm:as-quoted}, and then by integrating that bound against the
distribution of $y_i$. A union bound over $i$ then multiplies the bound by
$n$, and $nq_n=n^{-p}$.

For $N\geq1$, $y\in\mathcal M_\sigma$ and $R\in\mathcal M_\sigma^N$, let
\[
  g_N(y,R):=\bigl\|\mathcal A_{k,q_n}(y;R)-\pi(y)\bigr\|
\]
be the error of the algorithm at the query $y$ with dataset $R$,
and let
\[
  h_N(y):=\nu^{\otimes N}\bigl(\{R:\ g_N(y,R)>t\}\bigr)
\]
be the probability that this error exceeds $t$ when the query $y$ is fixed.
The function $g_N$ is measurable on $\mathcal M_\sigma\times\mathcal M_\sigma^N$,
because $\mathcal A_{k,q_n}$ is jointly measurable and $\pi$ is continuous
\cite[Theorem~4.8(4)]{federer}.
By Tonelli's theorem, $h_N$ is therefore measurable on $\mathcal M_\sigma$,
and $0\leq h_N\leq1$.

First fix a query $y\in\mathcal M_\sigma$. Since $\ln(1/q_n)=(p+1)\ln n$,
Theorem~\ref{thm:as-quoted} with $q=q_n$ gives, for every $N>N_0(y,q_n)$,
\[
  g_N(y,R)\;\leq\;C_2\,(p+1)\,\ln n\left(\frac{N}{(\ln\ln N)^2}\right)^{-r_k}
\]
for all $R$ outside a set of $\nu^{\otimes N}$-measure at most $q_n$. Since
$r_k>0$, the right-hand side tends to zero as $N\to\infty$, so it is at most
$t$ for all large $N$. Hence $h_N(y)\leq q_n$ for all sufficiently large
$N$.

Next we integrate over the query. By the previous paragraph,
$(h_N-q_n)_+\to0$ pointwise on $\mathcal M_\sigma$ as $N\to\infty$. Since
$0\leq(h_N-q_n)_+\leq1$ and $\nu$ is a probability measure, dominated
convergence gives $\int(h_N-q_n)_+\,d\nu\to0$. Hence there is an integer $N_1(n,t)$ such that, for every
$N\geq N_1(n,t)$,
\[
  \int_{\mathcal M_\sigma}h_N\,d\nu
  \;\leq\;q_n+\int_{\mathcal M_\sigma}(h_N-q_n)_+\,d\nu
  \;\leq\;2q_n.
\]

Now let $N\geq N_1(n,t)$ and fix $i$. By \eqref{eq:tube-model} and
\eqref{eq:as-denoiser}, $\|\hat x_i-x_i\|=g_N(y_i,R_N)$. Since $R_N$ is
independent of $y_i$, the pair $(y_i,R_N)$ has the product distribution
$\nu\otimes\nu^{\otimes N}$. Tonelli's theorem, applied to the indicator of
$\{(y,R):\ g_N(y,R)>t\}$, therefore gives
\[
\begin{aligned}
  \mathbb P\bigl(\|\hat x_i-x_i\|>t\bigr)
  &=\int_{\mathcal M_\sigma}\int_{\mathcal M_\sigma^N}
     \mathbf 1_{\{g_N(y,R)>t\}}\,d\nu^{\otimes N}(R)\,d\nu(y)\\
  &=\int_{\mathcal M_\sigma}h_N\,d\nu
  \;\leq\;2q_n.
\end{aligned}
\]
A union bound over $i$ gives
\[
  \mathbb P\bigl(\widehat\delta>t\bigr)
  \leq\sum_{i=1}^n\mathbb P\bigl(\|\hat x_i-x_i\|>t\bigr)
  \leq 2nq_n=2n^{-p}. \qedhere
\]
\end{proof}

\subsection{The denoised spectrum}\label{subsec:denoised-spectrum}

The displacements $\hat x_i-x_i$ of the denoised points are random, and they
depend on the observations and on the additional dataset. The results of
Section~\ref{sec:results} still apply to them, because we assume only that
the displacements are bounded by a deterministic $\delta$; the displacements
themselves may be random. Lemma~\ref{lem:denoise-uniform} makes this precise.

\begin{lemma}\label{lem:denoise-uniform}
Let $\delta\geq0$ be deterministic, and assume the hypotheses of one of the
results of Section~\ref{sec:results}. The proof of that result gives its bound
given an event that depends only on the clean points $x_1,\ldots,x_n$, namely
$\mathcal E_n$ for Theorem~\ref{thm:main} and
Corollary~\ref{cor:continuum-eig}, $\mathcal E_n^{\mathrm{proj}}$ for
Theorem~\ref{thm:projector-headline}, $\mathcal H_{n,l}$ for
Theorem~\ref{thm:eigenvector-headline} and $\mathcal E_n^{(2)}$ for
Remark~\ref{rem:index-one}. Given that event, the bound holds, with the same
constants, for the graph Laplacian of
every cloud $z_1,\ldots,z_n\in\mathbb R^d$ with $\|z_i-x_i\|\leq\delta$ for
all $i$, and for every choice of its spectral projectors and unit
eigenvectors.
\end{lemma}

\begin{proof}
In the proofs of the results of Section~\ref{sec:results}, we use the noisy
points only through Lemma~\ref{lem:weights}, which holds for every cloud as in
the statement. Everything else we use, namely $L^\pm$, the comparison operators
$A,B$, $P_l^+$, $\tau_l$, $\Gamma_l^+$ and the events, is built from
$x_1,\ldots,x_n$, $\varepsilon$ and $\delta$ alone, so the events depend only
on the clean points. Lemma~\ref{lem:fcpb} holds for every symmetric $M$ with
$A\preceq M\preceq B$, and Corollary~\ref{cor:cluster-vector} for every
choice of eigenbases. So each proof applies without change to the graph
Laplacian of $z_1,\ldots,z_n$.
\end{proof}

It is therefore enough that, with high probability, the residual is at most
a deterministic bound. We take this bound to be $\varepsilon_n^3$, where
$\varepsilon_n$ is the bandwidth, since this is the size of displacement at
which Theorem~\ref{thm:eigenvector-headline} gives the eigenvector rate
$O(\varepsilon_n)$ of \cite[Theorem~2.7(2)]{calder_gt}.
Proposition~\ref{prop:as-residual}
makes the residual at most $\varepsilon_n^3$ with probability at least
$1-2n^{-p}$ once the additional dataset is large enough.

\begin{theorem}
\label{thm:denoised-spectrum}
Assume the hypotheses of Proposition~\ref{prop:as-residual}, and run the
denoiser of Definition~\ref{def:as-denoiser} with $q=q_n$, as there. Let
$\eta$ satisfy the hypotheses of Section~\ref{sec:setup}, and fix $l\geq1$.
Use the notation and constants of Sections~\ref{sec:setup}
and~\ref{sec:results} for the density $\rho_\sigma$, and build $\widehat L$ at
$\varepsilon_n=C_\varepsilon(\log n/n)^{1/(m+4)}$, where $C_\varepsilon$ is
large enough depending on $p,l,\mathcal M,\rho_\sigma$ and $\eta$, as in
Section~\ref{sec:results}. For each $n\geq2$, let the size $N(n)$ of the additional
dataset be deterministic, with
\begin{equation}\label{eq:reference-size}
  N(n)\;\geq\;N_1\bigl(n,\varepsilon_n^3\bigr)
\end{equation}
for the size $N_1$ of Proposition~\ref{prop:as-residual}. Then, for $n$
large enough, each of the following holds with probability at least
$1-O(n^{-p})$, over the observations and the additional dataset.
\begin{enumerate}
  \item The $l$-th eigenvalues satisfy
        \begin{equation}\label{eq:denoised-eigenvalue}
          \bigl|\lambda_l-\hat\lambda_l\bigr|
            \;\leq\;K_l\,\varepsilon_n^2,
          \qquad
          \bigl|\hat\lambda_l-\sigma_\eta\lambda_l^{\mathcal M}\bigr|
            \;\leq\;K_l\,\varepsilon_n^2+C_l'\,\varepsilon_n.
        \end{equation}
  \item If $\gamma_l^{\mathcal M}>0$, then for every choice of $P_l$ and
        $\widehat P_l$,
        \begin{equation}\label{eq:denoised-projector}
          \bigl\|P_l-\widehat P_l\bigr\|_{\mathrm{HS},n}
          \;\leq\;D_l\,\varepsilon_n.
        \end{equation}
  \item If $l\geq2$, then every eigenvector $\hat u_l$ of $\widehat L$ with
        eigenvalue $\hat\lambda_l$ and $\|\hat u_l\|_n=1$ admits an
        eigenfunction $f_l$ of $\Delta_{\rho_\sigma}$ with eigenvalue
        $\lambda_l^{\mathcal M}$ and $\|f_l\|_{L^2(\mu_\sigma)}=1$ such that
        \begin{equation}\label{eq:denoised-eigenvector}
          \bigl\|\hat u_l-f_l|_{X_n}\bigr\|_n
          \;=\;O(\varepsilon_n).
        \end{equation}
\end{enumerate}
Here $\lambda_l$ and $P_l$ belong to the Laplacian $L$ built on the clean
points $x_i=\pi(y_i)$, and $f_l|_{X_n}$ is the vector with entries
$f_l(x_i)$. The clean points are not observed.
\end{theorem}

\begin{proof}
By Lemma~\ref{lem:tube-instance}, the clean points satisfy the hypotheses of
Section~\ref{sec:setup} with $\rho=\rho_\sigma$. We apply the results of
Section~\ref{sec:results} with the deterministic displacement bound
$\delta_n:=\varepsilon_n^3$. Since
$\delta_n/\varepsilon_n=\varepsilon_n^2\to0$, the condition of
Section~\ref{sec:results} that $\delta/\varepsilon$ be sufficiently small
holds at $\delta=\delta_n$ for $n$ large enough. For each item, let
$\mathcal G_n$ be the event that the proofs of the corresponding results use
at $\delta=\delta_n$, namely $\mathcal E_n$,
$\mathcal E_n^{\mathrm{proj}}$ and $\mathcal H_{n,l}$ of
Sections~\ref{sec:eigenvalue} and~\ref{sec:projector} for items 1, 2 and~3.

First we bound the probabilities. The event $\mathcal G_n$ depends only on
the clean points (Lemma~\ref{lem:denoise-uniform}), which are i.i.d.\ with
distribution $\mu_\sigma$ (Lemma~\ref{lem:tube-instance}). The probability
bounds proved in those sections therefore apply and give
$\mathbb P(\mathcal G_n)\geq1-O(n^{-p})$. By \eqref{eq:reference-size},
Proposition~\ref{prop:as-residual} with $t=\delta_n$ gives
$\mathbb P(\widehat\delta>\delta_n)\leq2n^{-p}$. A union bound then gives
probability at least $1-O(n^{-p})$ for
$\mathcal G_n\cap\{\widehat\delta\leq\delta_n\}$.

Given $\{\widehat\delta\leq\delta_n\}$, every denoised point satisfies
$\|\hat x_i-x_i\|\leq\delta_n$. By Lemma~\ref{lem:denoise-uniform}, the
bounds of Theorem~\ref{thm:main}, Corollary~\ref{cor:continuum-eig},
Theorem~\ref{thm:projector-headline} and
Theorem~\ref{thm:eigenvector-headline} at $\delta=\delta_n$ therefore hold
for $\widehat L$ given $\mathcal G_n\cap\{\widehat\delta\leq\delta_n\}$.
These bounds are $K_l\delta_n/\varepsilon_n$,
$K_l\delta_n/\varepsilon_n+C_l'\varepsilon_n$,
$D_l\sqrt{\delta_n/\varepsilon_n}$ and
$O(\sqrt{\delta_n/\varepsilon_n}+\varepsilon_n)$, and substituting
$\delta_n/\varepsilon_n=\varepsilon_n^2$ gives the asserted bounds.
\end{proof}

Item 3 excludes $l=1$, as Theorem~\ref{thm:eigenvector-headline} does. For
$l=1$, the same argument with the event $\mathcal E_n^{(2)}$ used for
Remark~\ref{rem:index-one} shows that, for $n$ large enough and with
probability at least $1-O(n^{-p})$, every eigenvector $\hat u_1$ of
$\widehat L$ with eigenvalue $\hat\lambda_1$ and $\|\hat u_1\|_n=1$ is
constant, so that $\|\hat u_1-f_1|_{X_n}\|_n=0$ for the matching constant
eigenfunction $f_1$.
\section{Hypotheses of the Spectral Convergence Results of Calder and Garc\'ia Trillos}
\label{app:cgt-hypotheses}

The constant $C_l'$ and the event $\mathcal E_{n,3}$ of
Section~\ref{subsec:hpe} come from \cite[Theorem~2.4(2)]{calder_gt}, and
Section~\ref{subsec:proof-thm-eigenvector} also uses the proof of
\cite[Theorem~2.7(2)]{calder_gt}. We state the hypotheses of these results
in our notation. They are of three kinds.

\emph{The manifold, the density and the kernel.}
\cite[Section~2]{calder_gt} assumes that $\mathcal M$ is a compact,
connected, orientable, smooth $m$-dimensional submanifold of $\mathbb R^d$,
that $\rho\in C^{2,\vartheta}(\mathcal M)$ with
$0<\rho_{\min}\leq\rho\leq\rho_{\max}$, and that $\eta$ is non-increasing,
supported in $[0,1]$, Lipschitz on $[0,1]$ and positive at $\tfrac12$.
Section~\ref{sec:setup} assumes all of this.

\emph{Deterministic restrictions on the parameters.} Besides $\varepsilon$,
the results of \cite{calder_gt} involve two parameters $\theta,\delta_e>0$
that do not enter the definition of $L$. The first controls the error in
approximating the density, the second the transport distance from that
approximation to the empirical measure; the event on which both are
controlled is described below. Here and below, $C,c>0$ denote constants of
\cite{calder_gt}; they do not depend on $n,\varepsilon,\theta$ or
$\delta_e$. The restrictions of \cite[Assumption~2.3]{calder_gt} are:
\begin{enumerate}
  \item[(A1)] $\varepsilon$ is small enough, and in particular
        $2\varepsilon<\min\{1,i_0,K^{-1/2},R/2\}$, where $i_0>0$ is a lower
        bound on the injectivity radius of $\mathcal M$, $K>0$ an upper
        bound on the absolute values of its sectional curvatures, and $R$
        its reach;
  \item[(A2)] $\delta_e\leq\varepsilon/4$;
  \item[(A3)] $\delta_e>n^{-1/m}$;
  \item[(A4)] $C(\theta+\delta_e)\leq\rho_{\min}/2$, where $C$ depends only
        on $m$, $\rho$ and $\mathcal M$.
\end{enumerate}
Each result also has a smallness condition at the eigenvalue indices it
uses. For an index $j$, let
$g_j:=\min\{|\lambda_k^{\mathcal M}-\lambda_j^{\mathcal M}|:
\lambda_k^{\mathcal M}\neq\lambda_j^{\mathcal M}\}>0$ be the distance from
$\lambda_j^{\mathcal M}$ to the other eigenvalues of $\Delta_\rho$, and let
$c_j>0$ be the smaller of the thresholds that
\cite[Theorem~2.4(1) and~(2)]{calder_gt} give at index $j$.
\begin{enumerate}
  \item[(S1$_j$)] \cite[Theorem~2.4(2)]{calder_gt} at index $j$ requires
        $\bigl(\varepsilon\sqrt{\lambda_j^{\mathcal M}}+\delta_e/\varepsilon
        +\varepsilon+\theta\bigr)\lambda_j^{\mathcal M}\leq c_j$.
  \item[(S2$_N$)] The proof of \cite[Theorem~2.7(2)]{calder_gt} fixes
        $N\in\mathbb N$ and requires, before its equation (4.13), that for
        every $j\leq N$ the hypothesis
        $\varepsilon\sqrt{\lambda_j^{\mathcal M}}+C(\theta+\delta_e)\leq c_j$
        of \cite[Theorem~2.4(1)]{calder_gt} holds and that the error bound
        of that theorem satisfies
        $C\bigl(\varepsilon\sqrt{\lambda_j^{\mathcal M}}
        +\delta_e/\varepsilon+\varepsilon+\theta\bigr)\lambda_j^{\mathcal M}
        \leq\tfrac12\sigma_\eta g_j$.
\end{enumerate}
The proof in \cite{calder_gt} writes $\tfrac12g_j$ in place of
$\tfrac12\sigma_\eta g_j$ and omits $\sigma_\eta$ in the displays after
(4.13). Restoring the factor changes only the constant $C$, and the argument
needs it, since (4.13) compares $\lambda_j$ with
$\sigma_\eta\lambda_j^{\mathcal M}$.

\emph{Probabilistic events.} The estimates of \cite{calder_gt} are proved
on the event of \cite[Proposition~2.12]{calder_gt}: under (A1)--(A4), with
probability at least $1-n\exp(-Cn\theta^2\delta_e^m)$ there are a
probability density $\widetilde\rho_n$ on $\mathcal M$ with
$\|\rho-\widetilde\rho_n\|_{L^\infty}\leq C(\theta+\delta_e)$ and a map
$T:\mathcal M\to X_n$ that pushes
$\widetilde\rho_n\,d\operatorname{Vol}_{\mathcal M}$ forward to the
empirical measure of $X_n$, with $d_{\mathcal M}(x,T(x))\leq\delta_e$ for
every $x$. The event $\mathcal E_{n,3}$, and the event of
Lemma~\ref{lem:cgt-cluster-correspondence}, are defined without $\theta$
and $\delta_e$. These parameters enter only the conditions under which
\cite{calder_gt} bounds the failure probabilities. Under (A1)--(A4) and
(S1$_l$), \cite[Theorem~2.4(2)]{calder_gt} gives
\[
  \mathbb P_n\bigl(\mathcal E_{n,3}^c\bigr)
  \;\leq\;2Cn\exp\bigl(-cn\varepsilon^{m+4}\bigr)
        +Cn\exp\bigl(-cn\theta^2\delta_e^m\bigr),
\]
and Lemma~\ref{lem:cgt-cluster-correspondence} gives a bound of the same
form under (A1)--(A4) and (S2$_N$). Lemma~\ref{lem:cgt-parameters} chooses
$\theta$ and $\delta_e$ so that these conditions hold for finitely many
indices at once and both terms are $O(n^{-p})$: the term with
$\exp(-cn\theta^2\delta_e^m)$ for every $C_\varepsilon>0$ once $n$ is large
enough, and the term with $\exp(-cn\varepsilon^{m+4})$ once $C_\varepsilon$
is large enough.

\begin{lemma}\label{lem:cgt-parameters}
Fix $N\geq1$ and $p\geq1$, and let
$\varepsilon=C_\varepsilon(\log n/n)^{1/(m+4)}$ with $C_\varepsilon>0$. There
is a $\theta_N\in(0,\tfrac14]$, depending only on $N,\mathcal M,\rho$ and
$\eta$, such that the choice
\[
  \theta:=\theta_N,
  \qquad
  \delta_e:=\theta_N\,\varepsilon
\]
satisfies (A1)--(A4), (S1$_j$) for every $j\leq N$, and (S2$_N$), for all
$n$ large enough. Moreover, for every constant $c>0$,
$n\exp(-cn\theta^2\delta_e^m)\leq n^{-p}$ for all $n$ large enough, and
$n\exp(-cn\varepsilon^{m+4})\leq n^{-p}$ for all $n$ if
$cC_\varepsilon^{m+4}\geq p+1$.
\end{lemma}

\begin{proof}
Each inequality in (A4), (S1$_j$) and (S2$_N$), for $j\leq N$, bounds a
linear combination of $\varepsilon$, $\theta$, $\delta_e/\varepsilon$ and
$\delta_e$ with nonnegative coefficients by a positive number, such as
$c_j$ or $\tfrac12\sigma_\eta g_j$. The coefficients and these numbers
depend only on $j$, $\mathcal M$, $\rho$ and $\eta$, and there are
finitely many such inequalities. Choose $\theta_N\in(0,\tfrac14]$ so small
that in each of them, with $\theta=\delta_e/\varepsilon=\theta_N$, the
terms in $\theta$ and $\delta_e/\varepsilon$ add up to at most half of the
right-hand side. The remaining terms are multiples of $\varepsilon$ and
$\delta_e=\theta_N\varepsilon$, which tend to zero as $n\to\infty$, so each
inequality holds for $n$ large. So does (A1), since $\varepsilon\to0$.
(A2) holds since $\theta_N\leq\tfrac14$, and (A3) holds for $n$ large since
\[
  n^{1/m}\delta_e
  =\theta_NC_\varepsilon(\log n)^{1/(m+4)}\,n^{1/m-1/(m+4)}
  \;\to\;\infty.
\]
For the first probability bound,
$n\theta^2\delta_e^m=\theta_N^{m+2}C_\varepsilon^m\,n^{4/(m+4)}(\log n)^{m/(m+4)}$,
which is at least $c^{-1}(p+1)\log n$ for $n$ large, and then
$n\exp(-cn\theta^2\delta_e^m)\leq n\cdot n^{-(p+1)}=n^{-p}$. For the second,
$n\varepsilon^{m+4}=C_\varepsilon^{m+4}\log n$, so
$n\exp(-cn\varepsilon^{m+4})=n^{1-cC_\varepsilon^{m+4}}\leq n^{-p}$.
\end{proof}

By Lemma~\ref{lem:cgt-parameters} with $N=l$, the hypotheses of
\cite[Theorem~2.4(2)]{calder_gt} at index $l$ hold for $n$ large, and both
terms in the bound on $\mathbb P_n(\mathcal E_{n,3}^c)$ are at most
$2Cn^{-p}$ for $n$ large once $cC_\varepsilon^{m+4}\geq p+1$. Hence, for
$C_\varepsilon$ large enough,
$\mathbb P_n(\mathcal E_{n,3})\geq1-O(n^{-p})$, which is the probability
bound in \eqref{eq:cgt-cited}. At several indices, each of these events has
this probability bound, and
Lemma~\ref{lem:cgt-parameters}, with $N$ the largest index, gives one
choice of $\theta$ and $\delta_e$ that serves all of them.

\end{document}